\documentclass[11pt,a4paper]{article}
\usepackage[T1]{fontenc}
\usepackage{lmodern}
\usepackage{authblk}

\usepackage[margin=25mm]{geometry}
\usepackage{amsmath,amssymb,amsthm,mathtools,microtype,xcolor}
\usepackage{needspace}
\usepackage{graphicx}
\makeatletter
\let\inputrows\@@input
\makeatother
\usepackage{float}
\AddToHook{env/table/begin}{%
  \setlength{\abovecaptionskip}{10pt}%
  \setlength{\belowcaptionskip}{10pt}%
}
\usepackage{placeins}
\usepackage[authoryear,round,semicolon,sort]{natbib}
\setcitestyle{aysep={,},yysep={,},notesep={, }}
\usepackage[unicode,colorlinks=true,linkcolor=blue!45!black,
  citecolor=blue,urlcolor=blue!45!black]{hyperref}
\usepackage{doi}

\hypersetup{
  pdftitle={Policy and Value Iteration for Entropy-Regularized Stochastic Control},
  pdfauthor={Yuling Jiao, Peiying Wu, Jerry Zhijian Yang, Ruizhe Zhang, Tao Zhou},
  pdfkeywords={soft HJB equations, approximate policy iteration, entropy-regularized control, Gibbs policy improvement, error propagation}
}
\allowdisplaybreaks[2]
\numberwithin{equation}{section}
\theoremstyle{definition}
\newtheorem{assumption}{Assumption}[section]
\newtheorem{definition}{Definition}[section]
\theoremstyle{plain}
\newtheorem{proposition}{Proposition}[section]
\newtheorem{lemma}{Lemma}[section]
\newtheorem{theorem}{Theorem}[section]
\newtheorem{corollary}{Corollary}[section]
\newcommand{\R}{\mathbb R}
\newcommand{\E}{\mathbb E}
\newcommand{\PP}{\mathcal P}
\newcommand{\LL}{\mathcal L}
\newcommand{\HH}{\mathcal H}
\newcommand{\KL}{D_{\mathrm{KL}}}
\newcommand{\dd}{\,\mathrm d}
\newcommand{\tr}{\operatorname{tr}}

\title{Policy and Value Iteration\\for Entropy-Regularized Stochastic Control}
\author[1]{Yuling Jiao}
\author[2]{Peiying Wu}
\author[2]{Jerry Zhijian Yang}
\author[2]{Ruizhe Zhang}
\author[3]{Tao Zhou}

\affil[1]{School of Artificial Intelligence, Wuhan University,\protect\\
Wuhan, Hubei 430072, China.\protect\\
\texttt{yulingjiaomath@whu.edu.cn}}

\affil[2]{School of Mathematics and Statistics, Wuhan University,\protect\\
Wuhan, Hubei 430072, China.\protect\\
\texttt{\{peiyingwu,zjyang.math,ruizhezhang\}@whu.edu.cn}}

\affil[3]{SKLMS \& Institute of Computational Mathematics and\protect\\
Scientific/Engineering Computing, Academy of Mathematics and Systems Science,\protect\\
Chinese Academy of Sciences, Beijing 100190, China.\protect\\
\texttt{tzhou@lsec.cc.ac.cn}}
\date{November 25, 2026}

\begin{document}
\maketitle
\vspace{-0.8em}

\begin{abstract}
We study finite-horizon entropy-regularized stochastic control with
action-independent diffusion. To solve the associated
Hamilton--Jacobi--Bellman (HJB) equation, we alternate between value and
policy updates. Each iteration approximately solves a linear backward
equation for the current policy, then approximates the Gibbs distribution
defined by the resulting value gradient. We prove a uniform bound on the
final policy's excess regularized cost. The bound assigns explicit,
summable weights to the evaluation residuals and the Kullback--Leibler
divergences from the updated policies to their Gibbs targets. Under the
stated assumptions and with exact terminal data, vanishing errors in both
steps imply uniform convergence to the optimal value.
Two finite-horizon execution problems test value accuracy against a known
solution and closed-loop performance when no exact solution is available.
We also evaluate the iteration on a discounted infinite-horizon routing
problem with twelve modes and an unregularized HJB benchmark in up to
10,000 dimensions.
\end{abstract}


\section{Introduction}\label{sec:introduction}

Finite-horizon stochastic control concerns sequential decisions under
uncertainty over a prescribed time interval. In pumped-storage scheduling,
reservoir levels and electricity prices guide decisions about when to pump
water and when to generate power \citep{PV2021PumpedHydro}. In multi-asset
execution, traders choose rates jointly over a fixed trading window based
on market conditions, remaining positions, and cross-asset price impact
\citep{BL1998Execution,CJ2016OrderFlow,BDG2022MultiAsset}. In adaptive
radiotherapy, doses are adjusted to the observed anatomy over a prescribed
number of treatment sessions \citep{RCBT2012AdaptiveFractionation}. In
multisatellite rendezvous, the aim is to reach target regions by a specified
time with sufficiently high probability \citep{PO2024Rendezvous}. In
stochastic thermodynamics, a control protocol can steer the distribution of
a Brownian particle to a target distribution in a fixed time while minimizing
dissipation \citep{SS2007FiniteTimeThermodynamics}. Classical stochastic control minimizes a prescribed cost, but this cost
may not fully capture practical considerations \citep{knox2023reward}.
Expert demonstrations can convey preferences and skills that are difficult
to specify explicitly, as illustrated in autonomous driving and robot
learning \citep{kuderer2015learning,argall2009survey}.
We therefore incorporate prior information from reliable historical actions
or tested baseline rules into the control objective alongside the modeled cost.

We study feedback control for finite-horizon stochastic systems. An
ordinary feedback policy selects an action based on the current time
and state. When two
actions exchange their relative merits as the state varies, the optimal
action can switch abruptly \citep{RZ2021Relaxed}. We adopt relaxed
feedback policies, which assign an action distribution
$\pi(\dd u\mid s,x)$ to each time and state and can respond to changes in
the state by adjusting the probabilities assigned to actions. Starting
from $(t,x)$, the controlled state evolves over $[t,T]$ according to
\begin{equation}\label{eq:relaxed-dynamics}
\dd X_s^\pi
=\left[\int_{\mathcal U}b(s,X_s^\pi,u)\pi(\dd u\mid s,X_s^\pi)\right]\dd s
+\sigma(s,X_s^\pi)\dd W_s,\qquad X_t^\pi=x.
\end{equation}
Here $X_s^\pi\in\R^d$ is the state and
$u\in\mathcal U\subseteq\R^m$ is an action. The drift $b$ is averaged
under the current action distribution. The diffusion coefficient
$\sigma\in\R^{d\times r}$ is independent of the action, and $W$ is an
$r$-dimensional Brownian motion.

Let $\mu_{s,x}$ be a reference action distribution encoding baseline
preferences specified by rules, estimated from reliable historical action
data, or informed by expert knowledge
\citep{TGH2022BehaviorPriors,SSB2020BehaviorPriors,SS2024GradientFlows,RR2014Advice}.
Writing $\pi_s=\pi(\cdot\mid s,X_s^\pi)$, we define the
KL-regularized cost by
\begin{equation}\label{eq:regularized-cost}
J_\tau(t,x;\pi)=\E_{t,x}^\pi\!\left[
g(X_T^\pi)+\int_t^T\left\{
\int_{\mathcal U}c(s,X_s^\pi,u)\pi_s(\dd u)
+\tau\KL(\pi_s\Vert\mu_{s,X_s^\pi})
\right\}\dd s\right],
\end{equation}
where $\tau>0$ balances modeled cost minimization against adherence
to the reference preferences.
The functions $c$ and $g$ are the running and terminal costs, and
$\E^\pi_{t,x}$ averages over trajectories starting from $(t,x)$ under
$\pi$. The reference distribution and the regularization weight
are fixed throughout the iteration. The optimal value over the
admissible control class $\mathcal A_\tau(t,x)$ is
\begin{equation}\label{eq:intro-optimal-value}
v_\tau(t,x)=\inf_{\pi\in\mathcal A_\tau(t,x)}J_\tau(t,x;\pi).
\end{equation}

Under the appropriate regularity and verification conditions, dynamic
programming gives the soft Hamilton--Jacobi--Bellman (HJB) equation for
the optimal value $v_\tau$ \citep{WZZ2020,TZZ2022}
\begin{equation}\label{eq:soft-hjb}
\partial_t v_\tau(t,x)+\HH_\tau[v_\tau](t,x)=0,\qquad
t<T,\qquad v_\tau(T,x)=g(x).
\end{equation}
Here the soft Hamiltonian $\HH_\tau[v_\tau](t,x)$ is the minimum of the
local entropy-regularized objective over action distributions.

Policy iteration is a classical dynamic-programming method
\citep{Howard1960,PB1979PolicyIteration}, with reinforcement-learning
\citep{LP2003LSPI} and entropy-regularized variants
\citep{HZAL2018SAC,GSP2019RegularizedMDP}.
For the soft HJB equation, it alternates linear policy evaluation with
Gibbs policy improvement \citep{HWZ2025,TWZ2025}.
For entropy-regularized control, one line of theory gives
algorithm-specific convergence or error bounds:
\citet{KCKK2026} allow errors in both steps and bound the infinite-horizon
$L^2$ value error under sufficiently large discounting.
Related analyses cover $q$-learning \citep{TZ2024Regret}, gradient-flow
convergence under convexity \citep{SS2024GradientFlows}, and Langevin
policy--value dynamics \citep{FW2025,FW2026Amortized}.
Another line assumes exact evaluation and Gibbs improvement, establishing
convergence in an exploratory linear--quadratic model
\citep{WZ2020Portfolio} and for infinite-horizon soft HJB equations
\citep{HWZ2025,TWZ2025}, as well as a superexponential finite-horizon
rate \citep{MWZ2025}.
For the undiscounted finite-horizon entropy-regularized problem considered
here, how do errors in approximate policy evaluation and Gibbs
improvement propagate across iterations? How much can the regularized
cost of the final deployed policy exceed the optimal cost?

We bound this cost gap uniformly over the initial time and state.
Using transition-density gradient estimates that hold with a common
constant across the policies \citep{MPZ2021}, we derive a recursion that
retains the time dependence of the value-gradient error. This yields
explicit, summable weights for each round's evaluation residual and the
KL divergence from the returned policy to that round's Gibbs target.
Under the assumptions of the main theorem, the regularized cost of the
deployed policy converges uniformly to the optimal cost as both errors
tend to zero.

The analysis shows how errors at each round affect the final policy.
However, running the iteration in high dimensions still poses
computational difficulties. Grid-based policy evaluation becomes
expensive as the state dimension grows, and sampling from the Gibbs law
can become a bottleneck as the action dimension grows. We approximate the
value using either the deep random difference method (DRDM) of
\citet*{CFZ2025} or a physics-informed neural network \citep{RPK2019PINN},
and generate actions
using conditional flow matching \citep{LCBNL2023} or Langevin sampling
\citep{Dalalyan2017LMC}. 
Numerical experiments demonstrate strong performance on finite-horizon
problems (Sections~\ref{subsec:experiment-exact}
and~\ref{subsec:experiment-no-exact}), as well as the framework's
applicability to infinite-horizon and unregularized HJB equations
(Sections~\ref{subsec:experiment-stationary-routing}
and~\ref{subsec:experiment-hjb3a}).

Section~\ref{sec:related-work} reviews related work.
Section~\ref{sec:preliminaries} formulates the control problem and
establishes the preliminary results. Section~\ref{sec:numerical-method}
presents exact and approximate policy iteration.
Section~\ref{sec:error-control} gives the main theorem and its proof.
Section~\ref{sec:experiments} reports the numerical experiments.
Section~\ref{sec:conclusion} summarizes the findings and discusses future work.

\section{Related work}\label{sec:related-work}

\citet{WZZ2020} formulated exploration in continuous-time control through
randomized actions and an entropy term in the objective.
\citet{TZZ2022} studied the resulting exploratory HJB equation and the
limit as the strength of regularization tends to zero. At fixed
temperature, the Gibbs formula turns policy improvement into a local
update of the action law. For the classical HJB equation,
\citet{KSS2020Howard} established convergence and stability of Howard's
policy improvement algorithm. For exploratory HJB equations,
\citet{HWZ2025,TWZ2025} established convergence of policy iteration
in discounted infinite-horizon settings.The probabilistic analysis of
\citet{MWZ2025} gives a superexponential convergence rate for exact
policy iteration on a finite horizon and, under sufficiently large
discounting, on an infinite horizon.

A second line of work updates the policy distribution continuously.
Mirror descent and entropy annealing give optimization procedures for
measure-valued controls \citep{KSSZ2025Mirror,SSZ2025}.
\citet{FW2025} use Langevin dynamics to update the policy and value
simultaneously in an infinite-horizon relaxed control problem.
Their amortized Score--Hamiltonian iteration shares a conditional policy
sampler across states and analyzes residuals of the score and the
value approximation \citep{FW2026Amortized}. In the iteration studied
here, the computed value fixes the Gibbs target for each policy update.
This permits separate control of the gradient error that perturbs the
target and the KL error of the policy returned for that target.

For numerical policy evaluation, martingale conditions characterize the
value of a fixed policy from trajectories \citep{JZ2022Evaluation}.
Related representations lead to continuous-time actor--critic algorithms
and to $q$-learning \citep{JZ2022Gradient,JZ2023q}.
\citet{MWZZ2026Learning} use probabilistic representations of values and
gradients to construct policy iteration with unknown drift and running
reward. When the control coefficients and costs are given, neural
solvers can instead approximate the frozen-policy PDE directly.
\citet{KCKK2026} use a physics-informed neural network in soft policy
iteration and bound the discounted $L^2$ value error in terms of
iteration, policy approximation, and PDE residual errors.
\citet{CFZ2025} approximate parabolic operators by random differences of
network values. Our finite-horizon bounds use the uniform residual of
the frozen-policy PDE and the KL divergence of the deployed action law
from its Gibbs target. They quantify how errors at different rounds
affect the final cost, the value gradient, and the KL divergence from
the optimal policy.

\section{Preliminary}\label{sec:preliminaries}

We make precise the entropy-regularized control problem introduced in
Section~\ref{sec:introduction}. We then use dynamic programming and
Gibbs minimization to derive the soft HJB equation and verify that its
solution gives the optimal value.

\subsection{From stochastic control to entropy-regularized control}
\label{subsec:control-setting}

\paragraph{State, actions, and ordinary controls.}
Fix a terminal time $T>0$. The state lies in $\R^d$, and the action
space $\mathcal U\subseteq\R^m$ is a nonempty Borel set equipped with
its trace $\sigma$-field
$\mathcal B(\mathcal U)=\{B\cap\mathcal U:B\in\mathcal B(\R^m)\}$.
The drift and running cost are Borel functions
$b:[0,T]\times\R^d\times\mathcal U\to\R^d$ and
$c:[0,T]\times\R^d\times\mathcal U\to\R$.
The diffusion coefficient
$\sigma:[0,T]\times\R^d\to\R^{d\times r}$ is Borel and independent
of the action, and $g:\R^d\to\R$ is the Borel terminal cost.
For an initial condition $(t,x)$ with $0\le t<T$, work on a filtered
probability space $(\Omega,\mathcal F,(\mathcal F_s)_{t\le s\le T},
\mathbb P)$ whose filtration is complete and right-continuous.
The process $W$ is an $r$-dimensional Brownian motion relative to this
filtration, started at $W_t=0$ on the shifted time interval.
States are continuous and adapted.

An ordinary, or strict, control is a progressively measurable
$\mathcal U$-valued process $(u_s)_{t\le s\le T}$. Its controlled
system and unregularized cost are
\begin{equation}\label{eq:strict-control-problem}
\begin{aligned}
\dd X_s^u
 &=b(s,X_s^u,u_s)\dd s+\sigma(s,X_s^u)\dd W_s,
 \qquad X_t^u=x,\\
J_0(t,x;u)
 &=\E^u_{t,x}\!\left[g(X_T^u)+\int_t^T c(s,X_s^u,u_s)\dd s\right].
\end{aligned}
\end{equation}
Admissibility requires a well-posed state equation and integrability
of the terminal cost and the absolute running cost. The expectation
is over the chosen controlled system starting from $(t,x)$.
An ordinary Markov feedback is a Borel map
$u:[0,T]\times\R^d\to\mathcal U$; its executed control
$u_s=u(s,X_s^u)$ is progressively measurable.

\paragraph{Relaxed controls and feedback kernels.}
Write $\PP(\mathcal U)$ for the probability measures on
$(\mathcal U,\mathcal B(\mathcal U))$, equipped with the $\sigma$-field
generated by the maps $\lambda\mapsto\lambda(B)$ for
$B\in\mathcal B(\mathcal U)$. An ordinary action $u$ is embedded as
the Dirac measure $\delta_u$, where $\delta_u(B)=\mathbf1_B(u)$.
For a bounded Borel function $f$ of the action,
$\int f(a)\delta_u(\dd a)=f(u)$. Relaxation permits a general
$\lambda\in\PP(\mathcal U)$ instead of this concentrated law
\citep[Section~3.1]{Lacker2015Controlled}.

A relaxed control is a progressively measurable
$\PP(\mathcal U)$-valued process $(\pi_s)_{t\le s\le T}$.
Explicitly, for every $a\in[t,T]$ and $B\in\mathcal B(\mathcal U)$,
$(s,\omega)\mapsto\pi_s(\omega)(B)$ restricted to $[t,a]\times\Omega$
is $\mathcal B([t,a])\otimes\mathcal F_a$-measurable.
Its dynamics are
\begin{equation}\label{eq:general-relaxed-dynamics}
\dd X_s=\left[\int_{\mathcal U}b(s,X_s,u)\pi_s(\dd u)\right]\dd s
       +\sigma(s,X_s)\dd W_s,\qquad X_t=x.
\end{equation}
Thus the action law determines the mean drift; the diffusion is unchanged
because it is action-independent. This averaged dynamics is also used
in exploratory stochastic control \citep[Section~2]{WZZ2020}.
Taking $\pi_s=\delta_{u_s}$ recovers the ordinary dynamics in
\eqref{eq:strict-control-problem}.
Relaxation alone does not require a spread of actions: whenever the
unregularized local action objective attains its minimum, a Dirac law
at a minimizing action also attains the minimum over distributions.
The relative-entropy penalty below introduces a preference over these
laws in addition to their average control cost.

A Markov feedback kernel $\pi(\dd u\mid s,y)$ assigns a probability
measure to each $(s,y)$, with $(s,y)\mapsto\pi(B\mid s,y)$ Borel for
each Borel $B$. It induces the relaxed control
$\pi_s=\pi(\cdot\mid s,X_s)$ and the feedback state equation
\eqref{eq:relaxed-dynamics}. Its averaged drift is
$\bar b^\pi(s,y)=\int_{\mathcal U}b(s,y,u)\pi(\dd u\mid s,y)$.
A general relaxed control also induces the random finite measure
\[
\Lambda^\pi(\dd s,\dd u)=\dd s\,\pi_s(\dd u)
\]
on the product Borel space $[t,T]\times\mathcal U$.
Its time marginal is Lebesgue measure. Versions agreeing
$\dd s\otimes\mathbb P$-almost everywhere induce the same random
control measure almost surely.

\paragraph{Relative entropy and admissibility.}
Fix $\tau>0$ and a Borel reference probability kernel $\mu_{s,y}$
with full support on $\mathcal U$. The reference may have atoms and
may depend on time and state; it and $\tau$ remain fixed throughout
the iteration. For two probability measures on $\mathcal U$,
\[
\KL(\pi\Vert\mu)
=\int_{\mathcal U}\log\!\left(\frac{\dd\pi}{\dd\mu}\right)\dd\pi
=\int_{\mathcal U}f\log f\,\dd\mu,
\qquad f=\frac{\dd\pi}{\dd\mu},
\]
when $\pi\ll\mu$, with $0\log0=0$, and the divergence is $+\infty$
otherwise. The derivative is the Radon--Nikodym derivative on
$(\mathcal U,\mathcal B(\mathcal U))$. For feedback kernels,
$\KL(\pi\Vert\rho)(s,y)$ compares their action laws at the same
$(s,y)$.

The regularized cost of a general relaxed controlled system is
\begin{equation}\label{eq:general-regularized-cost}
J_\tau(t,x;\pi)=\E_{\mathbb P}\!\left[
g(X_T)+\int_t^T\left\{
\int_{\mathcal U}c(s,X_s,u)\pi_s(\dd u)
+\tau\KL(\pi_s\Vert\mu_{s,X_s})\right\}\dd s\right].
\end{equation}
For Markov feedback, this is \eqref{eq:regularized-cost}.
In the weak formulation, an element of $\mathcal A_\tau(t,x)$
comprises the filtered probability space, Brownian motion, progressive
control, and compatible state solving \eqref{eq:general-relaxed-dynamics}.
The filtered space may vary across systems, and we denote each system
by its control $\pi$. Admissible control specifications have well-posed
dynamics; their drift and stochastic integrals are well defined and
the cost terms satisfy
\[
\E_{\mathbb P}\!\left[|g(X_T)|+
\int_t^T\int_{\mathcal U}|c(s,X_s,u)|\pi_s(\dd u)\dd s\right]<\infty,
\qquad
\E_{\mathbb P}\int_t^T\KL(\pi_s\Vert\mu_{s,X_s})\dd s<\infty.
\]
These requirements also apply to the closed-loop dynamics of the
feedback laws used for iteration. Such progressive-measurability and
integrability requirements are standard for measure-valued controls
\citep[Section~2.1]{KSSZ2025Mirror}. The last inequality implies
$\pi_s\ll\mu_{s,X_s}$ for $\dd s\otimes\mathbb P$-almost every
$(s,\omega)$. If the reference is atomless, a strict action
$\delta_{u_s}$ has infinite relative entropy. The Dirac embedding
therefore recovers the ordinary dynamics and original running cost,
but need not be admissible for the regularized problem at fixed
$\tau>0$. The KL penalty changes the objective; the optimum in
\eqref{eq:intro-optimal-value} is the regularized optimum.

\paragraph{Expectations under feedback policies.}
For a general adapted control, the expectation in
\eqref{eq:general-regularized-cost} averages over both the state and
the control. Under a Markov feedback, the terminal cost and accumulated
regularized running cost are functions of the state path. When the
feedback state equation is well posed in law, write
$\mathbb P^\pi_{t,x}=\mathbb P\circ(X^\pi)^{-1}$ for its law on
$C([t,T];\R^d)$, with the Borel $\sigma$-field for the uniform norm.
The coordinate process $X_s(\omega)=\omega(s)$ has natural filtration
generated by $(X_a)_{t\le a\le s}$, with the usual augmentation under
$\mathbb P^\pi_{t,x}$. For an integrable Borel path functional $F$,
\[
\E^\pi_{t,x}[F(X)]
=\int_{C([t,T];\R^d)}F(\omega)\mathbb P^\pi_{t,x}(\dd\omega).
\]
At a fixed time $s$, an auxiliary action can be sampled from
$\pi(\cdot\mid s,X_s)$ using fresh randomness independent of the path.
The joint law of the path and this action is
$\mathbb P^\pi_{t,x}(\dd\omega)\pi(\dd u\mid s,\omega(s))$.
This distinguishes the averaging over actions inside the running cost
from the outer averaging over trajectories.

\subsection{Derivation and verification of the soft HJB equation}
\label{subsec:soft-hjb-derivation}

For a smooth candidate value function
$w:[0,T]\times\R^d\to\R$, let $D_xw\in\R^d$ and
$D_x^2w\in\R^{d\times d}$ denote its gradient and Hessian.
Set $A(t,x)=\sigma(t,x)\sigma(t,x)^\top\in\R^{d\times d}$.
The fixed-action generator is
\begin{equation}\label{eq:generator}
\LL^u w(t,x)=b(t,x,u)\cdot D_xw(t,x)
+\frac12\tr\!\left(A(t,x)D_x^2w(t,x)\right).
\end{equation}
It describes the instantaneous expected change in $w$ caused by the
state dynamics, with the explicit time derivative $w_t$ kept separate.
Averaging under a feedback policy gives
\begin{equation}\label{eq:feedback-generator}
\LL^\pi w(t,x):=\int_{\mathcal U}\LL^u w(t,x)\pi(\dd u\mid t,x)
=\bar b^\pi(t,x)\cdot D_xw(t,x)
+\frac12\tr\!\left(A(t,x)D_x^2w(t,x)\right).
\end{equation}
The regularized running cost rate under the same policy is
\begin{equation}\label{eq:policy-running-cost}
\ell^\pi(t,x)=\int_{\mathcal U}c(t,x,u)\pi(\dd u\mid t,x)
+\tau\KL(\pi(\cdot\mid t,x)\Vert\mu_{t,x}).
\end{equation}
Thus $\LL^\pi w$ is the dynamics contribution to value change, whereas
$\ell^\pi$ is the running cost already present in
\eqref{eq:regularized-cost}; it does not depend on $w$ and does not
include the terminal cost. At a fixed $(t,x)$, we write
$\lambda\in\PP(\mathcal U)$ for a candidate action distribution such as
$\pi(\cdot\mid t,x)$, the distribution prescribed there by a feedback
policy $\pi$. The same integrals define $\LL^\lambda$ and $\ell^\lambda$
for such a law.
The Hamiltonian combines the running cost and dynamics:
\begin{equation}\label{eq:hamiltonian}
\HH[w](t,x,u):=c(t,x,u)+\LL^u w(t,x).
\end{equation}
The local objective for an action law $\lambda$ at $(t,x)$ is
\begin{equation}\label{eq:local-policy-objective}
\LL^\lambda w(t,x)+\ell^\lambda(t,x)
=\int_{\mathcal U}\HH[w](t,x,u)\lambda(\dd u)
+\tau\KL(\lambda\Vert\mu_{t,x}).
\end{equation}

At fixed $(t,x)$ and $w$, the diffusion contribution in $\HH[w]$ is
independent of the action. Define the action-dependent energy by
\begin{equation}\label{eq:action-energy}
E_w(t,x,u)=c(t,x,u)+b(t,x,u)\cdot D_xw(t,x).
\end{equation}
Then $\HH[w]=E_w+\tfrac12\tr(A D_x^2w)$.
The two terms in $E_w$ are the immediate running cost and the drift's
contribution to future value. Removing the common diffusion term changes
the local minimum value, but not its minimizing action distribution.
For a given candidate $w$ with a positive, finite normalizing integral,
define the Gibbs action law
\begin{equation}\label{eq:improvement-map}
q_w(\dd u\mid t,x)=
\frac{\exp[-E_w(t,x,u)/\tau]\,\mu_{t,x}(\dd u)}
{\displaystyle\int_{\mathcal U}\exp[-E_w(t,x,a)/\tau]\,\mu_{t,x}(\dd a)}.
\end{equation}
The normalized exponential is a density relative to $\mu_{t,x}$, so
the formula also applies when the reference measure has atoms.

\begin{lemma}[Gibbs identity]
\label{lem:gibbs-residual-link}
Fix $(t,x)$ and a smooth candidate $w$. Assume that the normalizing
integral in \eqref{eq:improvement-map} is positive and finite and that
$q_w$ has finite absolute Hamiltonian integral and finite relative
entropy with respect to $\mu_{t,x}$. For every action law $\lambda$
with the same integrability properties,
\begin{equation}\label{eq:gibbs-identity}
\int_{\mathcal U}\HH[w]\dd\lambda
+\tau\KL(\lambda\Vert\mu_{t,x})
=-\tau\log\int_{\mathcal U}\exp[-\HH[w]/\tau]\dd\mu_{t,x}
+\tau\KL(\lambda\Vert q_w).
\end{equation}
\end{lemma}

\begin{proof}
The action-independent diffusion factor cancels in the normalized
Gibbs law, so using $\HH[w]$ or $E_w$ gives the same $q_w$.
Write $Z=\int_{\mathcal U}\exp[-\HH[w]/\tau]\dd\mu_{t,x}$.
The Gibbs density is positive relative to $\mu_{t,x}$, so
$\lambda\ll\mu_{t,x}$ implies $\lambda\ll q_w$. The chain rule for
densities and the formula for $q_w$ give, $\lambda$-almost everywhere,
\[
\log\frac{\dd\lambda}{\dd q_w}
=\log\frac{\dd\lambda}{\dd\mu_{t,x}}
-\log\frac{\dd q_w}{\dd\mu_{t,x}}
=\log\frac{\dd\lambda}{\dd\mu_{t,x}}+\frac{\HH[w]}{\tau}+\log Z.
\]
Integrating with respect to $\lambda$, multiplying by $\tau$, and
rearranging gives \eqref{eq:gibbs-identity}
\citep[Proposition~2.5]{BD1998Variational}.
\end{proof}

Assume the dynamic programming principle for the admissible systems
of Section~\ref{subsec:control-setting}
\citep[Sections~III.7--III.8]{FS2006ControlledMarkov}.
Let $v_\tau$ and the verification candidate $w$ be $C^{1,2}$
for $t<T$, continuous up to $T$, and satisfy the hypotheses of
Lemma~\ref{lem:gibbs-residual-link} pointwise.
Along admissible systems, assume It\^o's formula is valid in
expectation, with all terms integrable, on short intervals for
$v_\tau$ and on $[t,T]$ for $w$.
For each $(t,x)$ with $t<T$, take the infimum of
$(v_\tau)_t+\LL^\lambda v_\tau+\ell^\lambda$
over the lemma's action laws.
Assume this infimum can be approached by locally admissible
constant laws.
Under each such law, the integrand converges in $L^1$
to its value at $(t,x)$ as $s\downarrow t$.
Along a family with short-time dynamic programming cost gap $o(h)$,
the infimum's expected time average over $[t,t+h]$ converges
to its value at $(t,x)$ as $h\downarrow0$.
For verification, also require the lemma's integrability
conditions for $\lambda=\pi_s$ to hold
$\dd s\otimes\mathbb P$-almost everywhere along every
admissible system.

\begin{theorem}
\label{thm:soft-hjb-derivation-verification}
For a candidate $w$ satisfying the hypotheses of
Lemma~\ref{lem:gibbs-residual-link}, define the soft Hamiltonian by
\begin{equation}\label{eq:soft-hamiltonian}
\HH_\tau[w](t,x):=\min_\lambda
\{\LL^\lambda w(t,x)+\ell^\lambda(t,x)\}
=-\tau\log\int_{\mathcal U}\exp[-\HH[w](t,x,u)/\tau]\mu_{t,x}(\dd u),
\end{equation}
where the minimum is over the action laws in the lemma and is attained
at $\lambda=q_w(\cdot\mid t,x)$.
Under the preceding conditions for the classical derivation, the optimal
value $v_\tau$ is a classical solution of the soft HJB terminal-value
problem
\[
\partial_t w(t,x)+\HH_\tau[w](t,x)=0,\qquad
t<T,\qquad w(T,x)=g(x).
\]
Conversely, let $w$ be a classical solution of this problem satisfying
the preceding verification conditions. If $q_w$ is an admissible
feedback policy, then $w=v_\tau$ and
$\pi^*=q_w$ is optimal.
\end{theorem}

\begin{proof}
\textbf{(I)}
For $0<h<T-t$, dynamic programming splits the cost into the cost on
$[t,t+h]$ and the optimal continuation value:
\begin{equation}\label{eq:short-time-dpp}
v_\tau(t,x)=\inf_{\pi\in\mathcal A_\tau(t,x)}
\E_{\mathbb P}\!\left[
\int_t^{t+h}\ell^{\pi_s}(s,X_s)\dd s+v_\tau(t+h,X_{t+h})\right].
\end{equation}
Apply It\^o's formula to the continuation value with the current
action law $\pi_s$ in the generator. Subtract $v_\tau(t,x)$ and
divide by $h$ to get
\[
0=\inf_{\pi\in\mathcal A_\tau(t,x)}
\E_{\mathbb P}\!\left[\frac1h\int_t^{t+h}
\big[(v_\tau)_t+\LL^{\pi_s}v_\tau+\ell^{\pi_s}\big](s,X_s)\dd s\right].
\]
Each admissible system gives a nonnegative expected average. The
constant action laws in the first short-time condition give
nonnegative $L^1$ limits and approach the pointwise infimum at $(t,x)$.
This infimum is therefore nonnegative.

For the near-optimal controls in the second short-time condition,
the $o(h)$ cost gap in \eqref{eq:short-time-dpp} becomes an $o(1)$
gap in the expected average above. Taking the pointwise infimum inside
the integral bounds that average from below:
\[
\E_{\mathbb P}\!\left[\frac1h\int_t^{t+h}
\inf_\lambda\big[(v_\tau)_t+\LL^\lambda v_\tau+\ell^\lambda\big](s,X_s)
\dd s\right]\le o(1).
\]
By the second short-time condition, the left side converges to the
pointwise infimum at $(t,x)$, which is therefore also nonpositive.
Thus $\inf_\lambda[(v_\tau)_t+\LL^\lambda v_\tau+\ell^\lambda](t,x)=0$.
Taking $(v_\tau)_t$ outside the infimum and using
\eqref{eq:local-policy-objective} with $w=v_\tau$ gives
\begin{equation}\label{eq:hjb-local-minimum}
\partial_t v_\tau(t,x)+
\inf_\lambda
\left\{\int_{\mathcal U}\HH[v_\tau](t,x,u)\lambda(\dd u)
+\tau\KL(\lambda\Vert\mu_{t,x})\right\}=0.
\end{equation}
In the Gibbs identity \eqref{eq:gibbs-identity}, the KL term is
nonnegative and vanishes at $\lambda=q_w(\cdot\mid t,x)$. Thus the
logarithmic term is the minimum in \eqref{eq:soft-hamiltonian}. Taking
$w=v_\tau$ and substituting this minimum into \eqref{eq:hjb-local-minimum}
gives $(v_\tau)_t+\HH_\tau[v_\tau]=0$ for $t<T$. At $t=T$, the cost
reduces to $g(x)$, which gives the terminal condition.

\textbf{(II)}
Let $w$ be a classical solution of the terminal-value problem in the
theorem. For any admissible system, It\^o's formula, the zero mean of
the stochastic integral, and $w(T,\cdot)=g$ give
\[
\E_{\mathbb P}[g(X_T)]-w(t,x)
=\E_{\mathbb P}\!\int_t^T
\big[w_t+\LL^{\pi_s}w\big](s,X_s)\dd s.
\]
Adding the expected running cost to both sides gives
$J_\tau(t,x;\pi)-w(t,x)$ on the left and adds $\ell^{\pi_s}$ to
the integrand on the right.
By \eqref{eq:local-policy-objective}, the Gibbs identity
\eqref{eq:gibbs-identity}, and \eqref{eq:soft-hamiltonian}, this
integrand equals $w_t+\HH_\tau[w]+\tau\KL(\pi_s\Vert q_w)$.
Since $w_t+\HH_\tau[w]=0$, we obtain
\begin{equation}\label{eq:verification-identity}
J_\tau(t,x;\pi)-w(t,x)
=\tau\E_{\mathbb P}\!\left[\int_t^T
\KL\!\left(\pi_s\Vert q_w(\cdot\mid s,X_s)\right)\dd s\right]\ge0.
\end{equation}
The expectation is under the given admissible system. Taking the
infimum over admissible systems gives $w(t,x)\le v_\tau(t,x)$.
For the admissible feedback $q_w$, the KL term vanishes along its
trajectories, so $v_\tau(t,x)\le J_\tau(t,x;q_w)=w(t,x)$.
Hence $w=v_\tau$ and $\pi^*=q_w=q_{v_\tau}$ attains the optimum.
\end{proof}

\begin{corollary}[Performance-gap identity]
\label{cor:performance-gap}
Under the verification hypotheses of
Theorem~\ref{thm:soft-hjb-derivation-verification}, every admissible
feedback policy $\pi$ satisfies
\begin{equation}\label{eq:performance-gap-identity}
J_\tau(t,x;\pi)-v_\tau(t,x)
=\tau\E^\pi_{t,x}\!\left[\int_t^T
\KL(\pi\Vert\pi^*)(s,X_s)\dd s\right].
\end{equation}
Here $\pi^*=q_{v_\tau}$, and the expectation is along the state
trajectories under $\pi$.
\end{corollary}

\begin{proof}
Set $w=v_\tau$ and $q_w=\pi^*$ in \eqref{eq:verification-identity}.
For a feedback policy $\pi$, the integrand is a function of the state
path, so the expectation is $\E^\pi_{t,x}$.
\end{proof}


\section{Policy and value iteration}
\label{sec:numerical-method}

The control formulation and Gibbs minimization in
Section~\ref{sec:preliminaries} lead to an iteration for solving the
soft HJB problem. Starting from an admissible feedback policy $\pi_0$,
each exact update takes the form
\begin{equation}\label{eq:pi-improvement}
\pi_k\xrightarrow{\text{policy evaluation}}v_k=v^{\pi_k}
\xrightarrow{\text{Gibbs improvement}}\pi_{k+1}=q_{v_k}.
\end{equation}
The evaluation step computes the cost of the current policy; the
improvement step uses that value to choose a new action law at each
time and state. The new policy is then evaluated in the next round.

We first establish the linear evaluation equation, use the Gibbs
identity from Section~\ref{subsec:soft-hjb-derivation} to identify the
improved policy, and show that the update decreases the total cost.
Section~\ref{subsec:numerical-iteration} then describes how numerical
PDE and policy solvers approximate these two steps.

\subsection{Exact policy and value iteration}\label{sec:iteration}

\paragraph{Fixed-policy evaluation.}
To evaluate $\pi_k$, hold its feedback law fixed and compute
$v_k(t,x)=J_\tau(t,x;\pi_k)$. The following lemma gives the equation
for this cost under a fixed policy $\pi$. It uses the averaged
generator $\LL^\pi$ and regularized running cost $\ell^\pi$ defined
in \eqref{eq:feedback-generator}--\eqref{eq:policy-running-cost}.

\begin{lemma}[Fixed-policy evaluation equation]
\label{lem:fixed-policy-evaluation}
Fix a Markov feedback policy $\pi$ that is admissible for every initial
time and state in the sense of Section~\ref{subsec:control-setting},
and whose state laws form a Markov family. Assume that
$\bar b^\pi$, $A$, and $\ell^\pi$ are continuous on
$[0,T)\times\R^d$, and that
$v^\pi(t,x):=J_\tau(t,x;\pi)$ belongs to
$C^{1,2}([0,T)\times\R^d)\cap C([0,T]\times\R^d)$.
Then
\begin{equation}\label{eq:pi-evaluation}
(v^\pi)_t+\LL^{\pi}v^\pi+\ell^{\pi}=0,\qquad t<T,
\qquad v^\pi(T,x)=g(x).
\end{equation}
\end{lemma}

\begin{proof}
Fix any $(t,x)$ with $t<T$, and let $X_t=x$ under the feedback state
law $\mathbb P^\pi_{t,x}$. Write $\mathcal F_s^X$ for the augmented
natural filtration from Section~\ref{subsec:control-setting}.
Admissibility makes the total cost integrable. Splitting that cost at
$s$ and using the Markov property gives
\[
\begin{aligned}
M_s:=v^\pi(s,X_s)+\int_t^s\ell^\pi(r,X_r)\dd r=\E^\pi_{t,x}\!\left[
g(X_T)+\int_t^T\ell^\pi(r,X_r)\dd r
\,\middle|\,\mathcal F_s^X\right].
\end{aligned}
\]
Thus the accumulated cost plus the remaining value is a continuous
martingale.

On any interval $[t,T']$ with $t<T'<T$, It\^o's formula under the
feedback state law yields
\[
M_s-M_t=\int_t^s
\big[(v^\pi)_t+\LL^\pi v^\pi+\ell^\pi\big](r,X_r)\dd r+N_s,
\]
where $N_t=0$ and $N$ is a continuous local martingale. Localization
on compact state sets justifies this decomposition without a global
integrability condition on the derivatives of $v^\pi$.
The time integral is therefore both a continuous finite-variation
process and a local martingale starting at zero, so it vanishes
almost surely, simultaneously for all $s\in[t,T']$.

The integrand is continuous in time and state by the hypotheses.
Since $X_r\to x$ as $r\downarrow t$, its short-time average satisfies,
almost surely,
\[
0=\lim_{h\downarrow0}\frac1h\int_t^{t+h}
\big[(v^\pi)_t+\LL^\pi v^\pi+\ell^\pi\big](r,X_r)\dd r
=\big[(v^\pi)_t+\LL^\pi v^\pi+\ell^\pi\big](t,x).
\]
This is a pathwise limit, so no interchange of expectation and limit
is needed. The initial pair $(t,x)$ was arbitrary, which proves the
PDE at every point with $t<T$. At $t=T$, the cost definition gives
$v^\pi(T,x)=g(x)$.
\end{proof}

With $\pi$ fixed, both $\LL^\pi$ and $\ell^\pi$ are prescribed, so
\eqref{eq:pi-evaluation} is a linear, generally inhomogeneous PDE
in $v^\pi$. In particular, the entropy term belongs to the prescribed
running cost and introduces no nonlinearity in the unknown value.
Applying the lemma with $\pi=\pi_k$ gives the evaluation step in
\eqref{eq:pi-improvement}.

For a candidate value function $w$, define the policy evaluation
residual by
\begin{equation}\label{eq:residual-definitions}
\mathcal R^\pi[w]:=w_t+\LL^\pi w+\ell^\pi.
\end{equation}
Exact evaluation gives $\mathcal R^{\pi_k}[v_k]=0$; the same residual
will measure the error of a numerical evaluation in
Section~\ref{subsec:numerical-iteration}.

\paragraph{Gibbs policy improvement.}
Once $v_k$ is available, hold it fixed and use its spatial gradient
$D_xv_k$ to form the action energy $E_{v_k}$ in
\eqref{eq:action-energy}. Set $\pi_{k+1}=q_{v_k}$ by the Gibbs formula
\eqref{eq:improvement-map}. Under the hypotheses of
Lemma~\ref{lem:gibbs-residual-link}, this choice minimizes
$\LL^\lambda v_k+\ell^\lambda$ over action laws $\lambda$ at each
$(t,x)$. The new feedback must be admissible for the next evaluation.

Under the same hypotheses, the Gibbs identity
\eqref{eq:gibbs-identity} rewrites the local objective
$\LL^\pi w+\ell^\pi$ as $\HH_\tau[w]+\tau\KL(\pi\Vert q_w)$.
Adding $w_t$ as in the residual definition
\eqref{eq:residual-definitions} gives
\begin{equation}\label{eq:residual-decomposition}
\mathcal R^\pi[w]=w_t+\HH_\tau[w]+\tau\KL(\pi\Vert q_w).
\end{equation}

We now evaluate $\pi_{k+1}$ to obtain $v_{k+1}=v^{\pi_{k+1}}$ and
compare the total costs of the two policies from the same initial
condition.

\begin{proposition}[Monotone improvement of the exact policy iteration]
\label{prop:exact-policy-improvement}
Suppose both policies are admissible, their values solve the respective
evaluation equations classically with terminal data $g$, and the Gibbs
identity applies to $w=v_k$ for both policies. For It\^o's formula applied
to $v_k-v_{k+1}$ along $\pi_{k+1}$ from any $(t,x)$, assume the drift is
absolutely integrable in time and expectation and the stochastic
integral extends to a martingale on $[t,T]$. Then
\begin{equation}\label{eq:pi-monotonicity}
v_k(t,x)-v_{k+1}(t,x)
=\tau\E^{\pi_{k+1}}_{t,x}\int_t^T
\KL(\pi_k\Vert\pi_{k+1})(s,X_s)\dd s\ge0,
\end{equation}
where the expectation is over state trajectories under $\pi_{k+1}$.
In particular, $v_\tau\le v_{k+1}\le v_k$.
\end{proposition}

\begin{proof}
First hold the old value $v_k$ fixed and compare the two policies in
\eqref{eq:residual-decomposition}. Since
$q_{v_k}=\pi_{k+1}$ and exact evaluation gives
$\mathcal R^{\pi_k}[v_k]=0$,
\[
\mathcal R^{\pi_{k+1}}[v_k]
=\mathcal R^{\pi_k}[v_k]-\tau\KL(\pi_k\Vert\pi_{k+1})
=-\tau\KL(\pi_k\Vert\pi_{k+1})\le0.
\]
Next hold the new policy $\pi_{k+1}$ fixed and compare $v_k$ with
its exact value $v_{k+1}$. Subtracting
$\mathcal R^{\pi_{k+1}}[v_{k+1}]=0$ uses the same generator and
cancels the same running cost $\ell^{\pi_{k+1}}$, yielding
\[
(\partial_t+\LL^{\pi_{k+1}})(v_k-v_{k+1})
=-\tau\KL(\pi_k\Vert\pi_{k+1}),\qquad
(v_k-v_{k+1})(T,\cdot)=0.
\]
It\^o's formula along the new policy gives
\[
0-(v_k-v_{k+1})(t,x)
=-\tau\E^{\pi_{k+1}}_{t,x}\int_t^T
\KL(\pi_k\Vert\pi_{k+1})(s,X_s)\dd s.
\]
The zero is the terminal value difference. Rearranging proves
\eqref{eq:pi-monotonicity}; the lower bound by $v_\tau$ follows from
its definition as the optimal cost.
\end{proof}

For convergence of exact policy iteration in a finite-horizon
setting under additional assumptions, see \citet[Theorem~2.4]{MWZ2025}.
Their proof derives a recursion for the value-gradient error on short
time intervals and extends the estimate backward from the terminal
time to the full horizon. Section~\ref{subsec:numerical-iteration}
describes how we approximate the solution of the evaluation PDE and
the Gibbs action law in practice.

\subsection{Numerical policy and value iteration}\label{subsec:numerical-iteration}

One numerical round starts from the policy $\pi_k$ used to generate
actions and returns the policy $\pi_{k+1}$ for the next round:
\[
\pi_k\xrightarrow{\text{approximate evaluation}}\hat v_k
\xrightarrow{\text{Gibbs target}}\hat q_k
\xrightarrow{\text{approximate sampling}}\pi_{k+1}.
\]

\paragraph{Approximate policy evaluation.}
Hold $\pi_k$ fixed and approximate its value $v_k$ by solving the
linear evaluation equation \eqref{eq:pi-evaluation} numerically. The
solver returns $\hat v_k$ with terminal data $\hat v_k(T,\cdot)=g$.
For the classical residual, take $\hat v_k$ to be $C^{1,2}$ for $t<T$.
To measure how far $\hat v_k$ departs from the evaluation equation,
define the uniform evaluation error by
\begin{equation}\label{eq:pde-solver-error}
E_{\mathrm{PDE}}^k
:=\sup_{(t,x)\in[0,T)\times\R^d}
  \left|\mathcal R^{\pi_k}[\hat v_k](t,x)\right|.
\end{equation}
The residual \eqref{eq:residual-definitions} uses the true generator
$\LL^{\pi_k}$ and running cost $\ell^{\pi_k}$, even if the solver
estimates the averaged drift or running cost. Errors in those
estimates are therefore included in the evaluation error.

The experiments in Section~\ref{sec:experiments} use two PDE solvers.
A physics-informed neural network (PINN) uses automatic differentiation
to form a residual loss at sampled time and state points
\citep{RPK2019PINN,KCKK2026}. The deep random difference method (DRDM)
approximates $(\partial_t+\LL^{\pi_k})\hat v_k$ by finite differences of
network values at the starting point and at random points one time step
later \citep{CFZ2025}. Appendix~\ref{app:solver-implementations} gives
the training losses and implementation details.

\paragraph{Approximate Gibbs improvement.}
Hold $\hat v_k$ fixed and construct $E_{\hat v_k}$ from its spatial
gradient using \eqref{eq:action-energy}. The Gibbs formula
\eqref{eq:improvement-map} defines the target
\begin{equation}\label{eq:critic-teacher}
\hat q_k=q_{\hat v_k}.
\end{equation}
The policy solver approximates this target and returns an admissible
action law $\pi_{k+1}$. An error in $\hat v_k$ can already shift the
target from $q_{v_k}$ to $\hat q_k$. To measure how far the returned
law departs from its given target, define the uniform policy solver
error by
\begin{equation}\label{eq:policy-solver-error}
E_{\mathrm{policy}}^k
:=\sup_{(t,x)\in[0,T)\times\R^d}\KL(\pi_{k+1}\Vert\hat q_k)(t,x).
\end{equation}
Under the hypotheses of Lemma~\ref{lem:gibbs-residual-link}, the
pointwise divergence multiplied by $\tau$ equals the excess in the
local objective for the fixed value $\hat v_k$.

The finite-horizon comparisons in Section~\ref{sec:experiments} use
Langevin sampling and conditional flow matching. When the reference
has a smooth positive density and the energy is differentiable,
Langevin sampling takes a finite number of noisy steps driven by the
gradient of the log Gibbs density \citep{Dalalyan2017LMC}. Conditional
flow matching fits a velocity field using base draws and approximate
Gibbs samples, then generates an action by solving the resulting ODE
from a fresh base draw \citep{LCBNL2023}.

\paragraph{Passing a policy to the next evaluation.}
The policy $\pi_{k+1}$ returned by the improvement step is the input
to the next evaluation. Its averaged drift and mean action cost can
be estimated from action samples. The regularized running cost
\eqref{eq:policy-running-cost} also needs the entropy term
$\KL(\pi_{k+1}\Vert\mu_{t,x})$, which can be evaluated using the
log-density ratio or estimated separately when only samples are
available.

After $N$ updates, the algorithm returns $\pi_N$. The last value
estimate $\hat v_{N-1}$ evaluates $\pi_{N-1}$, so obtaining an
approximate value for $\pi_N$ requires one more evaluation.
Section~\ref{sec:error-control} bounds how the two solver errors in
\eqref{eq:pde-solver-error} and \eqref{eq:policy-solver-error}
accumulate and affect the final policy and its cost.


\section{Theoretical results}\label{sec:error-control}

In the numerical iteration of Section~\ref{subsec:numerical-iteration},
inexact policy evaluation perturbs the gradient used to construct the
Gibbs target, and inexact policy improvement can make the deployed
policy deviate from that target. Under the assumptions below, we give
uniform bounds after $N$ updates on the policy error in KL divergence,
the value error, and the value-gradient error. These bounds quantify
the contribution of the initial gradient error and of both solver
errors at each update through explicit, summable weights. In particular,
if both solver errors tend to zero, the deployed policy's excess cost
converges uniformly to zero.

\subsection{Main results}\label{subsec:main-performance}

We use the numerical iteration and the uniform solver errors
$E_{\mathrm{PDE}}^k$ and $E_{\mathrm{policy}}^k$ from
Section~\ref{subsec:numerical-iteration}. Write $v_k=v^{\pi_k}$ for
the exact value of the deployed policy $\pi_k$. For randomized solvers,
we fix the realized sequence of value approximations and policy kernels.
Value norms are evaluated over $[0,T]\times\R^d$. Gradient and residual
norms and policy divergences are evaluated over $[0,T)\times\R^d$.
To control how the solver errors propagate across updates, we impose
the following regularity conditions on the model and the iterates.

\begin{assumption}[Regularity for the outer iteration]\label{ass:outer-regularity}
Consider the finite-horizon model in \eqref{eq:relaxed-dynamics} and
\eqref{eq:regularized-cost} with fixed horizon $T$, temperature $\tau$,
and reference kernel $\mu_{t,x}$. The coefficients are Borel
measurable. The action-independent diffusion $\sigma$ is bounded and
Lipschitz in space uniformly in time, with
$A=\sigma\sigma^\top\succeq a_0I_d$ for some $a_0>0$.
The drift and cost satisfy $|b(t,x,u)|\le B$ and
$|c(t,x,u)|\le C_c$, with $B>0$. The policies $\pi_k$ and
$\pi^*=q_{v_\tau}$ are admissible with well-posed closed-loop
dynamics. The optimal value $v_\tau$ and
the exact policy values $v_k$ solve their respective equations
classically. These values and $\hat v_k$ are bounded, of class
$C^{1,2}$ for $t<T$, have bounded spatial gradients, and extend
continuously to the terminal value $g$.
\end{assumption}

We use transition operators to estimate the spatial gradients of
value differences represented as expected integrals along state
trajectories.

\begin{definition}[Transition operator]\label{def:transition-operator}
For a fixed Markov policy $\pi$, times $0\le t<s\le T$, and a bounded
measurable function $f:\R^d\to\R$, define
\[
P^\pi_{t,s}f(x):=\E^\pi_{t,x}[f(X_s)]
=\int_{\R^d} f(y)\,\mathbb P^\pi_{t,x}(X_s\in\dd y).
\]
Here $\mathbb P^\pi_{t,x}(X_s\in\cdot)$ is the time-$s$ state marginal
of the path law started at $(t,x)$.
\end{definition}

\begin{lemma}[Uniform transition gradient estimate]
\label{lem:uniform-smoothing}
Under Assumption~\ref{ass:outer-regularity}, there is a common finite
$K_T>0$ such that for every policy $\pi$ in the iteration and every
bounded measurable $f$,
\begin{equation}\label{eq:uniform-transition-gradient}
\|D_xP^\pi_{t,s}f\|_\infty\le K_T(s-t)^{-1/2}\|f\|_\infty,
\qquad 0\le t<s\le T.
\end{equation}
The constant depends only on $T$, $d$, $B$, $a_0$, the uniform bound on
$\sigma$, and its spatial Lipschitz constant. The proof is in
Appendix~\ref{app:uniform-smoothing}.
\end{lemma}

\begin{theorem}[Error and performance bounds after $N$ updates]
\label{thm:inexact-policy-iteration}
Suppose Assumption~\ref{ass:outer-regularity} holds and both uniform
solver errors are finite. For every $N\ge1$, the exact value of the
deployed policy satisfies
\begin{equation}\label{eq:main-value}
\|v_N-v_\tau\|_\infty
\le 2BT\!\left(\lambda_{N-1}D_0+
       \sum_{k=0}^{N-2}\lambda_{N-2-k}\varepsilon_k\right)
       +\frac{\sqrt T}{2K_T}\,\varepsilon_{N-1}.
\end{equation}
Since $v_N(t,x)-v_\tau(t,x)=J_\tau(t,x;\pi_N)-J_\tau(t,x;\pi^*)\ge0$,
the left-hand side of \eqref{eq:main-value} is the supremum of the excess
cost over all initial times and states. The remaining individual bounds are
\begin{align}
\sup_{0\le t<T,\,x}\KL(\pi_N\Vert\pi^*)(t,x)
&\le \frac{2B}{\tau}\!\left(\lambda_{N-1}D_0+
       \sum_{k=0}^{N-2}\lambda_{N-2-k}\varepsilon_k\right)
       +\frac{\varepsilon_{N-1}}{2\tau K_T\sqrt T},
       \label{eq:main-policy}\\
\|D_xv_N-D_xv_\tau\|_\infty
&\le \lambda_ND_0+\sum_{k=0}^{N-1}\lambda_{N-1-k}\varepsilon_k .
       \label{eq:main-gradient}
\end{align}
Adding the value and gradient bounds to the policy KL bound scaled by
$\tau$ gives the joint estimate
\begin{equation}\label{eq:total-error-master}
\|v_N-v_\tau\|_\infty+\tau\sup_{t<T,\,x}\KL(\pi_N\Vert\pi^*)(t,x)
 +\|D_xv_N-D_xv_\tau\|_\infty
 \le \omega_N D_0+\sum_{k=0}^{N-1}\omega_{N-1-k}\varepsilon_k.
\end{equation}
Here $D_0:=\|D_x(v_0-v_\tau)\|_\infty$ is the initial gradient error,
and the solver contribution and propagation weights are
\begin{align}
&\varepsilon_k:=8BK_T^2T\,E_{\mathrm{PDE}}^k
             +2\tau K_T\sqrt T\,E_{\mathrm{policy}}^k,\qquad
\lambda_n:=\frac{(2BK_T\sqrt{\pi T})^n}{\Gamma(1+n/2)},\quad n\ge0,
\label{eq:solver-error-contribution}\\
&\omega_n:=\lambda_n+2B(T+1)\lambda_{n-1}
            +\frac{T+1}{2K_T\sqrt T}\,\mathbf1_{\{n=0\}},
\qquad n\ge0.
\label{eq:total-error-weights}
\end{align}
We set $\lambda_{-1}:=0$ and interpret an empty sum as zero.
\end{theorem}

\begin{theorem}[Convergence and persistent solver errors]
\label{thm:inexact-policy-convergence}
Under the hypotheses of Theorem~\ref{thm:inexact-policy-iteration},
the propagation weights satisfy
\begin{equation}\label{eq:propagation-summability}
\sum_{n=0}^{\infty}\lambda_n<\infty,
\qquad \sum_{n=0}^{\infty}\omega_n<\infty.
\end{equation}
If $E_{\mathrm{PDE}}^k\to0$ and $E_{\mathrm{policy}}^k\to0$ as
$k\to\infty$, the upper bounds in
\eqref{eq:main-value}--\eqref{eq:total-error-master} tend to zero.
Consequently,
\begin{equation}\label{eq:iteration-uniform-convergence}
\|v_N-v_\tau\|_\infty
+\tau\sup_{t<T,\,x}\KL(\pi_N\Vert\pi^*)(t,x)
+\|D_xv_N-D_xv_\tau\|_\infty\longrightarrow0,
\end{equation}
and the cost gap converges uniformly to zero on $[0,T]\times\R^d$.
No summability or prescribed decay rate of the solver errors is required.
More generally, if the solver error sequences are bounded, then
\begin{equation}\label{eq:persistent-error-limsup}
\limsup_{N\to\infty}\left(
\|v_N-v_\tau\|_\infty
+\tau\sup_{t<T,\,x}\KL(\pi_N\Vert\pi^*)(t,x)
+\|D_xv_N-D_xv_\tau\|_\infty\right)
\le \left(\sum_{n=0}^{\infty}\omega_n\right)
\limsup_{k\to\infty}\varepsilon_k.
\end{equation}
\end{theorem}

To approximate the value of the output policy $\pi_N$, perform one
additional PDE solve with $\pi_N$ fixed. If the solve uses exact
terminal data and returns $\hat v_N$ with the regularity in
Assumption~\ref{ass:outer-regularity}, then
\begin{equation}\label{eq:final-network-guarantee}
\begin{aligned}
\|\hat v_N-v_\tau\|_\infty
  &\le \|v_N-v_\tau\|_\infty+T E_{\mathrm{PDE}}^N,\\
\|D_x\hat v_N-D_xv_\tau\|_\infty
  &\le \|D_xv_N-D_xv_\tau\|_\infty+2K_T\sqrt T\,E_{\mathrm{PDE}}^N.
\end{aligned}
\end{equation}
The two estimates add the final evaluation error relative to
$v_N=v^{\pi_N}$, controlled by \eqref{eq:proof-critic-gradient-bound}
with $k=N$, to the errors relative to $v_\tau$ in
\eqref{eq:main-value} and \eqref{eq:main-gradient}.

The estimates apply on $\R^d$ and, by periodic lifting, on a torus.
Applying the estimates with a sampled loss requires a bound relating
that loss to the corresponding uniform residual or KL divergence.

\subsection{Proof of the main results}
\label{subsec:main-proofs}\label{subsec:iteration-mechanism}

\begin{proof}[Proof of Theorem~\ref{thm:inexact-policy-iteration}]
We first use the evaluation residual to control the gradient used in
the Gibbs update. We then bound the KL divergence of the new policy
and use the performance-gap identity to derive a recursion for the
gradient error. Iterating that recursion yields the finite-update
bounds.

\emph{Step 1. From the evaluation residual to the critic gradient.}
Define
\[
\begin{aligned}
d_k(t)&=\|D_x(v_k-v_\tau)(t,\cdot)\|_\infty,
&D_k&=\sup_{0\le t<T}d_k(t),\\
\gamma_k(t)&=\|D_x(\hat{v}_k-v_k)(t,\cdot)\|_\infty,
&R_k(t,x)&=\mathcal R^{\pi_k}[\hat{v}_k](t,x).
\end{aligned}
\]
Here $d_k(t)$ compares the exact policy-value gradient to the optimal
gradient, $D_k$ takes its supremum over physical time, and $\gamma_k(t)$
measures the fitted critic's gradient error relative to its exact policy
value. The quantity $R_k$ is the PDE evaluation residual, with
$\sup_{t<T}\|R_k(t,\cdot)\|_\infty=E_{\mathrm{PDE}}^k$.
Once $\pi_k$ and $\hat v_k$ are fixed, this residual is determined by
the frozen-policy equation.
Subtracting the exact evaluation equation
\eqref{eq:pi-evaluation} from the residual equation gives
\begin{equation}\label{eq:proof-critic-error-pde}
(\partial_t+\LL^{\pi_k})(\hat{v}_k-v_k)=R_k,
\qquad (\hat{v}_k-v_k)(T,\cdot)=0.
\end{equation}
Apply It\^o's formula to $\hat v_k-v_k$ along $\pi_k$ up to a
time $r<T$. The bounded diffusion and spatial gradients make the
stochastic integral a square-integrable martingale. After taking
expectations, let $r\uparrow T$. Boundedness and continuity at $T$
give convergence of the expected value difference at $r$ to zero,
and the bounded residual permits the limit in its integral. This yields
the Feynman--Kac representation \citep[Section~8.2]{Oksendal1998SDE}
\begin{equation}\label{eq:proof-critic-feynman-kac}
(\hat{v}_k-v_k)(t,\cdot)
=-\int_t^T P^{\pi_k}_{t,s}R_k(s,\cdot)\dd s.
\end{equation}
For this linear-PDE representation in policy iteration, see
\citet[p.~822, Eq.~(2.14)]{MWZ2025}. Here it is applied to
the residual equation \eqref{eq:proof-critic-error-pde}.
At each fixed $s>t$, apply Lemma~\ref{lem:uniform-smoothing} with
the bounded spatial function $f=R_k(s,\cdot)$. This gives
\[
\|D_xP^{\pi_k}_{t,s}R_k(s,\cdot)\|_\infty
\le\frac{K_T}{\sqrt{s-t}}\|R_k(s,\cdot)\|_\infty.
\]
The right-hand side is integrable in $s$, so truncation and dominated
convergence justify differentiating \eqref{eq:proof-critic-feynman-kac}.
Consequently,
\begin{equation}\label{eq:proof-critic-gradient-bound}
\begin{aligned}
&\gamma_k(t)\le K_T\int_t^T\frac{\|R_k(s,\cdot)\|_\infty}{\sqrt{s-t}}\dd s\le2K_T\sqrt{T-t}\,E_{\mathrm{PDE}}^k
\le2K_T\sqrt T\,E_{\mathrm{PDE}}^k,\\
&\|(\hat{v}_k-v_k)(t,\cdot)\|_\infty
\le\int_t^T\|R_k(s,\cdot)\|_\infty\dd s\le TE_{\mathrm{PDE}}^k.
\end{aligned}
\end{equation}
The value estimate uses the Markov contraction
$\|P^\pi_{t,s}f\|_\infty\le\|f\|_\infty$, which follows directly
by integration against the time-$s$ probability measure.

\emph{Step 2. A pointwise Gibbs density-ratio estimate.}
By \eqref{eq:improvement-map} and action-independent diffusion, the
Hessian term cancels from each normalized Gibbs kernel. Applying the
energy definition \eqref{eq:action-energy} to $\hat{v}_k$ and $v_\tau$,
the gradient triangle inequality and $|b|\le B$ imply
\begin{equation}\label{eq:proof-energy-perturbation}
|E_{\hat{v}_k}-E_{v_\tau}|(t,x,u)
\le B|D_x\hat{v}_k-D_xv_\tau|(t,x) \le B\big(d_k(t)+\gamma_k(t)\big).
\end{equation}
Comparing the two normalizing constants against the same reference
measure, as detailed in Appendix~\ref{app:gibbs-normalization}, gives
\begin{equation}\label{eq:gibbs-log-ratio-gradient}
\sup_x\left\|\log\frac{\dd\hat{q}_k(\cdot\mid t,x)}
                        {\dd\pi^*(\cdot\mid t,x)}\right\|_{L^\infty(\mu_{t,x})}
\le\frac{2B}{\tau}\big(d_k(t)+\gamma_k(t)\big).
\end{equation}
The factor two accounts separately for the energy and its normalizing
constant. This estimate is a direct consequence of the Gibbs formula.

\emph{Step 3. Isolate the policy improvement error.}
The nonnegative quantity $\tau\KL(\pi_{k+1}\Vert\pi^*)$
is the source in the performance-gap equation below.
Finiteness of $E_{\mathrm{policy}}^k$ and the
bounded log ratio in \eqref{eq:gibbs-log-ratio-gradient} justify the exact
change-of-reference identity
\begin{equation}\label{eq:proof-kl-reference-change}
\KL(\pi_{k+1}\Vert\pi^*)
={}\KL(\pi_{k+1}\Vert\hat{q}_k) +\int_{\mathcal U}\log\frac{\dd\hat{q}_k}{\dd\pi^*}(u\mid t,x)
\pi_{k+1}(\dd u\mid t,x).
\end{equation}
All terms are at the same $(t,x)$. Both Gibbs laws are equivalent to
$\mu_{t,x}$, and finite improvement error gives
$\pi_{k+1}(\cdot\mid t,x)\ll\hat q_k(\cdot\mid t,x)$.
Thus changing a density ratio on a $\mu_{t,x}$-null set leaves this
integral unchanged. Inserting
\eqref{eq:gibbs-log-ratio-gradient} yields
\begin{equation}\label{eq:outer-source-bound}
\begin{aligned}
\tau\sup_x\KL(\pi_{k+1}\Vert\pi^*)(t,x)
&\le2B\big(d_k(t)+\gamma_k(t)\big)+\tau E_{\mathrm{policy}}^k\\
&\le2Bd_k(t)+4BK_T\sqrt T\,E_{\mathrm{PDE}}^k
+\tau E_{\mathrm{policy}}^k.
\end{aligned}
\end{equation}

\emph{Step 4. Derive the time-dependent gradient recursion.}
The bounded energy $E_{v_\tau}$ makes the log-density ratio of
$\pi^*$ to $\mu_{s,X_s}$ uniformly bounded. For any admissible control,
its entropy relative to $\mu_{s,X_s}$ is time-integrable in expectation,
so the same holds for its entropy relative to $\pi^*$ by the change of
reference. At each time and state the Hamiltonian is bounded over actions,
so the Gibbs identity applies along the control. By that identity and
the soft HJB equation, the It\^o drift of $v_\tau(s,X_s)$ equals
\[
\tau\KL\!\left(\pi_s\Vert\pi^*(\cdot\mid s,X_s)\right)
-\ell^{\pi_s}(s,X_s).
\]
Both terms are time-integrable in expectation by the preceding bound
and admissibility. As in Step~1, the bounded diffusion and spatial
gradient make the stochastic integral a square-integrable martingale.
\ifcsname r@cor:performance-gap\endcsname
The verification conditions therefore hold, so Corollary~\ref{cor:performance-gap}
applies with $\pi=\pi_{k+1}$.
\else
Taking expectations in It\^o's formula along $\pi_{k+1}$ and
using the terminal condition gives the performance-gap identity.
\fi
Since $v_{k+1}(t,x)=J_\tau(t,x;\pi_{k+1})$, writing the trajectory
expectation in terms of the transition operator gives
\begin{equation}\label{eq:outer-feynman-kac}
v_{k+1}(t,x)-v_\tau(t,x) =\tau\int_t^T P^{\pi_{k+1}}_{t,s}
\big[\KL(\pi_{k+1}\Vert\pi^*)(s,\cdot)\big](x)\dd s\ge0.
\end{equation}
By \eqref{eq:outer-source-bound}, the KL source is bounded uniformly
in time and state. The factor $(s-t)^{-1/2}$ in
Lemma~\ref{lem:uniform-smoothing} is integrable on $(t,T)$, so the
same truncation and dominated-convergence argument as in Step~1
justifies differentiating \eqref{eq:outer-feynman-kac} under the
integral. Applying the gradient estimate and
\eqref{eq:outer-source-bound} gives
\[
\begin{aligned}
d_{k+1}(t)
&\le K_T\int_t^T\frac{2Bd_k(s)+4BK_T\sqrt T\,E_{\mathrm{PDE}}^k+\tau E_{\mathrm{policy}}^k}{\sqrt{s-t}}\dd s\\
&=2BK_T\int_t^T\frac{d_k(s)}{\sqrt{s-t}}\dd s+2K_T\sqrt{T-t}\left(4BK_T\sqrt T\,E_{\mathrm{PDE}}^k+\tau E_{\mathrm{policy}}^k\right).
\end{aligned}
\]
Using $\sqrt{T-t}\le\sqrt T$ and
\eqref{eq:solver-error-contribution} gives
\begin{equation}\label{eq:outer-volterra-recursion}
\begin{aligned}
d_{k+1}(t)&\le(\mathsf Kd_k)(t)+\varepsilon_k,\\
(\mathsf Kh)(t)&=2BK_T\int_t^T(s-t)^{-1/2}h(s)\dd s.
\end{aligned}
\end{equation}
Here $h$ is a nonnegative function of time. Iterating the time integral
produces the summable propagation weights below.

\emph{Step 5. Accumulate the history and recover policy performance.}
The Abel convolution identity \eqref{eq:abel-convolution}, proved
in Appendix~\ref{app:abel-convolution}, evaluates the iterated kernel.
For the weakly singular Volterra iteration, see also the proof of \citet[Lemma~7.1.1]{Henry1981}.
Integrating the $n$-fold kernel once more gives
\begin{equation}\label{eq:volterra-powers}
(\mathsf K^n\mathbf1)(t)
=\frac{(2BK_T\sqrt\pi)^n(T-t)^{n/2}}{\Gamma(1+n/2)},
\qquad n\ge0.
\end{equation}
Because $\mathsf K$ preserves nonnegativity, repeated substitution
in \eqref{eq:outer-volterra-recursion}, followed by
$d_0(t)\le D_0$ and \eqref{eq:volterra-powers}, gives
\begin{equation}\label{eq:two-error-gradient-sum}
D_N\le\lambda_ND_0+\sum_{k=0}^{N-1}\lambda_{N-1-k}\varepsilon_k,
\end{equation}
which is \eqref{eq:main-gradient}.

To estimate the value and policy, use \eqref{eq:outer-source-bound}
at $k=N-1$ and the Markov contraction in \eqref{eq:outer-feynman-kac}.
Since $4BK_T\sqrt T\,E_{\mathrm{PDE}}^{N-1}
+\tau E_{\mathrm{policy}}^{N-1}=\varepsilon_{N-1}/(2K_T\sqrt T)$,
\begin{equation}\label{eq:proof-value-policy-from-gradient}
\begin{aligned}
&\tau\sup_{0\le t<T,\,x}\KL(\pi_N\Vert\pi^*)(t,x)
\le 2B D_{N-1}+\frac{\varepsilon_{N-1}}{2K_T\sqrt T},\\
&\|v_N-v_\tau\|_\infty
\le T\left(2B D_{N-1}+\frac{\varepsilon_{N-1}}{2K_T\sqrt T}\right).
\end{aligned}
\end{equation}
For $N\ge2$, insert \eqref{eq:two-error-gradient-sum} at $N-1$.
For $N=1$, use $D_{N-1}=D_0$ and $\lambda_0=1$.
This gives \eqref{eq:main-value} and \eqref{eq:main-policy}, with the
last update's contribution explicit in each case. By
\eqref{eq:outer-feynman-kac}, $v_N-v_\tau$ is the nonnegative cost gap,
so its supremum norm is the uniform cost gap.

Adding \eqref{eq:proof-value-policy-from-gradient} and the gradient
error $D_N$ bounds the left side of \eqref{eq:total-error-master} by
\[
 D_N+2B(T+1)D_{N-1}+\frac{T+1}{2K_T\sqrt T}\,\varepsilon_{N-1}.
\]
By \eqref{eq:two-error-gradient-sum}, the coefficient of $D_0$ is
$\lambda_N+2B(T+1)\lambda_{N-1}=\omega_N$. The coefficient of
$\varepsilon_k$ is
\[
 \lambda_{N-1-k}+2B(T+1)\lambda_{N-2-k}
 +\frac{T+1}{2K_T\sqrt T}\,\mathbf1_{\{k=N-1\}}
 =\omega_{N-1-k},
\]
including the last update since $\lambda_{-1}=0$. This proves
\eqref{eq:total-error-master}.
\end{proof}

\begin{proof}[Proof of Theorem~\ref{thm:inexact-policy-convergence}]
The gamma recurrence gives
\[
 \frac{\lambda_{n+2}}{\lambda_n}
 =\frac{(2BK_T\sqrt{\pi T})^2}{1+n/2}\longrightarrow0.
\]
Both parity subsequences are summable, as is $\omega_n$ by
\eqref{eq:total-error-weights}. This proves
\eqref{eq:propagation-summability}, and in particular both weights
tend to zero. For a fixed index $m$, the initial-error term and the
finitely many terms with $k<m$ in the right side of
\eqref{eq:total-error-master} tend to zero as $N\to\infty$.
For $N>m$, the remaining terms satisfy
\[
\sum_{k=m}^{N-1}\omega_{N-1-k}\varepsilon_k
\le\left(\sum_{n=0}^{\infty}\omega_n\right)
\sup_{k\ge m}\varepsilon_k.
\]
First take the limit superior as $N\to\infty$, then let $m\to\infty$.
For bounded solver errors this proves \eqref{eq:persistent-error-limsup}.
If both solver errors vanish, then $\varepsilon_k\to0$, so the joint
upper bound and the cost upper bound tend to zero. Each individual
upper bound in \eqref{eq:main-value}--\eqref{eq:main-gradient} is
bounded by this joint upper bound, with the policy bound scaled by
$\tau$. Their convergence proves \eqref{eq:iteration-uniform-convergence}
and uniform convergence of the cost gap.
\end{proof}

\section{Numerical experiments}
\label{sec:experiments}

\noindent\hspace*{2em}We study how the choice of PDE solver and action sampler affects
value accuracy and the cost of a deployed controller. Two finite-horizon
execution problems test the methods with and without an exact solution.
We then examine multimodal routing under infinite-horizon discounting
and compare with martingale solvers on an unregularized HJB benchmark.

\subsection*{Experimental setup}

In the two finite-horizon execution comparisons, we pair DRDM
\citep{CFZ2025} or a physics-informed neural network (PINN) with
conditional flow matching \citep{LCBNL2023} or unadjusted Langevin
sampling. Each pair starts from the same reference policy. We evaluate
the fitted value function $\hat v$ on held-out time--state pairs and
the learned policy $\hat\pi$ on closed-loop trajectories with shared
initial states and Brownian increments. The baselines are continuous-time q-learning
\citep{JZ2023q}, Deep Splitting
\citep{BBCJN2021}, projected policy mirror descent \citep{KSSZ2025Mirror},
and PINN-SPI \citep{KCKK2026}.
Appendix~\ref{app:experiment-details} gives the training, evaluation,
and plotting settings for each experiment.

In each policy plot, a panel shows actions drawn from one policy at
the same test state as all other panels. Contours mark highest-density
regions of the projected comparison law.
Above each sampled policy we give the empirical sliced Wasserstein
distance ($\mathrm{SW}_2$) to the first panel. The comparison law is the optimal policy $\pi^*$ when known
and the prescribed reference policy $\mu$ otherwise.

When the optimal value $v_\tau$ is known, we compare the numerical value
$\hat v$ and its full state gradient with $v_\tau$ and $D_x v_\tau$ on
the same time--state pairs. The relative errors are
\begin{equation}
\begin{aligned}
 E_V&=\frac{\bigl(\sum_i|\hat v(t_i,x_i)-v_\tau(t_i,x_i)|^2\bigr)^{1/2}}
              {\bigl(\sum_i|v_\tau(t_i,x_i)|^2\bigr)^{1/2}},\\
 E_{\nabla}&=\frac{\bigl(\sum_i\|D_x\hat v(t_i,x_i)-D_x v_\tau(t_i,x_i)\|^2\bigr)^{1/2}}
                    {\bigl(\sum_i\|D_x v_\tau(t_i,x_i)\|^2\bigr)^{1/2}}.
\end{aligned}
\label{eq:experiment-gradient-metric}
\end{equation}
For both execution problems, we report the root-mean-square soft HJB
residual over $M$ held-out time--state pairs,
\begin{equation}
 \mathcal R_{\rm HJB}[\hat v]=\partial_t\hat v+\HH_\tau[\hat v],\qquad
 \mathrm{HJB\ RMS}=\left(\frac1M\sum_{i=1}^{M}
 |\mathcal R_{\rm HJB}[\hat v](t_i,x_i)|^2\right)^{1/2}.
\label{eq:experiment-hjb-metric}
\end{equation}
When the optimal policy $\pi^*$ is known, we average the reverse KL
of the deployed action law $\hat\pi$ over $N$ equally weighted test pairs,
\begin{equation}
 \mathcal K(\hat\pi)=\frac1N\sum_{i=1}^{N}
 \KL\bigl(\hat\pi(\cdot\mid t_i,x_i)\Vert
                 \pi^*(\cdot\mid t_i,x_i)\bigr).
\label{eq:experiment-policy-metric}
\end{equation}
For closed-loop evaluation, let $X_i$ be the Euler state at $t_i=ih$,
$Lh=T$, and $\hat\pi_i$ the action law conditional on the history through
$X_i$. For a Markov controller, $\hat\pi_i=\hat\pi(\cdot\mid t_i,X_i)$.
At step size $h$, the Euler approximations to the ordinary cost $J_0$
and the entropy-regularized cost $J_\tau$ are
\begin{equation}
\begin{aligned}
 J_{0,h}(\hat\pi)&=\E_h^{\hat\pi}\!\left[g(X_L)
     +h\sum_{i=0}^{L-1}c(t_i,X_i,U_i)\right],\\
 J_{\tau,h}(\hat\pi)&=J_{0,h}(\hat\pi)+\tau h\,\E_h^{\hat\pi}
     \!\left[\sum_{i=0}^{L-1}\KL(\hat\pi_i\Vert\mu_{t_i,X_i})\right].
\end{aligned}
\label{eq:experiment-cost-metric}
\end{equation}
Here $\E_h^{\hat\pi}$ averages over the initial state, sampled actions
and Brownian increments.

\subsection{Experiment 1: Execution with a known solution}
\label{subsec:experiment-exact}

\paragraph{Problem definition.}
Following the entropy-regularized control framework of \citet{WZZ2020},
consider $n=2500$ assets with state $x=(q,\alpha,y,z)\in\R^{4n}$
(inventory, signal, impact and liquidity) and action
$u=(u^M,u^P)\in\R^{2n}$ (market and passive execution speeds).
Starting from $X_t=x$, the controller chooses a feedback law
$\pi_s=\pi(\cdot\mid s,X_s)$ over $s\in[t,1]$. With reference law
$\mu_x$ and regularization weight $\tau=0.08/n$, the control problem is
\begin{equation}
 v_\tau(t,x)=\inf_\pi \E^\pi_{t,x}\!\left[
 \int_t^1\!\left(\int_{\R^{2n}}c(s,X_s,u)\,\pi_s(\dd u)
 +\tau\KL(\pi_s\Vert\mu_{X_s})\right)\dd s+g(X_1)\right].
\label{eq:execution-control}
\end{equation}
The infimum is over admissible feedback laws with finite expected cost
and KL penalty. The terminal cost is $g(x)=\tfrac12\langle q^2\rangle$,
where $\langle r\rangle=n^{-1}\sum_{i=1}^n r_i$.
Writing $\nu=u^M+u^P$ and $\bar\nu_s=\int_{\R^{2n}}\nu\,\pi_s(\dd u)$,
the state evolves according to
\begin{equation}
\begin{aligned}
 \dd q_s&=\bar\nu_s\,\dd s,&
 \dd\alpha_s&=-2\alpha_s\,\dd s+0.4\,\dd B_s^\alpha,\\
 \dd y_s&=(-4y_s+G\bar\nu_s)\,\dd s+0.2\,\dd B_s^y,&
 \dd z_s&=-1.5z_s\,\dd s
 +0.35\bigl(0.25\,\dd B_s^\alpha+\sqrt{1-0.25^2}\,\dd B_s^z\bigr).
\end{aligned}
\label{eq:experiment-dynamics}
\end{equation}
Here $B^\alpha,B^y,B^z$ are independent standard $n$-dimensional Brownian
motions, and $G=0.30I+0.08\Pi_{\rm sec}$, where $\Pi_{\rm sec}$ averages
over each of ten equal sectors. Conditional on $x$, the reference law
$\mu_x$ has independent five-dimensional Gaussian-mixture blocks in
$\xi=(u^M+u^P)/\sqrt2$ and a Gaussian law in $\zeta=(u^M-u^P)/\sqrt2$.
Let $\varrho_\mu(\cdot\mid x)$ be its Lebesgue density. To obtain a
benchmark with a known solution, we prescribe the optimal value $v_\tau$ and an action density
$p^*$ with the same block structure as $\mu_x$, and set
\begin{equation}
\begin{aligned}
 v_\tau(t,x)&=g(x)+(1-t)f(x),\\
 c(t,x,u)&=-\partial_t v_\tau(t,x)-\LL^u v_\tau(t,x)
 -\tau\log\frac{p^*(u\mid t,x)}{\varrho_\mu(u\mid x)}.
\end{aligned}
\label{eq:experiment-known-cost}
\end{equation}
Here $\LL^u$ is the spatial generator \eqref{eq:generator} with the
drift in \eqref{eq:experiment-dynamics} evaluated at the fixed speed
$\nu=u^M+u^P$. The soft HJB equation for \eqref{eq:execution-control} is
\begin{equation}
 0=\partial_t v_\tau(t,x)-\tau\log\int_{\R^{2n}}
 \exp\!\left[-\frac{c(t,x,u)+\LL^u v_\tau(t,x)}{\tau}\right]\mu_x(\dd u),
 \qquad v_\tau(1,x)=g(x).
\label{eq:execution-control-hjb}
\end{equation}
Substituting \eqref{eq:experiment-known-cost} into \eqref{eq:execution-control-hjb}
reduces the integral to $\exp(\partial_t v_\tau/\tau)$ because $p^*$ integrates
to one. The corresponding Gibbs feedback has density $p^*$.
Appendix~\ref{app:execution-model} specifies the initial-state law and
sector projection. Appendix~\ref{app:benchmark-costs} defines $f$, $p^*$
and $\varrho_\mu$. The solvers use the exact solution only for evaluation.

\paragraph{Evaluation criteria.}
We assess value and full-gradient accuracy using $E_V$ and
$E_{\nabla}$ on shared held-out time--state pairs, and equation
accuracy using the soft HJB residual RMS. Projected action samples
are compared with the exact marginal of $\pi^*$ at a shared test state.
The costs $J_{0,h}$ and $J_{\tau,h}$ measure performance along each
controller's own state trajectories. Appendix~\ref{app:experiment-evaluation}
gives the evaluation settings.

\begin{figure}[H]
\centering
\includegraphics[width=\linewidth]{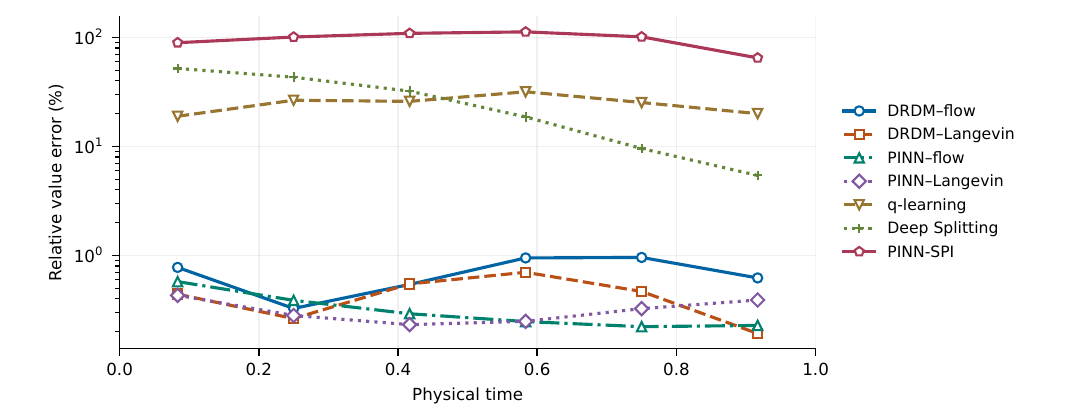}
\caption{Relative $L^2$ error in the value function by physical time.
Each marker gives $E_V$ in percent for the held-out pairs in one time
interval. The vertical axis is logarithmic.}
\label{fig:known-time-errors}
\end{figure}

\begin{figure}[H]
\centering
\includegraphics[width=\linewidth]{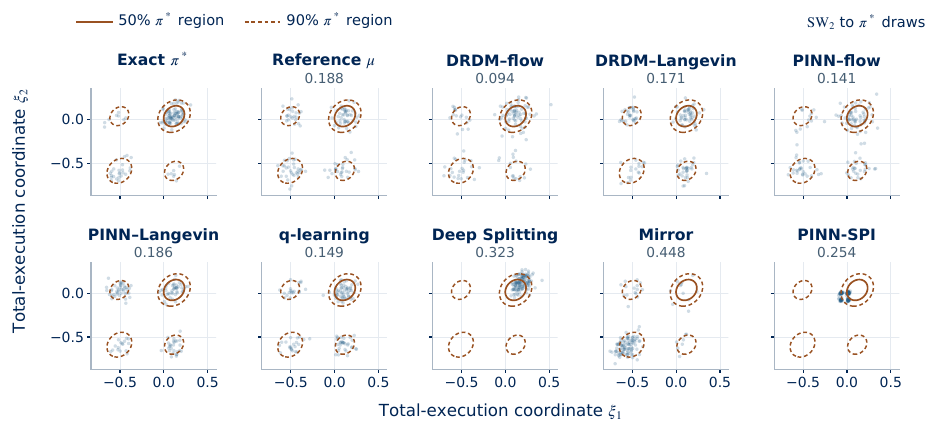}
\caption{Projected actions at a shared test pair.
Each panel shows draws in the first two total-execution coordinates
$(\xi_1,\xi_2)$ of the first block. Contours and distances use the
exact optimal marginal $\pi^*$ as the comparison.}
\label{fig:known-policy-scatter}
\end{figure}

\paragraph{Results and analysis.}
All four policy-iteration combinations approximate $v_\tau$ more accurately
than every baseline in Table~\ref{tab:exact-results}, with PINN--Langevin
attaining the smallest $E_V$. Their advantage over the plotted baselines
persists in every time interval (Figure~\ref{fig:known-time-errors}).
At the illustrated state, samples from all four combinations and q-learning
cover the clusters of the optimal action marginal, while Deep Splitting
concentrates near a single cluster (Figure~\ref{fig:known-policy-scatter}).
Both flow controllers achieve lower ordinary execution cost $J_{0,h}$ than
every baseline, with PINN--flow best overall. Both Langevin controllers incur
higher costs than Mirror (Table~\ref{tab:known-deployed-costs}).

\begin{table}[H]
\centering
\caption{Accuracy on the known-solution problem, measured by
\eqref{eq:experiment-gradient-metric} and \eqref{eq:experiment-hjb-metric}.
Mirror's value estimate uses conditional Monte Carlo.}
\label{tab:exact-results}
\begin{tabular}{lrrr}
\hline
Method & $E_V$ (\%) & $E_{\nabla}$ (\%) & $\mathrm{HJB\ RMS}$\\
\hline
DRDM--flow & $0.6823$ & $27.3894$ & $0.0049442$\\
DRDM--Langevin & $0.4460$ & $25.8983$ & $0.0040353$\\
PINN--flow & $0.4644$ & $\boldsymbol{19.0089}$ & $\boldsymbol{0.0027148}$\\
PINN--Langevin & $\boldsymbol{0.3568}$ & $19.3810$ & $0.0036499$\\
q-learning & $23.4263$ & $2619.6305$ & $7.0248504$\\
Deep Splitting & $44.0642$ & $99.8453$ & $0.0316631$\\
PINN-SPI & $97.3630$ & $321.9715$ & $0.5041311$\\
Mirror & $11.7040$ & N/A & N/A\\
\hline
\end{tabular}
\end{table}

\begin{table}[H]
\centering
\caption{Rollout costs \eqref{eq:experiment-cost-metric} on the known-solution
problem, reported as mean $\pm$ Monte Carlo standard error.}
\label{tab:known-deployed-costs}
\begingroup\small
\setlength{\tabcolsep}{5pt}
\begin{tabular}{lrr}
\hline
Method & $J_{0,h}$ & $J_{\tau,h}$\\
\hline
DRDM--flow & $0.338669\pm0.000405$ & $0.370928\pm0.000404$\\
DRDM--Langevin & $0.378911\pm0.000400$ & ---\\
PINN--flow & $\boldsymbol{0.337435\pm0.000402}$ & $\boldsymbol{0.369192\pm0.000407}$\\
PINN--Langevin & $0.378963\pm0.000397$ & ---\\
q-learning & $0.376266\pm0.000400$ & $0.376266\pm0.000400$\\
Deep Splitting & $0.361417\pm0.000401$ & $0.457001\pm0.000464$\\
PINN-SPI & $0.359976\pm0.000398$ & $0.670065\pm0.000488$\\
Mirror & $0.344230\pm0.000430$ & $0.392223\pm0.000386$\\
\hline
\end{tabular}
\par\smallskip
{\footnotesize For both execution problems, Appendix~\ref{app:langevin-missing-density}
reports estimated bounds on $J_\tau$ for the finite Langevin chains, whose
endpoint densities are unavailable for evaluating the KL penalty.
The Deep Splitting rollouts use Gibbs policies derived
from the fitted value gradients.\par}
\endgroup
\end{table}

\subsection{Experiment 2: Execution with an unknown solution}
\label{subsec:experiment-no-exact}

\paragraph{Problem definition.}
As in \eqref{eq:execution-control} of Experiment~1, the controller chooses
an admissible feedback law $\pi_s=\pi(\cdot\mid s,X_s)$ for the market and
passive execution speeds $u=(u^M,u^P)\in\R^{2n}$ of $n=2500$ assets.
The conditional mean total speed $\bar\nu_s=\int (u^M+u^P)\,\pi_s(\dd u)$
enters the inventory and impact dynamics in \eqref{eq:experiment-dynamics},
with $\Pi_{\rm sec}$ now averaging within eight sectors. Given
$X_t=x=(q,\alpha,y,z)\in\R^{4n}$, the value on $[t,1]$ is
\begin{equation}
 v_\tau(t,x)=\inf_\pi \E^\pi_{t,x}\!\left[
 \int_t^1\!\left(\int_{\R^{2n}}c_2(X_s,u)\,\pi_s(\dd u)
 +\tau\KL(\pi_s\Vert\mu_{X_s})\right)\dd s+g(X_1)\right].
\label{eq:coupled-execution-control}
\end{equation}
Here $g(x)=\tfrac12\langle q^2\rangle$, $\tau=0.08/n$, and $\mu_x$ is
the block-mixture reference law in \eqref{eq:app-half-reference}.
In the same action coordinates $(\xi,\zeta)$, neighboring total-execution
coordinates within each five-asset block are coupled by the last term of
the running cost
\begin{equation}
\begin{aligned}
 c_2(x,u)={}&-\langle q\alpha\rangle
 +\frac{0.225}{n}q_z^\top(I+0.30\Pi_{\rm sec})q_z\\
 &+\left\langle(y-q)\sqrt2\xi+\lambda(z)(\xi^c)^2+0.15\zeta^2\right\rangle
 +\frac{0.035}{n}\sum_b(\xi_b^c)^\top L_5\xi_b^c.
\end{aligned}
\label{eq:experiment-c-cost}
\end{equation}
Here $\widetilde z=\operatorname{clip}(z,-2,2)$,
$q_z=q\odot e^{\widetilde z/2}$, $\lambda(z)=0.12+0.03\cosh(\widetilde z)$,
and $\xi^c=\xi-h(x)/\sqrt2$ with
$h(x)=\operatorname{clip}(-0.80q+0.15\alpha-0.08y,-1.15,1.15)$.
The matrix $L_5$ is the cyclic five-site graph Laplacian and $b$ indexes
the blocks. With $\LL^u$ the generator for the eight-sector dynamics,
the soft HJB equation for \eqref{eq:coupled-execution-control} is
\begin{equation}
 0=\partial_t v_\tau(t,x)-\tau\log\int_{\R^{2n}}
 \exp\!\left[-\frac{c_2(x,u)+\LL^u v_\tau(t,x)}{\tau}\right]\mu_x(\dd u),
 \qquad v_\tau(1,x)=g(x).
\label{eq:coupled-execution-control-hjb}
\end{equation}
The optimal feedback is the Gibbs law $q_{v_\tau}$ in \eqref{eq:improvement-map},
with $c_2$ as its running cost. Its density depends on the unknown value gradient.

\paragraph{Evaluation criteria.}
We assess equation accuracy with the nonlinear HJB residual RMS in
\eqref{eq:experiment-hjb-metric} and controller performance with the
achieved costs in \eqref{eq:experiment-cost-metric}. We also compare
projected action samples with the reference law $\mu$ at a shared test
state. Appendix~\ref{app:experiment-two-evaluation} gives the evaluation settings.

\begin{figure}[H]
\centering
\includegraphics[width=\linewidth]{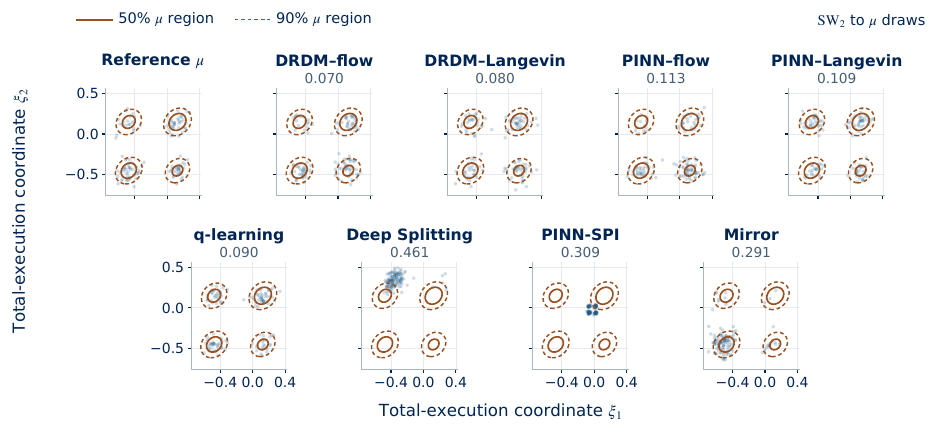}
\caption{Projected actions in Experiment~2, in the first two
total-execution coordinates $(\xi_1,\xi_2)$ of the first block at the
first shared test state.
Contours and distances use the reference distribution $\mu$.}
\label{fig:problem-c-policy-scatter}
\end{figure}

\paragraph{Results and analysis.}
Flow yields a lower ordinary execution cost $J_{0,h}$ than Langevin with
either PDE solver when the running cost couples neighboring assets
(Table~\ref{tab:deployed-finite-costs}). Switching PDE solvers with the
sampler fixed has a much smaller effect on $J_{0,h}$, and the two flow
controllers give nearly identical estimates of $J_{\tau,h}$. DRDM--flow
attains the smallest soft HJB residual on the test states
(Table~\ref{tab:no-exact-quality}). Despite their different closed-loop
costs, all four combinations retain the clusters of the reference policy
$\mu_x$ in the displayed action coordinates at the illustrated state
(Figure~\ref{fig:problem-c-policy-scatter}).
Interpreting $\mu$ as an expert prior, we seek to avoid excessive deviation
from it \citep{GJH2019DefaultPolicy,TGH2022BehaviorPriors};
Figure~\ref{fig:problem-c-policy-scatter} shows that DRDM--flow and
DRDM--Langevin attain the two smallest empirical $\mathrm{SW}_2$ distances
to $\mu$ among the learned policies.

\begin{table}[H]
\centering
\caption{Soft HJB residual on the held-out states of Experiment~2,
measured by \eqref{eq:experiment-hjb-metric}.
Mirror has no differentiable scalar critic.}
\label{tab:no-exact-quality}
\small
\begin{tabular}{lr}
\hline
Method & $\mathrm{HJB\ RMS}$\\
\hline
DRDM--flow & $\boldsymbol{0.0020702}$\\
DRDM--Langevin & $0.0109895$\\
PINN--flow & $0.0022276$\\
PINN--Langevin & $0.0078722$\\
q-learning & $6.9359038$\\
Deep Splitting & $0.4631421$\\
PINN-SPI & $0.4636509$\\
Mirror & N/A\\
\hline
\end{tabular}
\end{table}

\begin{table}[H]
\centering
\caption{Rollout costs \eqref{eq:experiment-cost-metric} in Experiment~2,
reported as mean $\pm$ Monte Carlo standard error.}
\label{tab:deployed-finite-costs}
\begingroup\small
\setlength{\tabcolsep}{5pt}
\begin{tabular}{lrr}
\hline
Method & $J_{0,h}$ & $J_{\tau,h}$\\
\hline
DRDM--flow & $0.315189\pm0.000531$ & $0.327255\pm0.000531$\\
DRDM--Langevin & $0.337573\pm0.000550$ & ---\\
PINN--flow & $\boldsymbol{0.314673\pm0.000502}$ & $\boldsymbol{0.327186\pm0.000502}$\\
PINN--Langevin & $0.337419\pm0.000510$ & ---\\
q-learning & $0.331638\pm0.000520$ & $0.331638\pm0.000520$\\
Deep Splitting & $0.379246\pm0.000577$ & $0.530242\pm0.000958$\\
PINN-SPI & $0.370841\pm0.000587$ & $0.680930\pm0.000661$\\
Mirror & $0.347399\pm0.000517$ & $0.392691\pm0.000495$\\
\hline
\end{tabular}
\endgroup
\end{table}

\FloatBarrier

\subsection{Experiment 3: Execution over an infinite horizon}
\label{subsec:experiment-stationary-routing}

\paragraph{Problem definition.}
The state $x\in\R^{1024}$ follows \eqref{eq:experiment-dynamics} with
$n=256$ assets in two sectors of $128$ and action $u\in\R^{512}$.
Starting from $X_0=x$, we minimize the discounted cost over admissible
stationary feedback laws $\pi_x=\pi(\cdot\mid x)$. With discount rate
$\beta=0.5$ and regularization weight $\tau=0.08/n$, the value is
\begin{equation}
 v_\tau(x)=\inf_\pi \E^\pi_x\!\left[
 \int_0^\infty e^{-\beta s}\!\left(\int_{\R^{512}} c_\infty(X_s,u)\,
 \pi_{X_s}(\dd u)+\tau\KL(\pi_{X_s}\Vert\mu_{X_s})\right)\dd s\right].
\label{eq:stationary-routing-control}
\end{equation}
The Gaussian reference $\mu_x=\mathcal N(0,0.20^2 I_{512})$ has density
$\varrho_\mu$. We prescribe the optimal value $v_\tau(x)=g(x)+f(x)$
from \eqref{eq:app-manufactured-value} and the optimal density $p^*(u\mid x)$.
Under $p^*$, $\xi$ follows a state-dependent mixture of twelve Gaussians,
and the conditional mean of $\zeta$ depends nonlinearly on dense projections
of $\xi$. The running cost is
\[
 c_\infty(x,u)=\beta v_\tau(x)-\LL^u v_\tau(x)
 -\tau\log\frac{p^*(u\mid x)}{\varrho_\mu(u\mid x)}.
\]
With the two-sector generator $\LL^u$, the stationary soft HJB equation is
\begin{equation}
 \beta v_\tau(x)=-\tau\log\int_{\R^{512}}
 \exp\!\left[-\frac{c_\infty(x,u)+\LL^u v_\tau(x)}{\tau}\right]\mu_x(\dd u).
\label{eq:experiment-stationary-main}
\end{equation}
On high-dimensional targets with well-separated modes, a finite Langevin
chain can stay near one mode and undersample the others \citep{CFF2025Langevin}.
For this routing problem we use a projected Langevin update based on
Gaussian moment matching \citep{RSGE2013}. At each state we represent the
action law by a mixture with one diagonal Gaussian per identified mode,
so actions can be sampled directly from the fitted mixture.
We call the resulting implementation of continuous policy and value
iteration \citep{FW2025} \emph{projected CPVI}. Its Langevin update
evolves each component's mean and coordinatewise variances at fixed
mixture weights. In DRDM--projected Langevin and PINN--projected Langevin,
the mixture is fitted to the endpoints of finite Langevin chains. \emph{Projected mirror descent}
fits the same mixture family to the policy obtained by the mirror update
\citep{KSSZ2025Mirror}.
The four solver--sampler combinations and projected CPVI share an initial
mixture fitted to the public Hamiltonian and a Monte Carlo value fit.
Each combination evaluates its own policy and improves it using either
a conditional affine-coupling flow or projected Langevin. Projected CPVI
updates its value and policy simultaneously. Score--Hamiltonian PI
\citep{FW2026Amortized} couples the two updates starting from its own initial
policy. The remaining baselines are Gaussian q-learning and Deep
Splitting. Appendix~\ref{app:stationary-routing} details the solvers
and separately evaluates the finite chain.

\begin{figure}[!htb]
\centering
\includegraphics[width=\linewidth]{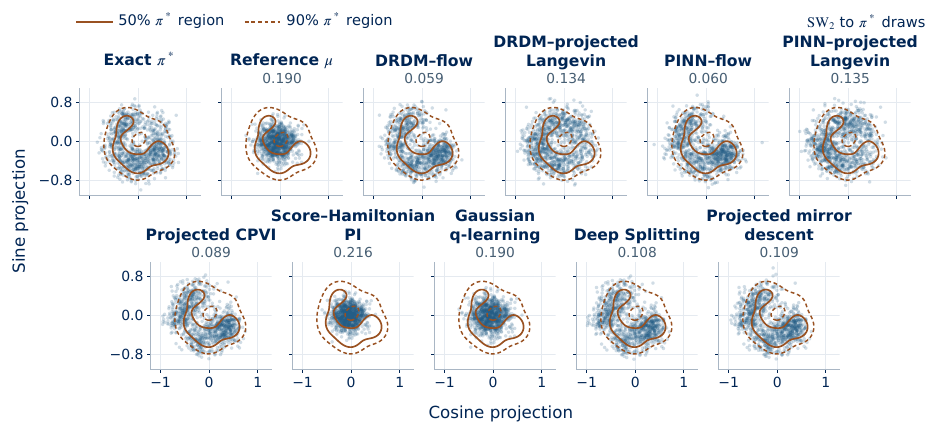}
\caption{Projected total-execution samples at a shared held-out state.
The first two show the exact optimal law $\pi^*$ and the Gaussian
reference $\mu$, followed by the deployed policies. All panels use
the same centered projection.
Contours and distances use the projected $\pi^*$.
Appendix~\ref{app:routing-evaluation} defines the projection and plotting settings.}
\label{fig:stationary-action-scatter}
\end{figure}

\paragraph{Evaluation criteria.}
We report relative value and full-gradient errors, reverse KL to
$\pi^*$ at held-out states, and discounted excess cost
$\Delta_T$ from \eqref{eq:routing-closed-loop-excess}. The last
quantity integrates the discounted policy KL along the controller's
own trajectories up to the evaluation horizon, weighted by $\tau$.
Results summarize variation across independent training runs, with
paired trajectories used for the cost comparison.
Appendix~\ref{app:routing-evaluation} gives the evaluation settings.

\paragraph{Results and analysis.}
DRDM--flow attains the lowest mean policy KL and the lowest mean discounted
excess cost among the routing policies in Table~\ref{tab:stationary-routing-comparison}.
Both measures favor flow over projected Langevin with either PDE solver and
DRDM over PINN with either sampler. DRDM--projected Langevin gives the most
accurate value estimate. Both flow policies recover the curved shape of the
projected optimal action law at the illustrated state
(Figure~\ref{fig:stationary-action-scatter}), although both select some
routes too frequently across the test states (Table~\ref{tab:stationary-mode-bias}).

\begin{table}[!htb]
\centering
\caption{Stationary routing metrics from \eqref{eq:experiment-gradient-metric},
\eqref{eq:experiment-policy-metric} and \eqref{eq:routing-closed-loop-excess},
reported as mean $\pm$ sample standard deviation over independent training seeds.
The evaluation horizon for $\Delta_T$ is specified in Appendix~\ref{app:routing-evaluation}.}
\label{tab:stationary-routing-comparison}
\vspace{.4em}
\small
\setlength{\tabcolsep}{4pt}
\begin{tabular}{lrrrr}
\hline
Method & $E_V$ (\%) & $E_{\nabla}$ (\%) & $\mathcal K(\hat\pi)$ & $10^3\Delta_T(\hat\pi)$\\
\hline
DRDM--flow & $0.576\pm0.121$ & $8.13\pm0.74$ & $\mathbf{5.95}\pm1.04$ & $\mathbf{3.824}\pm0.683$\\
DRDM--projected Langevin & $\mathbf{0.202}\pm0.016$ & $7.08\pm0.24$ & $7.34\pm0.68$ & $4.576\pm0.400$\\
PINN--flow & $0.927\pm0.123$ & $9.42\pm0.66$ & $7.05\pm0.95$ & $4.545\pm0.620$\\
PINN--projected Langevin & $0.467\pm0.096$ & $8.85\pm0.49$ & $8.52\pm0.67$ & $5.384\pm0.415$\\
Projected CPVI & $0.746\pm0.053$ & $10.65\pm0.45$ & $10.33\pm0.38$ & $6.598\pm0.224$\\
Score--Hamiltonian PI & $3.045\pm0.161$ & $35.95\pm0.87$ & $15.88\pm0.06$ & $10.061\pm0.038$\\
Gaussian q-learning & $8.959\pm0.058$ & $23.68\pm0.30$ & $28.12\pm0.03$ & $17.575\pm0.014$\\
Deep Splitting & $1.292\pm0.017$ & $13.18\pm0.09$ & $13.02\pm0.68$ & $8.408\pm0.526$\\
Projected mirror descent & $0.392\pm0.121$ & $\mathbf{6.51}\pm0.90$ & $11.10\pm0.07$ & $7.184\pm0.039$\\
\hline
\end{tabular}
\end{table}

\FloatBarrier

\subsection{Experiment 4: Unregularized HJB equation}
\label{subsec:experiment-hjb3a}

\paragraph{Problem definition.}
We adapt HJB-3a from \citet[Section~4.4]{CFZZ2024DFMartNet} by dividing
its quadratic action cost by the state dimension. With $T=1$ and
$d=m\in\{1000,10{,}000\}$, the dynamics and costs are
\begin{equation}
\begin{aligned}
 \dd X_t &= (\mathbf{1}_d+2U_t)\,\dd t+\sqrt{0.08}\,\dd B_t,
 & c(t,x,u)&=\frac{25}{d}\|u\|^2,\\
 g(x)&=\frac{1}{d}\sum_{i=1}^d\phi(x_i-1),
 & \phi(r)&=-\cos r+\sin\bigl((0.1\pi+r^2)^{-1}\bigr).
\end{aligned}
\label{eq:experiment-hjb3a-problem}
\end{equation}
Here $B$ is standard $d$-dimensional Brownian motion, $\mathbf{1}_d$
is the vector of ones and $X_0$ is uniform on $[-1,1]^d$, independent
of $B$. For adapted controls $U_s\in\R^d$ of finite expected energy,
the unregularized optimal value is
\[
 v(t,x)=\inf_U \E_{t,x}\!\left[
       \int_t^1\frac{25}{d}\|U_s\|^2\,\dd s+g(X_1)\right].
\]
Minimizing the quadratic Hamiltonian gives the classical HJB equation
\begin{equation}
 0=\partial_t v+\boldsymbol{1}_d^\top D_x v+0.04\Delta_x v
      -\frac{d}{25}\|D_x v\|^2,\qquad v(1,x)=g(x).
\label{eq:hjb3a-classical}
\end{equation}
Separation in the coordinates $x_i-t$ gives the Cole--Hopf reference
\begin{equation}
 v(t,x)=-\frac{1}{d}\sum_{i=1}^d
 \log\E_Z\!\left[\exp\!\left\{-\phi\!\left(x_i-t+\sqrt{0.08(1-t)}\,Z\right)\right\}\right],
 \qquad Z\sim\mathcal N(0,1).
\label{eq:experiment-hjb3a-reference}
\end{equation}
We pair each PDE solver (DRDM or PINN) with each policy update
(flow or Langevin). Each controller uses the conditional mean action
$U_t=\int_{\R^d}u\,\hat\pi(\dd u\mid t,X_t)$ in both the dynamics and
the quadratic running cost. For Langevin, the Gaussian target gives an
affine score, so we compute the finite-chain mean exactly. DF-MartNet
uses weak martingale and finite-difference control updates
\citep{CFZZ2026DFMartNet}, and SOC-MartNet uses gradient-based
Hamiltonian updates \citep{CFZ2025SOCMartNet}. All methods share
scalar network parameters across coordinates and train using only
the problem data. Appendix~\ref{app:hjb3a-solver-pairs} gives the
critic residual and finite-chain update.

\paragraph{Evaluation criteria.}
We compute $E_V$ in \eqref{eq:experiment-gradient-metric} against the
Cole--Hopf reference $v$ on shared held-out time--state pairs and
evaluate $J_{0,h}$ on paired closed-loop trajectories. Training settings
are fixed across dimensions, and the number of value-network updates
is matched across methods; the control updates and work per update
remain method dependent. Results summarize variation across independent
training runs. We compare the two flow controllers using a paired
$t$-interval for their differences in mean cost.
Appendix~\ref{app:hjb3a-evaluation} gives the evaluation settings.

\begin{table}[!htb]
\centering
\caption{HJB-3a metrics from \eqref{eq:experiment-gradient-metric} and
\eqref{eq:experiment-cost-metric}, reported as mean $\pm$ sample standard
deviation over independent training seeds at each dimension.}
\label{tab:experiment-hjb3a}
\begingroup\footnotesize
\setlength{\tabcolsep}{2.4pt}
\begin{tabular}{@{}lrrrr@{}}
\hline
 & \multicolumn{2}{c}{$d=m=1000$} & \multicolumn{2}{c}{$d=m=10{,}000$}\\
\cline{2-3}\cline{4-5}
Method & $E_V$ (\%) & $J_{0,h}$ & $E_V$ (\%) & $J_{0,h}$\\
\hline
DRDM--flow & $0.147\pm0.057$ & $-0.2229799\pm0.0000023$ & $0.143\pm0.057$ & $-0.2247566\pm0.0000004$\\
DRDM--Langevin & $0.163\pm0.059$ & $-0.2229785\pm0.0000006$ & $0.137\pm0.084$ & $-0.2247515\pm0.0000017$\\
PINN--flow & $\mathbf{0.095}\pm0.027$ & $\mathbf{-0.2229802}\pm0.0000014$ & $\mathbf{0.057}\pm0.021$ & $\mathbf{-0.2247569}\pm0.0000002$\\
PINN--Langevin & $0.119\pm0.025$ & $-0.2229792\pm0.0000003$ & $0.083\pm0.015$ & $-0.2247532\pm0.0000002$\\
DF-MartNet & $1.386\pm0.241$ & $-0.2224078\pm0.0001138$ & $1.285\pm0.186$ & $-0.2241795\pm0.0000339$\\
SOC-MartNet & $5.039\pm2.947$ & $-0.2222019\pm0.0005186$ & $3.519\pm1.067$ & $-0.2239153\pm0.0006676$\\
\hline
\end{tabular}
\endgroup
\end{table}

\paragraph{Results and analysis.}
All four policy-iteration combinations achieve lower mean $E_V$ and $J_{0,h}$
than DF-MartNet and SOC-MartNet at both dimensions in Table~\ref{tab:experiment-hjb3a}.
PINN--flow gives the most accurate value estimates at the lower dimension
and retains this lead at the higher dimension with its training settings
fixed. The paired comparison gives no clear difference in closed-loop cost
between the two flow controllers at either dimension. DF-MartNet
overestimates the reference value most at early test times, while DRDM--flow
stays close to it across the tested times
(Figures~\ref{fig:hjb3a-values} and~\ref{fig:hjb3a-value-errors}).

\begingroup
\setlength{\intextsep}{4pt}
\begin{figure}[H]
\centering
\includegraphics[width=\linewidth]{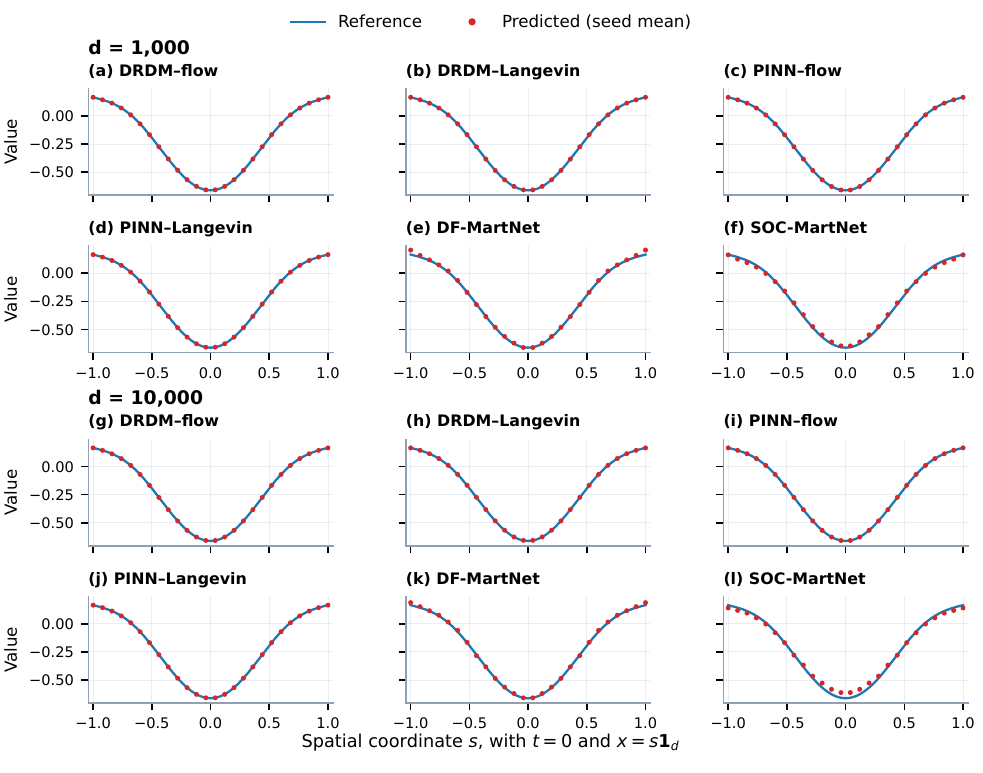}
\caption{Spatial profiles $s\mapsto v(0,s\mathbf1_d)$ for HJB-3a,
following the slice used by \citet[Figure~5]{CFZZ2024DFMartNet}.
Panels (a)--(f) use $d=m=1000$, and panels (g)--(l) use $d=m=10{,}000$.
Blue curves show the Cole--Hopf reference
\eqref{eq:experiment-hjb3a-reference}; red markers show the pointwise
mean prediction over three training seeds. All panels use the same
spatial interval and vertical scale.}
\label{fig:hjb3a-values}
\end{figure}

\begin{figure}[H]
\centering
\includegraphics[width=\linewidth]{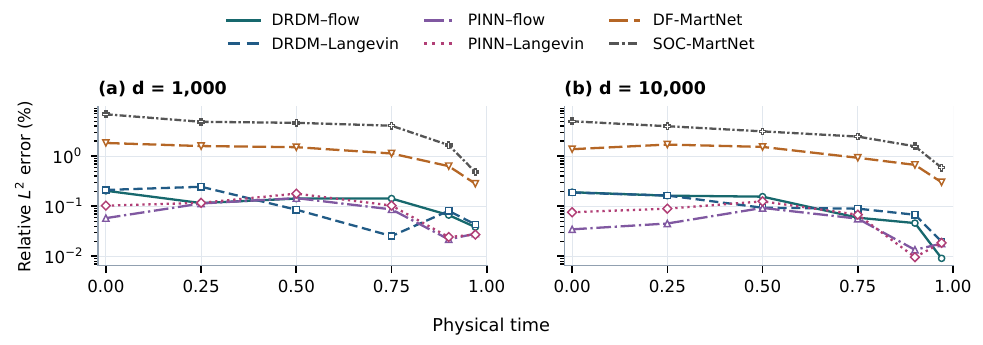}
\caption{Relative $L^2$ value error by physical time for all six methods
in HJB-3a. At each time, compute $E_V$ on the eight held-out states
separately for each training seed, then average over the three seeds.
Panels (a) and (b) use $d=m=1000$ and $d=m=10{,}000$.
Errors are shown in percent on a common logarithmic vertical scale.}
\label{fig:hjb3a-value-errors}
\end{figure}
\endgroup

\FloatBarrier

\section{Conclusion}\label{sec:conclusion}

For finite-horizon entropy-regularized stochastic control, we quantified
how errors in two numerical solves affect the policy ultimately deployed.
The residual of the frozen-policy PDE controls the value gradient used to
form the next Gibbs distribution. The KL divergence from the returned policy
to that distribution measures the error in the improvement step. Their
weighted history bounds the errors in the policy value, its gradient, and
the final policy, as well as the excess regularized cost. On the
finite-horizon benchmark with a known solution, held-out relative
value errors ranged from \(0.36\%\)
for PINN--Langevin to \(0.68\%\) for DRDM--flow. PINN--flow gave the
smallest relative gradient error (\(19.0\%\)) and the smallest nonlinear
HJB residual (Table~\ref{tab:exact-results}).

Two questions follow from the analysis. First, we would like to replace
the uniform norms in the iteration bound by \(L^2\) norms under a
specified distribution of times and states. An \(L^2\) estimate for the
final value and its gradient would connect the theorem to mean-square
error estimates for the PDE solver. This requires tracking those errors
along the trajectories of successive policies. Second, when the
running cost is unknown, it must be estimated from observations. The
estimate enters both the frozen-policy equation and the energy defining
the Gibbs target. Quantifying these two effects would extend the analysis
to policies computed from observed costs.

\appendix

\section{Experimental details}
\label{app:experiment-details}

This appendix first describes the common numerical methods, then gives
one subsection for each of Experiments~1--4 in the order used in
Section~\ref{sec:experiments}. Each experiment is organized into problem
setting, training protocol, and evaluation and reporting.
In the flow construction, $t$ is physical time and $s$ is generative
time; physical trajectory integrals elsewhere retain their stated
integration variables. The index $k$ counts outer policy iterations.
The four solver combinations in Experiments~1 and~2
each perform two policy improvements. Their tables report
the second policy and the value approximation fitted with its samples.
The q-learning, Deep Splitting, projected mirror descent, and PINN-SPI
rows use their respective training protocols on the same problems.
The reported errors reflect time discretization, finite evaluation samples,
policy approximation and variation across training seeds. Stabilization of
the validation residual alone does not establish convergence or optimality.

\subsection{Common numerical methods}\label{app:solver-implementations}

\subsubsection{Value updates by random differences}

The Deep Random Difference Method (DRDM) implements the PDE solve in
Section~\ref{sec:numerical-method}. It approximates the differential
operator using network values at randomly displaced points and trains
the network to make the mean residual small. We first describe the
evaluation of a frozen policy, then explain the compensated equation
used in Experiments~1 and~2.

\paragraph{The equation and the unknown network.}
Fix the deployed policy $\pi_k$. Using the notation of
\eqref{eq:pi-evaluation} and \eqref{eq:residual-definitions}, its value
solves
\[
\mathcal R^{\pi_k}[v^{\pi_k}]
=\partial_t v^{\pi_k}+\LL^{\pi_k}v^{\pi_k}+\ell^{\pi_k}=0,
\qquad v^{\pi_k}(T,x)=g(x).
\]
Here $\LL^{\pi_k}$ accounts for the averaged drift and state diffusion,
and $\ell^{\pi_k}$ is the running cost including relative entropy.
We approximate this value by $\hat v_\psi$, with trainable parameters
$\psi$. For smooth terminal data, the parameterization
\begin{equation}\label{eq:terminal-network}
\hat v_\psi(t,x)=g(x)+(T-t)N_\psi(t,x)
\end{equation}
enforces the terminal condition. During this solve, the policy and its
density remain fixed; only $\psi$ is updated. The density-based
implementation takes actions in $\R^m$, with densities $p_k$ and
$\varrho_\mu$ for $\pi_k$ and $\mu_{t,x}$, respectively, and assumes
$\pi_k(\cdot\mid t,x)\ll\mu_{t,x}$.

\paragraph{Step 1: sample starting points.}
Draw a minibatch $(t_n,x_n)$, $n=1,\ldots,M_B$, from the training
distribution $\nu_k$. Depending on the experiment, these points come
from current-policy trajectories or a prescribed state pool. Uniform
sampling from a fixed pool means that $\nu_k$ assigns equal mass to
its entries. At each point choose a physical difference step $h>0$
with $t_n+h\le T$. This step is distinct from the optimizer learning
rate and from the generative-time step used by the flow policy.

\paragraph{Step 2: construct random neighboring points.}
At a starting point $(t,x)$, draw two independent groups, indexed by
$i=1,2$, each containing $M_E$ independent pairs
\[
U_{i,j}\sim\pi_k(\cdot\mid t,x),\qquad
\xi_{i,j}\sim\mathcal N(0,I_r),\qquad j=1,\ldots,M_E.
\]
The actions and Gaussian vectors are independent. For each pair set
\begin{equation}\label{eq:drdm-successors}
Y_{i,j}^{\pm}
=x+h\,b(t,x,U_{i,j})
\pm\sqrt h\,\sigma(t,x)\xi_{i,j}.
\end{equation}
The two points share the same action and use opposite noise increments.
They provide local samples for the differential operator; they are
constructed afresh from the starting point, rather than requiring
another complete state trajectory.

\paragraph{Step 3: compute two residual estimates.}
Evaluate the same network at $(t,x)$ and at all the points
$(t+h,Y_{i,j}^{\pm})$. For each group form
\begin{equation}\label{eq:drdm-block}
\begin{aligned}
R_i(\psi;t,x)=\frac1{M_E}\sum_{j=1}^{M_E}\Bigg[
&\frac{\tfrac12\bigl[
\hat v_\psi(t+h,Y_{i,j}^{+})+
\hat v_\psi(t+h,Y_{i,j}^{-})\bigr]-\hat v_\psi(t,x)}{h}\\
&+c(t,x,U_{i,j})
+\tau\log\frac{p_k(U_{i,j}\mid t,x)}
{\varrho_\mu(U_{i,j}\mid t,x)}\Bigg].
\end{aligned}
\end{equation}
The first term is the average next-step value minus the current value,
divided by $h$. Its expectation approximates
$(\partial_t+\LL^{\pi_k})\hat v_\psi$; the remaining terms estimate
$\ell^{\pi_k}$. Thus $R_i$ estimates the policy-evaluation residual.
The positive and negative noise increments cancel the first-order
random term in the value difference.

\paragraph{Step 4: update the value network.}
Use the minibatch loss
\[
\widehat{\mathcal J}_{E,k}(\psi)
=\frac1{M_B}\sum_{n=1}^{M_B}
R_1(\psi;t_n,x_n)R_2(\psi;t_n,x_n).
\]
Backpropagate through both factors and update $\psi$, keeping the
sampled points and frozen policy fixed for that update. Repeat with
new minibatches and random samples. No exact value labels are needed.
This random-difference construction, motivated by \citet{CFZ2025},
uses network values and parameter gradients without explicitly
evaluating $\partial_t\hat v_\psi$, $D_x\hat v_\psi$, or
$D_x^2\hat v_\psi$ in the critic loss. After fitting, write
$\hat v_k=\hat v_{\psi_k}$. The subsequent Gibbs improvement uses
$D_x\hat v_k$, as in Section~\ref{sec:numerical-method}.

\paragraph{Why use two independent groups?}
All expectations conditional on $(t,x)$ below average only the sampled
actions and Gaussian vectors, with the network and policy fixed.
Write the common mean of the two residual estimates as
\[
\mathcal R_h^{\pi_k}[\hat v_\psi](t,x)
=\E[R_i(\psi;t,x)\mid t,x].
\]
The PDE requires this mean to be small. Squaring a single estimate
would also penalize its sampling variance: whenever the second moment
is finite,
\[
\E[R_1^2\mid t,x]
=\bigl(\mathcal R_h^{\pi_k}[\hat v_\psi](t,x)\bigr)^2
+\operatorname{Var}(R_1\mid t,x).
\]
The variance generally depends on $\psi$. Two independent groups at
the same starting point instead give
\begin{equation}\label{eq:drdm-product}
\E[R_1R_2\mid t,x]
=\E[R_1\mid t,x]\E[R_2\mid t,x]
=\bigl(\mathcal R_h^{\pi_k}[\hat v_\psi](t,x)\bigr)^2.
\end{equation}
Consequently, the population objective corresponding to the minibatch
loss is
\begin{equation}\label{eq:drdm-objective}
\begin{aligned}
\mathcal J_{E,k}(\psi)
&=\int \E[R_1R_2\mid t,x]\,\nu_k(\dd t,\dd x)\\
&=\int\bigl(\mathcal R_h^{\pi_k}[\hat v_\psi](t,x)\bigr)^2
\nu_k(\dd t,\dd x).
\end{aligned}
\end{equation}
Individual products can be negative even though this population
objective is nonnegative. The two sampling devices have different
roles: opposite noise increments cancel odd noise terms within a
pair, whereas independent groups remove the extra variance term from
the expected loss. Increasing $M_E$ reduces sampling variance; reducing
$h$ reduces the finite-step bias under the conditions stated below.

\paragraph{Computed densities, action caches, and stationary problems.}
Equation~\eqref{eq:drdm-block} defines the residual using the true
policy density. When the implementation supplies the computed log
density $\widetilde{\log p}_k$, it replaces $\log p_k$ in that equation.
The resulting loss then concerns the residual with this computed cost;
the density approximation remains part of the evaluation error.
For a fixed action cache, the groups draw indices independently with
replacement and use independent Gaussian increments. Conditional on
the cache, their product targets the square of the cache-averaged
residual. The finite cache is an additional approximation to the
policy expectation; the log-density term still uses the policy that
generated the cache. In the stationary discounted problems, omit the
time argument and add $-\beta\hat v_\psi(x)$ to each residual estimate.

\paragraph{Compensated updates versus true policy evaluation.}
The runs for Experiments~1 and~2 use the protocol in
Appendix~\ref{app:finite-training}. They evaluate
the initial reference policy with its actual running cost. After an
improvement, however, they use a compensated cost while fitting the
next value. With the preceding guide held fixed, the cost used in
this fit is
\begin{equation}\label{eq:drdm-compensated-cost}
C_{\hat v_k}(t,x,u)
=\HH_\tau[\hat v_k](t,x)-\LL^u\hat v_k(t,x).
\end{equation}
Both the guide $\hat v_k$ and the deployed law $\pi_{k+1}$ are frozen
during this fit. To implement it, sample actions from $\pi_{k+1}$ in
\eqref{eq:drdm-successors}--\eqref{eq:drdm-block} and replace the entire
running-cost and log-density contribution by $C_{\hat v_k}$.
The same random differences and two-group product are used, but the
limiting equation for a candidate $w$ is now
\begin{equation}\label{eq:drdm-compensated-equation}
\partial_t w
+\LL^{\pi_{k+1}}(w-\hat v_k)
+\HH_\tau[\hat v_k]=0,
\qquad w(T,x)=g(x).
\end{equation}
This follows by averaging \eqref{eq:drdm-compensated-cost} under
$\pi_{k+1}$.

The Gibbs identity \eqref{eq:gibbs-identity}, with
$\hat q_k=q_{\hat v_k}$, gives the precise difference from the
true running cost:
\begin{equation}\label{eq:drdm-compensation-gap}
\ell^{\pi_{k+1}}(t,x)
-\int_{\R^m}C_{\hat v_k}(t,x,u)\pi_{k+1}(\dd u\mid t,x)
=\tau\KL(\pi_{k+1}\Vert\hat q_k)(t,x).
\end{equation}
Here $\hat q_k$ is the exact Gibbs law of the fixed guide, and the
identity assumes its normalizing constant is finite and positive and
the displayed costs and KL divergence are finite. Therefore, for any
smooth candidate $w$,
\[
\mathcal R^{\pi_{k+1}}[w]
=\partial_t w+\LL^{\pi_{k+1}}(w-\hat v_k)
+\HH_\tau[\hat v_k]
+\tau\KL(\pi_{k+1}\Vert\hat q_k).
\]
The compensated equation coincides with true policy evaluation when
$\pi_{k+1}=\hat q_k$. Otherwise, even an exact solution of
\eqref{eq:drdm-compensated-equation} retains the displayed KL term in
its true policy-evaluation residual. We therefore distinguish the
value approximation obtained from this equation from the actual cost
of the deployed policy.

\paragraph{Relation to the theoretical evaluation error.}
For clarity, the two solver errors are
\[
\begin{aligned}
E_{\mathrm{PDE}}^k
&=\|\mathcal R^{\pi_k}[\hat v_k]\|_\infty,\\
E_{\mathrm{policy}}^k
&=\sup_{0\le t<T,\,x}\KL(\pi_{k+1}\Vert\hat q_k)(t,x).
\end{aligned}
\]
The norm and supremum are over $[0,T)\times\R^d$.
Thus $E_{\mathrm{PDE}}^k$ uses the true policy-evaluation residual,
while $E_{\mathrm{policy}}^k$ measures the error of the next policy
relative to $\hat q_k$. For a compensated fit at $k\ge1$, the
identity above gives
\begin{equation}\label{eq:compensated-evaluation-error}
\begin{aligned}
E_{\mathrm{PDE}}^k
\le{}&\bigl\|\partial_t\hat v_k
+\LL^{\pi_k}(\hat v_k-\hat v_{k-1})
+\HH_\tau[\hat v_{k-1}]\bigr\|_\infty\\
&+\tau E_{\mathrm{policy}}^{k-1}.
\end{aligned}
\end{equation}
The correction uses the preceding policy-improvement error. Even if
the continuous compensated equation is solved exactly, the true
evaluation error is $E_{\mathrm{PDE}}^k=\tau E_{\mathrm{policy}}^{k-1}$.
The norm here concerns the continuous equation; a sampled training
loss or validation RMS is not that uniform norm. Experiments~1 and~2
also report the nonlinear HJB residual and the deployed closed-loop
cost separately. The initial reference-policy fit uses the true
running cost and has no preceding compensation term.

\paragraph{Consistency of the random difference and product loss.}
The following proposition justifies the construction for true
policy evaluation. Replacing the sampled cost by a fixed integrable
compensated cost gives the corresponding consistency statement for
\eqref{eq:drdm-compensated-equation}.

\begin{proposition}[Random-difference consistency and the double-sample loss]
\label{prop:drdm-consistency-product}
Fix $\pi_k$, a candidate $\hat v_\psi$, and $0<h\le T-t$.
Suppose the two groups in \eqref{eq:drdm-block} are independent
conditional on $(t,x)$ and have finite first absolute moments.
Then \eqref{eq:drdm-product} holds. If, in addition,
\begin{equation}\label{eq:drdm-product-integrability}
\int\E[|R_1R_2|\mid t,x]\,\nu_k(\dd t,\dd x)<\infty,
\end{equation}
then \eqref{eq:drdm-objective} holds.

For weak consistency, assume also that $b$ and $\sigma$ are bounded,
the running-cost and log-density terms are integrable, and
$(\hat v_\psi)_{tt}$, $D_x^j\hat v_\psi$ for $1\le j\le4$, and
$D_x^j(\hat v_\psi)_t$ for $0\le j\le2$ are continuous and bounded.
Then, as $h\downarrow0$,
\begin{equation}\label{eq:drdm-consistency}
\mathcal R_h^{\pi_k}[\hat v_\psi](t,x)
=\mathcal R^{\pi_k}[\hat v_\psi](t,x)+O(h).
\end{equation}
The bound is uniform when the coefficient and derivative bounds are
uniform. The product identity does not require these additional
smoothness assumptions.
\end{proposition}

\begin{proof}
Conditional independence and identical sampling laws give
\[
\E[R_1R_2\mid t,x]
=\E[R_1\mid t,x]\E[R_2\mid t,x].
\]
Both means equal $\mathcal R_h^{\pi_k}[\hat v_\psi](t,x)$, proving
\eqref{eq:drdm-product}. Integrating under
\eqref{eq:drdm-product-integrability} gives
\eqref{eq:drdm-objective}.

For consistency, consider one action--Gaussian pair $(U,\xi)$ and
the associated $Y^\pm$. Conditional on $U=u$, freeze $b$ and
$\sigma$ at $(t,x,u)$ and $(t,x)$, respectively. Each endpoint has
the one-step law of the diffusion with these constant coefficients.
Applying the space--time Dynkin formula twice gives
\[
\E\!\left[
\frac{\tfrac12[\hat v_\psi(t+h,Y^+)+\hat v_\psi(t+h,Y^-)]
-\hat v_\psi(t,x)}{h}
\,\middle|\,t,x,U=u\right]
=\partial_t\hat v_\psi(t,x)+\LL^u\hat v_\psi(t,x)+O(h).
\]
The stated bounds control the remainder uniformly in $u$.
Equivalently, in the Taylor expansion the opposite noise increments
cancel the first-order random term, while
$\E[\xi\xi^\top]=I_r$ supplies the diffusion term
$\tfrac12\tr(A D_x^2\hat v_\psi)$. Averaging over
$U\sim\pi_k(\cdot\mid t,x)$ gives
$\partial_t\hat v_\psi+\LL^{\pi_k}\hat v_\psi+O(h)$.
The mean of the remaining terms is
\[
\int_{\R^m}\left[c(t,x,u)
+\tau\log\frac{p_k(u\mid t,x)}{\varrho_\mu(u\mid t,x)}\right]
\pi_k(\dd u\mid t,x)=\ell^{\pi_k}(t,x).
\]
Averaging $M_E$ copies leaves this mean unchanged, proving
\eqref{eq:drdm-consistency}. If each summand has finite variance,
its independent average has variance divided by $M_E$.
\end{proof}

\subsubsection{Constructing a Gibbs improvement teacher}
The input to policy improvement is the fitted value $\hat v_k$,
held fixed during this step. Its spatial gradient defines the energy
$E_{\hat v_k}$ in \eqref{eq:action-energy}; together with the reference
$\mu_{t,x}$ and temperature $\tau$, this energy defines the Gibbs
target $\hat q_k$ in \eqref{eq:critic-teacher}. The teacher's task is
to produce actions for training a policy toward this target. It uses
the energy and reference, without requiring an optimal-policy oracle.
For the density-based implementation below, actions lie in $\R^m$,
the reference density is positive, and the Gibbs normalizing constants
are finite and positive.

\paragraph{Step 1: draw and score candidate actions.}
At a training condition $(t,x)$, draw $M_A$ independent candidates
$U^{(j)}\sim\mu_{t,x}$. Evaluate their energies and normalize the
exponential weights:
\begin{equation}\label{eq:teacher-weights}
w_j(t,x)=\frac{\exp[-E_{\hat{v}_k}(t,x,U^{(j)})/\tau]}
{\displaystyle\sum_{i=1}^{M_A}\exp[-E_{\hat{v}_k}(t,x,U^{(i)})/\tau]}.
\end{equation}
The reference is already represented by the candidate draws, so these
weights contain only the energy factor. If candidates instead come
from a proposal density $p_{\mathrm{prop}}(\cdot\mid t,x)$, with support
covering the target, include the importance ratio
$\varrho_\mu(U^{(j)}\mid t,x)/p_{\mathrm{prop}}(U^{(j)}\mid t,x)$
in each unnormalized weight.

\paragraph{Step 2: resample a teacher action.}
Select candidate $j$ with probability $w_j(t,x)$ and set
$U^{\rm tar}=U^{(j)}$. Conditional on this candidate bank, the teacher
law is the discrete measure
$\sum_{j=1}^{M_A}w_j(t,x)\delta_{U^{(j)}}$, where $\delta_a$ denotes
unit mass at $a$. Refreshing the bank at each draw averages these
discrete laws over the candidate distribution
$\mu_{t,x}^{\otimes M_A}$. Neither construction is assumed to give
exact samples from the continuous Gibbs law at finite $M_A$.

For a factorized reference and an energy additive across groups, draw
and resample candidates independently within each group. For coupled
energies, a finite-run MCMC teacher is another option. The experiment
descriptions below specify the teacher used in each case. In every
case, its output is a teacher action $U^{\rm tar}$; any approximation
in that action's law remains part of the policy-improvement error.

\subsubsection{Conditional flow matching and policy density}
Flow matching takes the teacher actions as training data and learns
an $\R^m$-valued conditional velocity field $f_\theta(s,u;t,x)$.
The trainable object is this velocity field. Integrating it transforms
a base draw into a policy action, so the fitted policy can generate
new actions without constructing a teacher bank at every decision.
Here $t$ remains physical time, while $s\in[0,1]$ is generative time;
the condition $(t,x)$ stays fixed as $s$ evolves.

Choose a known positive base density $p_{\mathrm{base}}$ on $\R^m$
from which samples and log-density values are readily available.
A standard Gaussian $\mathcal N(0,I_m)$ is one possible choice.
This base initializes the flow. Its role differs from the reference
$\mu_{t,x}$ in the control objective and from the deployed-policy
density $p_k$, including the initial density $p_0$.

\paragraph{Step 1: form interpolation samples.}
Draw $(t,x)$ from $\nu_k$ and obtain a teacher action $U^{\rm tar}$
at that condition. Independently draw $Z\sim p_{\mathrm{base}}$ and
$s\sim\operatorname{Unif}[0,1]$, with all three draws conditionally
independent given $(t,x)$. Form
\[
Y_s=(1-s)Z+sU^{\rm tar}.
\]
This point lies on the straight line from the base sample to the
teacher action. The velocity of this training line is
$U^{\rm tar}-Z$, which can be computed directly from the endpoints.

\paragraph{Step 2: fit the velocity network.}
Feed $(s,Y_s;t,x)$ into $f_\theta$ and fit its output to the velocity
$U^{\rm tar}-Z$ by minimizing
\begin{equation}\label{eq:fm-objective}
\mathcal J_{A,k}(\theta)
=\E\left[\left\|f_\theta(s,Y_s;t,x)-(U^{\rm tar}-Z)\right\|^2\right].
\end{equation}
The expectation follows the sampling procedure in Step~1: first the
condition, then its teacher action, an independent base draw, and an
independent uniform generative time. Training uses minibatch averages
of this squared error and updates $\theta$, with the guide and teacher
procedure fixed. It requires evaluations of the velocity network at
the interpolated points, without integrating an ODE in this matching
objective. This is the conditional flow-matching regression of
\citet{LCBNL2023}.

\paragraph{Step 3: generate actions with the fitted flow.}
After training, freeze $\theta$. At a requested condition $(t,x)$,
draw a fresh $Z\sim p_{\mathrm{base}}$ and solve
\begin{equation}\label{eq:flow-ode}
\frac{\dd U_s}{\dd s}=f_\theta(s,U_s;t,x),\qquad U_0=Z.
\end{equation}
The endpoint $U_1$ is an action. Independent base draws give independent
actions at the same condition, and their law defines the exact ODE
policy. The numerical solver repeatedly evaluates $f_\theta$ to advance
$U_s$ from $s=0$ to $s=1$; no teacher endpoint is needed at deployment.
The resulting trajectory need not be a straight training line, because
the network combines information from many endpoint pairs.

\paragraph{Step 4: return the deployed law and density information.}
Under the regularity needed for an invertible flow, the exact ODE
policy has density $p_\theta$ satisfying
\begin{equation}\label{eq:flow-density-method}
\log p_\theta(U_1\mid t,x)
=\log p_{\mathrm{base}}(Z)-\int_0^1\nabla_u\cdot f_\theta(s,U_s;t,x)\dd s.
\end{equation}
Thus the base log density and the accumulated divergence provide
the log density of the generated action. The implementation evaluates
the divergence and integrates the action and log-density equations
numerically.

The deployed law includes this numerical integration. Let
$\Phi^{\Delta_k}_{t,x}(z)$ be the action returned from base input $z$
with generative step $\Delta_k$. For each Borel set $B\subseteq\R^m$,
\[
\pi_{k+1}(B\mid t,x)=
\int_{\R^m}\mathbf 1_{B}(\Phi^{\Delta_k}_{t,x}(z))p_{\rm base}(z)\,\dd z .
\]
When this law is absolutely continuous with respect to Lebesgue
measure, denote its density by $p_{k+1}$. The integrated log-density
output $\widetilde{\log p}_{k+1}$ approximates that density's logarithm;
\eqref{eq:flow-density-method} itself describes the exact ODE density
$p_\theta$. These densities need not coincide at a finite numerical
step. A subsequent true policy evaluation samples from $\pi_{k+1}$
and uses $\widetilde{\log p}_{k+1}-\log\varrho_\mu$ for its entropy
contribution. The compensated updates in Experiments~1 and~2 instead use the
cost specified in \eqref{eq:drdm-compensated-cost}.

\paragraph{Why the matching objective transports a distribution.}
Different endpoint pairs can pass through the same intermediate point
with different velocities. Squared-error regression learns their
conditional mean velocity. The following standard property
\citep{LCBNL2023} explains when this velocity transports the base law
to the teacher law.

\begin{lemma}[Conditional-mean velocity and ideal transport]
\label{lem:conditional-flow-transport}
Fix a control condition $(t,x)$ and let $U^{\rm tar},Z$ be the independent
teacher and base samples above, with finite conditional second moments.
Among unrestricted square-integrable velocity fields, the population
regression in \eqref{eq:fm-objective} is minimized by
\[
\E[U^{\rm tar}-Z\mid Y_s=u],
\]
where the expectation is under the teacher and base laws at this fixed
$(t,x)$, with $s$ held fixed. The equality holds almost everywhere under
the generative-time and interpolation laws. If this conditional-mean
velocity admits a sufficiently regular representative whose ODE flow
exists on the full interval $[0,1]$, and the associated continuity
equation has a unique probability-measure solution from the base law,
then the ideal flow has the same law as $Y_s$ at each $s$.
In particular, its endpoint distribution is the teacher law.
\end{lemma}

\begin{proof}
The conditional expectation minimizes squared prediction error, giving
the stated regression solution. To identify the transported law, take
a smooth compactly supported test function $\varphi$. At the fixed
$(t,x)$,
\[
\frac{\dd}{\dd s}\E[\varphi(Y_s)]
=\E[D\varphi(Y_s)\cdot(U^{\rm tar}-Z)].
\]
Conditioning the right-hand side on the interpolated position replaces
$U^{\rm tar}-Z$ by its conditional mean. This is the weak continuity
equation for the interpolation law. The ODE flow driven by that same
velocity satisfies this equation with the same initial base law;
uniqueness identifies the two laws on $[0,1]$.
\end{proof}

The endpoint law in this ideal statement is the teacher law. It equals
$\hat q_k$ when the teacher supplies exact Gibbs samples. For the
implemented policy, teacher approximation, velocity fitting, and ODE
integration each affect the final action law $\pi_{k+1}$.

\subsection{Experiment 1: Execution with a known solution}
\label{app:experiment-one}

\subsubsection{Problem setting}
\label{app:execution-model}

Write $x=(q,\alpha,y,z)\in\R^{4n}$, with inventory, signal, impact, and
liquidity blocks. The action is $u=(u^M,u^P)\in\R^{2n}$, with total speed
$\nu=u^M+u^P$. For $r\in\R^n$, angle brackets denote the asset average
$\langle r\rangle=n^{-1}\sum_i r_i$. Let $F$ have normalized sector
indicator columns, $\Pi_{\rm sec}=FF^\top$, and
$G=0.30I+0.08\Pi_{\rm sec}$. At a given execution speed, the dynamics are
\begin{equation}
\begin{aligned}
 \dd q_t&=\nu_t\,\dd t,&
 \dd\alpha_t&=-2\alpha_t\,\dd t+0.4\,\dd B_t^\alpha,\\
 \dd y_t&=(-4y_t+G\nu_t)\,\dd t+0.2\,\dd B_t^y,&
 \dd z_t&=-1.5z_t\,\dd t
 +0.35\bigl(0.25\,\dd B_t^\alpha+\sqrt{1-0.25^2}\,\dd B_t^z\bigr).
\end{aligned}
\label{eq:app-execution-dynamics}
\end{equation}
Under a feedback law $\pi$, averaging the drift over the action gives
\eqref{eq:experiment-dynamics}. The Brownian drivers are mutually independent
$n$-dimensional processes.
The initial state is independent of the Brownian drivers. Its
coordinates are independent, with
$q_i\sim\operatorname{Unif}(0.3,1)$,
$\alpha_i\sim\mathcal N(0,0.25^2)$,
$y_i\sim\mathcal N(0,0.10^2)$, and $z_i\sim\mathcal N(0,0.20^2)$.
This experiment has $n=2500$, $d=10{,}000$, $m=5000$, $T=1$,
terminal cost $g(x)=\tfrac12\langle q^2\rangle$, and temperature
$\tau=0.08/n$. The state has ten equal physical-state sectors,
which are distinct from the action blocks below.
Sampled-action Euler transitions approximate the controlled dynamics;
the diffusion remains action independent.

\paragraph{Reference action law.}
\label{app:benchmark-costs}

Define the orthogonal action coordinates
\[
 \xi=(u^M+u^P)/\sqrt2,\qquad \zeta=(u^M-u^P)/\sqrt2.
\]
These correspond to the implementation's action coordinates $(r,v)$;
the notation $\xi,\zeta$ distinguishes them from the value functions.
The public template is
$h_i(x)=\operatorname{clip}(-0.80q_i+0.15\alpha_i-0.08y_i,-1.15,1.15)$.
Partition $\xi$ into $500$ blocks of size five. With
$R_5=0.8I_5+0.2\boldsymbol1\boldsymbol1^\top$, the normalized reference density is
\begin{equation}
 \varrho_\mu(u\mid x)=
 \prod_b\left[\sum_{\varepsilon\in\{-1,1\}^5}
 w_b^{\rm ref}(\varepsilon\mid x)
 \varphi_{0.075^2 R_5}\bigl(\xi_b-h_b/\sqrt2-0.30\varepsilon\bigr)\right]
 \prod_i\varphi_{0.075^2}(\zeta_i),
\label{eq:app-half-reference}
\end{equation}
where $\varphi_\Sigma$ is a centered Gaussian density with covariance
$\Sigma$. The weights are normalized over all sign patterns, with
\[
 w_b^{\rm ref}(\varepsilon\mid x)\ \propto\
 \exp\!\left\{\sum_{i\in b}(0.10q_i+0.20\alpha_i+0.15z_i)\varepsilon_i
       +0.20\sum_{j=1}^5\varepsilon_j\varepsilon_{j+1}\right\},
 \qquad \varepsilon_6=\varepsilon_1.
\]
Thus half the action coordinates are multimodal and half are Gaussian.
Both latent sign interactions and within-component covariance couple the
multimodal coordinates inside each block. The action density still
factorizes between blocks; this is structured block coupling, rather
than arbitrary global dense coupling. The reference kernel
$\mu_{t,x}(\dd u)=\varrho_\mu(u\mid x)\dd u$ is fixed during training.

\paragraph{Known value, optimal policy, and supplied cost.}
The prescribed optimal value is
$v_\tau(t,x)=g(x)+(1-t)f(x)$, where
\begin{equation}
\begin{aligned}
 f(x)={}&0.06+0.06\langle q^2\rangle+0.04\langle\alpha^2\rangle
 +0.025\langle y^2\rangle+0.025\langle z^2\rangle
 -0.05\langle q\alpha\rangle\\
 &+0.018\langle\sin(3q)\rangle
 +0.03\langle(\Pi_{\rm sec}q)^2\rangle.
\end{aligned}
\label{eq:app-manufactured-value}
\end{equation}
The exact optimal policy has the same block layout as
\eqref{eq:app-half-reference}, with covariance $0.072^2R_5$,
Gaussian-half variance $0.072^2$, multimodal amplitude
$0.30[1+0.08\sin(2\pi t)]$, and normalized sign weights
\[
 w_b^*(\varepsilon\mid t,x)\ \propto\
 \exp\!\left\{\sum_{i\in b}
 [-0.25+0.45\sin(2\pi t)+0.35\alpha_i+0.40z_i+0.15q_i]\varepsilon_i
 +0.65\sum_{j=1}^5\varepsilon_j\varepsilon_{j+1}\right\}.
\]
The supplied running cost is
\begin{equation}
 c(t,x,u)=-\partial_t v_\tau(t,x)-\LL^u v_\tau(t,x)
 -\tau\log\frac{p^*(u\mid t,x)}{\varrho_\mu(u\mid x)}.
\label{eq:app-manufactured-cost}
\end{equation}
Exact values, gradients, normalized optimal-policy parameters, and
optimal-policy sampling are reserved for evaluation. Training queries the
supplied public cost and energy/score interfaces.

\subsubsection{Training protocol}
\label{app:finite-training}

\paragraph{Inputs, outer loop, and reported outputs.}
The four combinations pair a value solver, DRDM or PINN, with a policy
sampler, Flow or Langevin. Each starts from the reference policy and
fits its value approximation $\hat v_0$. For $k=0,1$, one improvement
and the following value fit proceed as follows:
\begin{enumerate}
\item Freeze $\hat v_k$ and construct the Gibbs target
$\hat q_k=q_{\hat v_k}$ from its spatial gradient.
\item Train or run the chosen sampler for this target, obtaining the
deployed law $\pi_{k+1}$. Use this law to generate new state
trajectories and action samples.
\item Freeze both $\pi_{k+1}$ and the guide $\hat v_k$, and fit
$\hat v_{k+1}$ using the compensated equation described below.
\end{enumerate}
Thus two policy improvements require three value fits. The reported
value approximation is $\hat v=\hat v_2$, whereas the deployed policy
is $\hat\pi=\pi_2$, constructed from $\hat v_1$.

Both value solvers use the raw-state representation
\begin{equation}
 \hat v_\psi(t,x)=g(x)+0.1(1-t)N_\psi(t,x),
\label{eq:experiment-raw-critic}
\end{equation}
where $N_\psi$ has two $192$-wide SiLU hidden layers and takes only
unscaled time and state inputs. The terminal envelope enforces
$\hat v_\psi(1,x)=g(x)$. The solvers use different inner training
budgets, specified below.

\paragraph{The equation used in each value fit.}
The initial fit evaluates the reference policy, for which the
relative-entropy penalty is zero. The two subsequent fits use the
compensated cost \eqref{eq:drdm-compensated-cost}, with the previous
value $\hat v_k$ held fixed. In the notation of
Section~\ref{sec:numerical-method}, the residual fitted at these stages is
\[
\partial_t\hat v_{k+1}
+\LL^{\pi_{k+1}}(\hat v_{k+1}-\hat v_k)
+\HH_\tau[\hat v_k].
\]
For the benchmark drift $b=b_0+B\nu$, where $\nu=u^M+u^P$, the
control Hamiltonian is
\[
H_{\rm ctrl}(B^\top D_x\hat v_k)
=-\tau\log\int_{\R^{2n}}
\exp\!\left[-\frac{c(t,x,u)+(B^\top D_x\hat v_k)\cdot\nu}{\tau}\right]
\mu_{t,x}(\dd u),
\]
with $B^\top D_x\hat v_k=D_q\hat v_k+G^\top D_y\hat v_k$.
Since $\HH_\tau[\hat v_k]=\LL^0\hat v_k+
H_{\rm ctrl}(B^\top D_x\hat v_k)$, the compensated cost can be
computed as
\[
C_{\hat v_k}(t,x,u)
=H_{\rm ctrl}(B^\top D_x\hat v_k)
-(B^\top D_x\hat v_k)\cdot\nu.
\]
Here $\LL^0$ contains the uncontrolled drift and state diffusion;
these action-independent terms cancel from the compensated cost.
This is the benchmark implementation of
\eqref{eq:drdm-compensated-equation}. The true evaluation residual
contains the additional term
$\tau\KL(\pi_{k+1}\Vert\hat q_k)$ from
\eqref{eq:drdm-compensation-gap}. Thus $\hat v_2$ is the value
approximation produced by the compensated solve, and is not identified
with the true soft value of $\pi_2$. We measure its nonlinear HJB
residual and the achieved cost of $\pi_2$ separately. Throughout this
appendix, $p_k$ denotes the policy density.

\paragraph{Constructing the Gibbs target.}
With $\hat v_k$ fixed, the public Hamiltonian implementation reconstructs the declared
common-covariance block-mixture energy from public log-energy and score
queries, completes squares, and checks reconstruction at interior points.
Its teachers sample the reconstructed public Gibbs target. It does not
read the evaluation-only value or normalized optimal-policy helper.
This integration procedure uses the benchmark's declared mixture
structure; it is not a general integration method for arbitrary energies.

\paragraph{Turning the target into a deployed policy.}
Flow fits the conditional velocity field to teacher samples using
\eqref{eq:fm-objective}; a fresh base sample is then transported by
the learned field to generate an action. Langevin instead updates an
action repeatedly using the current target's score, rebuilt from
$\hat v_k$ at each improvement. The resulting flow map or finite chain
defines $\pi_{k+1}$, which supplies the samples for the next value fit.

Both samplers start from $\mathcal N(0,0.6^2I)$. Each Langevin chain
takes $512$ unadjusted steps with action step $10^{-4}$; no target-mode
initialization is supplied. Flow uses local and block-conditioned SiLU
networks of widths $96$ and $64$, respectively, with dimensions
$11\to96\to96\to8$ and $7\to64\to64\to3$. Its deployment uses
$16$ RK4 steps and complete analytic divergence. The Langevin action
step, flow generative-time step, and physical Euler step are distinct.

\paragraph{Two-round DRDM policy iterations.}
At each value stage, DRDM first builds a state pool from the frozen
policy's trajectories. Each minibatch then uses the random neighbors
and two-group product loss described above, with the actual cost in
the initial stage and $C_{\hat v_k}$ after each improvement. The
selected critic is passed to the next sampler update; the third fit
produces the reported $\hat v_2$.

Each of the three fits simulates $64$ fresh paths on a $96$-step
physical grid and takes $2400$ AdamW updates with a global minibatch
of $32$. Each residual group contains two action draws and
antithetic physical noise. Every $100$ updates we evaluate the residual
on $128$ held-out time--state pairs and retain the critic with the
lowest validation RMS.

For each Flow improvement, the frozen critic supplies two teacher
actions at each of $256$ contexts. Actor training uses minibatches
of $128$ and starts at learning rate $4\times10^{-4}$. We check
validation KL every $500$ updates, reduce the learning rate on
stagnation, and stop after two stable $3000$-update windows at the
minimum rate $1.25\times10^{-5}$. The two fits take $30{,}500$ and
$21{,}500$ updates in this experiment. Each improvement returns the flow with the lowest
validation KL. Each DRDM run uses seed $20261014$ and four V100 GPUs;
the final policy is rolled out with $96$ Euler steps.

\paragraph{DRDM learning-rate selection and refinement.}
We selected a common critic learning rate for the four DRDM runs in
Experiments~1 and~2. On each problem we tried $2\times10^{-4}$, $5\times10^{-5}$
and $2\times10^{-5}$ for $1200$ updates under the reference policy,
with the same initialization, draws and validation pairs across rates.
For each rate we divided the best validation RMS by its initial value
on each problem. The geometric mean of these two ratios was $0.0558$,
$0.0547$ and $0.0618$, respectively. All four runs therefore use
$5\times10^{-5}$.

The earlier two-round runs used $24$
physical steps, a critic learning rate of $2\times10^{-4}$, and $180$
updates per fit. On the same known-solution test states, the revised
settings reduced the relative value error from $2.80\%$ to $0.68\%$
for flow and from $1.49\%$ to $0.45\%$ for Langevin. Their full-gradient
errors increased from $26.1\%$ to $27.4\%$ and from $22.8\%$ to
$25.9\%$, respectively.

\paragraph{Two-round PINN policy iterations.}
\label{app:pinn-policy-iteration}
The PINN variants follow the same outer loop and use
\eqref{eq:experiment-raw-critic}. At a sampled time--state point,
automatic differentiation supplies the value derivatives, including
the diffusion trace and its parameter gradient exactly for the SiLU
network. Independent action samples form two estimates of the
strong-form residual; their product supplies the critic loss.
The initial stage uses the reference policy's running cost, and the
two later stages use $C_{\hat v_k}$ with the guide held fixed.

In each problem, the initial PINN fit starts from random weights and
provides the common $\hat v_0$ for PINN--flow and PINN--Langevin.
At each improvement, PINN--flow fits the same block-conditioned flow
as DRDM--flow to samples from its current critic's public block Gibbs
law. PINN--Langevin runs the $512$-step chain using its own critic's
score. Each returned policy then generates fresh trajectories and
actions for the next PINN fit.

Each evaluation uses $32$ training paths and $8$ independent validation
paths on a $24$-step grid. A training batch of $32$ states uses two
independent action pairs drawn from separate banks of four pairs per
state. The $128$ validation states use $16$ draws per block.
Each critic stage takes $2000$ AdamW updates and selects the checkpoint
with the smallest validation residual RMS. Each Flow improvement uses
$256$ teacher contexts with two targets per context and the DRDM--flow
stopping rule above. Reverse KL to the public Gibbs target is measured
every $500$ updates on independent validation states. The two actor fits
stop at $38{,}500$ and $25{,}000$ updates in this experiment. Each round deploys the checkpoint with
the lowest validation KL, and the second fit starts from the first
selected actor. Both PINN variants use seed $20261014$ on an RTX~3090.
We evaluate $\hat v_2$ on the common held-out states and trace probes and
deploy the second improved policy on the common $64$ initial states
and Brownian paths with $96$ Euler steps.

\paragraph{Baseline comparison scope.}

\label{app:retained-baselines}

The baseline methods use the same half-mixture problem formulation
and test pools as the four two-round runs. Known-case errors use the same $256$
state indices, retaining all $10{,}000$ gradient coordinates whenever a
scalar value network is available. Different architectures, budgets, state-sampling
protocols, and numerical grids are retained; this is not an equal-capacity,
equal-query, or common-training-operator comparison.

\paragraph{Continuous-time q-learning.}
The q-learning implementation \citep{JZ2023q} learns a value and a
normalized stochastic policy together. Each update samples
trajectories from the current policy and forms martingale increments
from consecutive reward-value estimates and the running reward minus
the temperature times the policy log density. These increments are
held fixed in the stochastic-approximation updates of the value and
policy parameters. The next-state value is also held fixed. This implementation
works in the reward convention; the reported value is converted to
the cost convention by changing its sign. Both networks take the raw
time and state and have two $64$-wide tanh hidden layers. The trained
actor is the deployed policy.

\paragraph{Deep Splitting.}
Deep Splitting \citep{BBCJN2021} represents the value by a separate
network at each of $24$ time layers. Starting from the terminal data,
it fits these networks backward in time. At each layer, the next
network and its gradient are frozen while forming the regression
target; the current network is then fitted to that target before
moving to the preceding layer. Each network has two $(d+10)$-wide
ReLU hidden layers and batch normalization. Evaluation interpolates
the native layers and uses the recorded right-limit time derivative
at grid knots. The fitted value gradient defines the
normalized Gibbs policy used for policy KL and rollout costs.

\paragraph{Reference-policy initialization of q-learning.}
\label{app:mu-baseline-initialization}
For Experiments~1 and~2, the actor uses
$U=\delta_\theta(t,x)+\exp(\ell_\theta(t,x))\odot Z$, $Z\sim\mu_x$,
where $\delta_\theta,\ell_\theta\in\R^m$ are the shift and log-scale
outputs of the original two-$64$-wide tanh network. Its final affine
layer is initialized to zero, giving $\pi_0=\mu$ exactly; the remaining
actor layers and reward-value network retain their seeded random
initialization. The normalized density is evaluated by change of
variables, including the log-Jacobian $-\sum_{j=1}^m\ell_{\theta,j}$.
The martingale stochastic-approximation update is unchanged, using
running reward $-c+\tau\log\varrho_\mu$, little-q equal to
$\tau\log\varrho_{\pi_\theta}$, and step size $0.01/i^{0.55}$ at
stochastic-approximation update $i\geq1$. This adaptation changes the policy family as well as
its initialization. Training uses seed $20261014$, $20{,}000$ updates,
and global batch size $32$.

\paragraph{Reference-policy initialization of Deep Splitting.}
Since Deep Splitting parameterizes the value alone, we first approximate
$v^\mu$ in a separate backward linear policy-evaluation pass with
terminal condition $g$. The Euler reference process uses the drift
averaged under $\mu$, and its regression target uses the next fitted
value and $\E_{U\sim\mu_x}[c(t,x,U)]$, estimated with $32$ action
draws per state. The relative-entropy term vanishes under $\mu$.
In Experiment~1, this pass fits $24$ time-layer networks with $600$ updates per layer
and global batch size $256$. Before fitting each layer in the nonlinear
soft-HJB pass, we load the corresponding $v^\mu$ network weights and
batch-normalization buffers and initialize a fresh optimizer. The
nonlinear pass retains the original architecture, regression update,
and $24\times600$ updates. Thus the initialization adds $14{,}400$ updates,
giving $28{,}800$ updates in total. Both passes use seed $20261014$; the
deployed policy is the Gibbs law of the final fitted value gradient.

\paragraph{Training budgets for q-learning and Deep Splitting.}
The Experiment~1 update counts, including the separate reference-policy
evaluation, and global batch sizes are
\begin{center}\small
\begin{tabular}{lrr}
\hline
Method & Updates & Global batch\\
\hline
q-learning & $20{,}000$ & $32$\\
Deep Splitting: $\mu$ evaluation & $24\times600$ & $256$\\
Deep Splitting: soft HJB & $24\times600$ & $256$\\
\hline
\end{tabular}
\end{center}
\paragraph{Projected mirror descent.}
Mirror descent \citep{KSSZ2025Mirror} updates a policy represented
through its Brownian history. In each outer round, it first simulates
paths under the current policy and regresses a causal adjoint on that
history. This fitted adjoint supplies the next log-geometric density
update. The update is projected into a common-covariance
$32$-component mixture per five-dimensional block, while the Gaussian
half is represented exactly. The resulting method is therefore a
projected adaptation of the density update.

The adjoint network has two $32$-wide tanh layers. Training uses
$10$ outer rounds, $256$ paths per round, and $2000$
adjoint-regression updates per round on a native $24$-step history
grid. Deployment retains this grid. For conditional value estimates,
the history is restarted at rounded native time nodes, and the
known-case value error is computed from the saved independent Monte
Carlo policy-value estimates. Policy KL uses the first action of the
same restarted history. These outputs do not provide a scalar value
network, its gradient, an HJB residual, or a feedback law for arbitrary
histories, so those quantities are not reported for this method.

\paragraph{Deployed-cost evaluation.}
The q-learning and Deep Splitting deployed costs use their recorded
$96$-step evaluations, $64$ paths, and the same physical-noise seed as
the four two-round controller evaluations. Mirror's cost uses its native
$24$-step history grid and $64$ paths.

\paragraph{Released PINN-SPI.}
\label{app:released-pinn-spi}
PINN-SPI \citep{KCKK2026} alternates actor updates and frozen-policy
value updates. We append time to the network inputs and impose the
critic's terminal condition through
$\hat v_\psi(t,x)=g(x)+(1-t)N_\psi(t,x)$. The release's reward-value
convention maps to $-\hat v_\psi$, with running reward
$-c+\tau\log\varrho_\mu$.

For Experiments~1 and~2, the actor starts exactly from the reference
law $\mu_x$. To represent its multimodal structure, we replace the
released diagonal-Gaussian policy family by
\[
 U=\delta_\theta(t,x)+\exp(\ell_\theta(t,x))\odot Z,
 \qquad Z\sim\mu_x,
\]
where $\delta_\theta,\ell_\theta\in\R^m$ are the shift and log-scale
outputs. The final affine layers of both heads are initialized to
zero, giving $\pi_0=\mu$. The log density is evaluated by the
change-of-variables formula, including the log-Jacobian
$-\sum_{j=1}^m\ell_{\theta,j}$. This adaptation changes the policy
family while retaining the released network architecture and update rule.

One actor update follows the
\href{https://github.com/tomatofromsky/pinn-spi-aaai2026/blob/a1b573a637ef77271bbbc23346cdbabe68d7a98d/pinn_spi/algs/pinn_spi.py\#L635-L675}{authors' released policy update}:
draw reparameterized actions at each training state, evaluate the
public unnormalized Hamiltonian on those actions, and form softmax
weights with their gradients detached. Adam then minimizes the
weighted negative log-likelihood. Terms independent of the action
cancel in the softmax, and the action-to-log-density gradient follows
the released computation graph. After the actor updates, freeze the
policy and fit the sampled policy-evaluation residual
$\mathcal R^\pi[\hat v_\psi]$, using the benchmark's diffusion
covariance and the release's detached trace gradient.

The $100$-wide residual SiLU actor and critic have depths $4$ and
$2$ and use the paper's Adam learning rates $3\times10^{-5}$ and
$3\times10^{-4}$. The critic and the actor's remaining layers retain
their random initialization. Each
problem uses seed $20261014$, $100$ outer rounds, and $10$ actor
updates followed by $10$ critic updates per round. Each batch uses
$32$ states from the training pool and $100$ actions per
state. The paper's compact LQR experiment uses $1000$ rounds and
$100$ states per batch.

A seeded check reproduces one complete Adam policy step of the
authors' unmodified routine. The policy loss printed with
Algorithm~2 of \citet{KCKK2026} is a candidate KL; the reported runs
use the loss in the released code. Under the above reparameterization,
the sampled log density is
$\log\varrho_\mu(x,Z)-\sum_{j=1}^m\ell_{\theta,j}$, so the shift
head receives zero direct gradient and its zero output is preserved.
In both reported runs, the residual scale reaches $e^{-2}$, the
lower limit of the released log-scale clamp $[-2,5]$.

\subsubsection{Evaluation and reporting}
\label{app:action-numerics}
\label{app:experiment-evaluation}

The evaluation uses three different objects: the fitted value network,
the deployed action distribution at fixed test states, and the
closed-loop trajectories generated by that distribution. Their errors
and costs are measured separately.

\paragraph{Value and PDE residual.}
Known-problem value and full-gradient errors use the same frozen
$256$-state test pool and all $10{,}000$ gradient coordinates,
with Mirror's native-time convention described
above. To evaluate the soft HJB residual in
\eqref{eq:experiment-hjb-metric}, substitute the fitted network into
the nonlinear equation at $M=32$ held-out time--state pairs and estimate
the diffusion trace with eight random probes per pair. This diagnostic
is distinct from the linear residual minimized during a value fit.

\paragraph{Policy error at fixed states.}
For the known problem, the policy score is
\[
 \mathcal K(\pi)=\frac1{64}\sum_{i=1}^{64}
 \KL\bigl(\pi(\cdot\mid t_i,x_i)\Vert
               \pi^*(\cdot\mid t_i,x_i)\bigr).
\]
The $64$ states are the unchanged prefix of a historical test pool
generated from reference-policy trajectories and are correlated within
trajectories. For Flow, draw $128$ actions at each state, compute their
log densities using the flow divergence, and average the log-density
ratios against the known optimal density. The reported uncertainty is
conditional action MCSE, with no IID-state standard error. For the
finite Langevin sampler, the unknown endpoint density requires the
bounds described below.

\paragraph{Closed-loop costs.}
Freeze the trained controller and simulate $64$ paths per method and
problem. Use $96$ physical Euler steps for all methods except Mirror,
which uses $24$. The unregularized cost $J_0$ adds the running and
terminal costs along these trajectories; $J_\tau$ also adds the
cumulative reference-KL penalty in \eqref{eq:regularized-cost}.
These quantities evaluate the actual deployed sampler. The state
weights in the fixed-state policy score need not match the occupation
weights along a controller's own trajectories.

\paragraph{Plotting protocol.}
In Figure~\ref{fig:known-time-errors}, group the held-out pairs into
six time intervals of width $1/6$ and compute $E_V$ within each
interval. Figure~\ref{fig:known-policy-scatter} uses the first shared
test pair, at $t=1/12$, and $128$ action draws per panel, projected
onto $(\xi_1,\xi_2)$ of the first block. The comparison law is the
corresponding marginal of $\pi^*$. Solid and dashed contours enclose
its $50\%$ and $90\%$ highest-density regions. The empirical
$\mathrm{SW}_2$ is computed against the first panel by averaging
squared one-dimensional Wasserstein distances over $360$ equally
spaced projection directions and taking the square root. Experiments~2
and~3 use the same contour levels and distance calculation with their
respective comparison laws.

\paragraph{Density supplied by a flow.}
To evaluate an action's log density, propagate its known base density
along the same ODE used to generate that action. For
$\dd U_s/\dd s=f_\theta(s,U_s;t,x)$, the density calculation is
\[
 \log p_\theta(U_1\mid t,x)
 =\log p_{\rm base}(U_0)-\int_0^1
 \operatorname{div}_u f_\theta(s,U_s;t,x)\,\dd s,
\]
as in \eqref{eq:flow-density-method}. The implementation integrates this
identity numerically and returns the action together with its computed
log density. Policy KL then combines this numerical-density
approximation with Monte Carlo action averaging. The reported action
MCSE does not include ODE integration bias.

\paragraph{The missing density for a finite Langevin chain.}
\label{app:langevin-missing-density}
Every Langevin step adds fresh noise. The known transition
densities determine the joint density of a full sampling path, while
the final action density requires integrating over all intermediate
actions. Write $\pi_K(\cdot\mid t,x)$ for the finite-step endpoint law
and $p_K(\cdot\mid t,x)$ for its Lebesgue density; here $K$ counts
Langevin steps within one frozen sampler. At fixed $(t,x)$,
the known-solution policy error is
\[
 \KL(\pi_K\Vert\pi^*)
 =\int_{\R^m}\log\frac{p_K(u\mid t,x)}{p^*(u\mid t,x)}
                 \,\pi_K(\dd u\mid t,x).
\]
The unavailable term is $\log p_K$, even when $\log p^*$ is known.

\paragraph{Construct the lower and upper scores.}
An auxiliary normalized density gives a lower bound on expected log
density, and independent witnesses tighten the lower score. For the
DRDM--Langevin policy KL, split the $512$ actions at each state into
$256$ witness-fit draws and $256$ independent test draws for the
Donsker--Varadhan lower score. Apply a composed Hoeffding penalty to
obtain a one-sided $95\%$ lower confidence limit, conditional on the
fixed states, guide, and fitted witnesses.

For the upper population bound, use the common-precision mixture
structure to propagate a lower entropy certificate through each
Langevin block by the entropy-power recursion. Evaluate the
independent Gaussian half exactly and estimate the remaining score
by Monte Carlo. This is an estimate of an upper population bound,
not an upper confidence limit. Neither calculation substitutes the
Gibbs target density for the density of the finite sampler.

\paragraph{Apply the bounds to the deployed soft cost.}
The same unknown density enters $\KL(\pi_K\Vert\mu_{t,x})$ along the
closed-loop paths. Here the comparison law is the reference $\mu_{t,x}$;
the policy-error metric instead compares with the optimal law $\pi^*$.
Using the entropy certificates in the trajectory cost gives population
bounds $L_\tau\leq J_\tau(\pi_K)\leq U_\tau$ for this deployed
controller. These bounds concern its achieved cost, rather than the
unknown optimal value. The displayed path MCSEs measure sampling
uncertainty in the estimated bounds and do not turn them into a
calibrated two-sided confidence interval.

\begin{table}[!htbp]
\centering
\caption{Experiment~1: finite-Langevin bounds for policy error \eqref{eq:experiment-policy-metric}
and soft cost evaluated by rollout \eqref{eq:experiment-cost-metric}.
The lower entry for $\mathcal K(\pi_K)$ is a one-sided $95\%$ confidence
limit. All other entries estimate population bounds. Standard errors
use action draws for $\mathcal K(\pi_K)$ and trajectories for $J_\tau(\pi_K)$.}
\label{tab:finite-langevin-bounds}
\begingroup\small
\begin{tabular}{lllrr}
\hline
Method & Experiment & Quantity & Lower & Upper\\
\hline
DRDM--Langevin & 1 & $\mathcal K(\pi_K)$ & $877.8692$ & $2633.4829\pm0.3398$\\
DRDM--Langevin & 1 & $J_\tau(\pi_K)$ & $0.364174\pm0.000407$ & $0.439109\pm0.000396$\\
PINN--Langevin & 1 & $J_\tau(\pi_K)$ & $0.364555\pm0.000402$ & $0.438975\pm0.000395$\\
\hline
\end{tabular}
\endgroup
\end{table}

Since $\mu_{t,x}$ is normalized and $\tau>0$, KL nonnegativity also gives
\[
 J_\tau(\pi_K)\geq J_0(\pi_K),\qquad
 J_\tau(\pi_K)\geq\max\{L_\tau,J_0(\pi_K)\}.
\]
For the two-round DRDM controllers in this experiment, the mean flow
soft cost is $0.370928$, below the Langevin ordinary-cost estimate
$0.378911$. The corresponding Experiment~2 results are reported in
Appendix~\ref{app:experiment-two}.
This comparison uses Monte Carlo point estimates. The displayed standard
errors quantify trajectory sampling for each trained controller.

This experiment uses the recorded training seed $20261014$, so the
Monte Carlo errors do not quantify training-seed variability.
The same density-access and entropy-bound procedures are used for
Experiment~2, whose cost bounds appear in
Table~\ref{tab:finite-langevin-bounds-experiment-two}.

\subsection{Experiment 2: Execution with an unknown solution}
\label{app:experiment-two}

\subsubsection{Problem setting}
This experiment keeps the state, action coordinates, initial-state law,
and form of the dynamics in Appendix~\ref{app:execution-model}, with $n=2500$,
$d=10{,}000$, $m=5000$, $T=1$, $g(x)=\tfrac12\langle q^2\rangle$,
and $\tau=0.08/n$. The sector projection now uses eight sectors of
sizes $313,313,313,313,312,312,312,312$, which also changes the impact
coefficient $G=0.30I+0.08\Pi_{\rm sec}$. The reference action density
is still \eqref{eq:app-half-reference}, with the same public template
and five-dimensional action blocks. Besides this sector partition,
the running cost is changed as specified below. No exact optimal value
or policy is supplied.

\paragraph{Running cost.}
Set $\widetilde z=\operatorname{clip}(z,-2,2)$,
$q_z=q\odot e^{\widetilde z/2}$,
$\lambda(z)=0.12+0.03\cosh(\widetilde z)$, and $\xi^c=\xi-h(x)/\sqrt2$.
Let $L_5$ be the cyclic five-site graph Laplacian. The running cost is
\begin{equation}
\begin{aligned}
 c_2(x,u)={}&-\langle q\alpha\rangle
 +\frac{0.225}{n}q_z^\top(I+0.30\Pi_{\rm sec})q_z\\
 &+\left\langle(y-q)\sqrt2\xi+\lambda(z)(\xi^c)^2+0.15\zeta^2\right\rangle
 +\frac{0.035}{n}\sum_b(\xi_b^c)^\top L_5\xi_b^c.
\end{aligned}
\label{eq:app-coupled-cost}
\end{equation}
The term $-\langle q\nu\rangle$ is a public terminal-potential anchor.
The reference-KL penalty is added separately. No exact value or optimal
policy is prescribed for Experiment~2, and its true value, gradient, and
optimal-policy KL errors are not reported.

\subsubsection{Training protocol}
\paragraph{Value and policy updates.}
Use the outer loop, raw-state networks, Gibbs construction, and deployed
samplers described in Appendix~\ref{app:finite-training}, replacing
$c$ by $c_2$ throughout the action energy and control Hamiltonian.
First evaluate the reference policy. Then perform two sampler updates,
each followed by a value fit with the deployed law and previous guide
held fixed. The later fits use $C_{\hat v_k}$, so the distinction
between the compensated equation and true policy evaluation in
\eqref{eq:drdm-compensation-gap} applies here as well. The outputs are
the final approximation $\hat v_2$ and the policy $\pi_2$ constructed
from $\hat v_1$.

\paragraph{Critic and sampler settings.}
The critic architecture remains \eqref{eq:experiment-raw-critic}.
Each DRDM fit uses $64$ paths on a $96$-step grid, $2400$ AdamW
updates, a global minibatch of $32$, and learning rate $5\times10^{-5}$.
Each PINN fit uses $32$ training and $8$ validation paths on a
$24$-step grid, $2000$ AdamW updates, and batches of $32$ states.
Both retain their residual sampling and validation rules from
Appendix~\ref{app:finite-training}. The common DRDM learning rate was
selected jointly over the two finite-horizon problems as described there.

Flow and Langevin retain the source $\mathcal N(0,0.6^2I)$.
Langevin takes $512$ steps of size $10^{-4}$, and Flow uses the same
block-conditioned networks, $16$ RK4 steps, and analytic divergence.
The teacher contexts, actor minibatches, learning-rate schedule, and
validation stopping rules are unchanged. The two Flow training stages
take $42{,}500$ and $21{,}500$ updates for DRDM--flow, and $40{,}500$
and $23{,}500$ for PINN--flow. All runs use seed $20261014$;
the DRDM runs use four V100 GPUs and the PINN runs use an RTX~3090.

\paragraph{Baseline settings.}
The q-learning, Deep Splitting, mirror-descent, and PINN-SPI
algorithms and their deployment conventions are those described in
Appendix~\ref{app:retained-baselines}, now using $c_2$.

q-learning and Deep Splitting use the reference-policy initialization in
Appendix~\ref{app:mu-baseline-initialization}, with $c_2$ in place of $c$.
q-learning retains $20{,}000$ updates. Deep Splitting uses $30$ updates
per time layer in each pass: $720$ updates for the separate $v^\mu$
evaluation and $720$ for the nonlinear soft-HJB solve, giving $1440$
updates in total.

The baseline training budgets for this experiment are
\begin{center}\small
\begin{tabular}{lrr}
\hline
Method & Updates & Global batch\\
\hline
q-learning & $20{,}000$ & $32$\\
Deep Splitting: $\mu$ evaluation & $24\times30$ & $256$\\
Deep Splitting: soft HJB & $24\times30$ & $256$\\
\hline
\end{tabular}
\end{center}
Mirror retains $10$ outer rounds, $256$ paths per round, and $2000$
adjoint-regression updates per round on its $24$-step history grid.
PINN-SPI retains $100$ outer rounds, with $10$ actor and $10$ critic
updates per round, $32$ states per batch, and $100$ actions per state.
Its running reward is $-c_2+\tau\log\varrho_\mu$.

\subsubsection{Evaluation and reporting}
\label{app:experiment-two-evaluation}
\paragraph{Residual and closed-loop costs.}
Evaluate the nonlinear soft HJB residual of the fitted value on the
same $32$ held-out diagnostic states for every method that supplies a
scalar value network, using eight diffusion-trace probes per state.
Since no exact solution is prescribed, true value, gradient, and
optimal-policy KL errors are not reported. The action plots instead
compare deployed samples with the reference law at the same test state.

Evaluate each deployed controller on $64$ closed-loop paths, using the
common initial states and physical Brownian increments. All methods
use $96$ Euler steps except Mirror, which retains its native $24$-step
history grid. The ordinary cost uses $c_2$ and $g$; the soft cost also
includes the reference-KL penalty. Density calculations and finite-chain
entropy bounds follow Appendix~\ref{app:langevin-missing-density}.

\paragraph{Plotting protocol.}
Figure~\ref{fig:problem-c-policy-scatter} uses $128$ action draws per
panel at the first shared test pair, $t\approx0.9321$, projected onto
$(\xi_1,\xi_2)$ of the first block. Mirror uses the same state at the
nearest native time node, $t=22/24$: the state is treated as a new initial
condition with no preceding Brownian history, and the trained updates
produce the first-action distribution without retraining.
The comparison law is the
reference $\mu$. Contour levels and the sliced Wasserstein calculation
follow Appendix~\ref{app:experiment-evaluation}.

\begin{table}[!htbp]
\centering
\caption{Experiment~2: finite-Langevin bounds for the deployed soft cost.
Both entries estimate population bounds; the displayed standard errors
use closed-loop trajectories and do not define a calibrated confidence interval.}
\label{tab:finite-langevin-bounds-experiment-two}
\begingroup\small
\begin{tabular}{lllrr}
\hline
Method & Experiment & Quantity & Lower & Upper\\
\hline
DRDM--Langevin & 2 & $J_\tau(\pi_K)$ & $0.327172\pm0.000556$ & $0.396401\pm0.000545$\\
PINN--Langevin & 2 & $J_\tau(\pi_K)$ & $0.327206\pm0.000516$ & $0.396169\pm0.000508$\\
\hline
\end{tabular}
\endgroup
\end{table}

\paragraph{Interpretation of the estimates.}
For the two-round DRDM controllers, the mean flow soft cost is
$0.327255$, below the Langevin ordinary-cost estimate $0.337573$.
This is a comparison of Monte Carlo point estimates. The displayed
standard errors quantify trajectory sampling for each trained
controller; they do not quantify training-seed variability or turn the
estimated population bounds into a calibrated confidence interval.
The residuals, policy plots, and rollout results appear in
Section~\ref{subsec:experiment-no-exact}.

\subsection{Experiment 3: Execution over an infinite horizon}
\label{app:stationary-routing}

\subsubsection{Problem setting}

The state dynamics are \eqref{eq:app-execution-dynamics} with $n=256$
assets in two sectors of $128$, so $d=1024$ and $m=512$. The discount
is $\beta=0.5$, the temperature is $\tau=0.08/n$, and the reference
law is $\mathcal N(0,0.20^2I_{512})$. We prescribe
$v_\tau(x)=g(x)+f(x)$ using \eqref{eq:app-manufactured-value} and
construct the twelve-mode optimal density below. Together they define
the supplied cost
\begin{equation}
 c_\infty(x,u)=\beta v_\tau(x)-\LL^u v_\tau(x)
 -\tau\log\frac{p^*(u\mid x)}{\varrho_\mu(u\mid x)}.
\label{eq:experiment-stationary-manufacture}
\end{equation}
The objective and stationary HJB equation are
\eqref{eq:stationary-routing-control} and
\eqref{eq:experiment-stationary-main}; there is no terminal boundary.

\paragraph{Manufactured law.}
Use the action coordinates $(\xi,\zeta)$ of
Section~\ref{subsec:experiment-stationary-routing}, with each block in
$\R^{256}$. Let $A_x\in\R^{1024\times8}$ and
$A_\xi,A_\zeta\in\R^{256\times12}$ have orthonormal columns and let
$M\in\R^{8\times12}$. The matrices are drawn once with construction
seed $26111101$. The three orthonormal matrices are obtained by QR
factorization of standard Gaussian matrices, and the entries of $M$
are Gaussian with variance $1/8$. Write $A_{\xi,8}$ for the first
eight columns of $A_\xi$, $e_j$ for coordinate vector $j$ in $\R^{12}$,
and set
\[
 \chi(x)=\tanh(2A_x^\top x),\qquad \kappa(x)=M^\top\chi(x).
\]
The mixture for $\xi$ uses covariance $0.17^2I_{256}$ in every
component. Its means and unnormalized log weights are
\[
 m_{\xi,j}(x)=0.08A_{\xi,8}\chi(x)
 +1.12[1+0.12\tanh\kappa_j(x)]A_\xi e_j,\qquad
 1.4\tanh\kappa_j(x)+b_j ,
\]
where the $b_j$ are evenly spaced from $-0.35$ to $0.35$. Given
$(x,\xi)$, the routing coordinate $\zeta$ is Gaussian with covariance
$0.17^2I_{256}$ and mean $A_\zeta a(x,\xi)$, with
\[
 r=A_\xi^\top \xi/0.17,\qquad
 a_j(x,\xi)=0.105\sin(0.8r_j+0.35\kappa_j(x))
             +0.065\tanh(r_jr_{j-1}),\qquad r_0=r_{12}.
\]
The orthogonal transformation between $(u^M,u^P)$ and
$(\xi,\zeta)$ has unit absolute Jacobian. The joint density is the
product of the mixture density of $\xi$ and the conditional density of
$\zeta$. This density and $v_\tau$ define the supplied cost via
\eqref{eq:experiment-stationary-manufacture}. Training uses the
supplied cost and dynamics, and the analytic law is used in evaluation.

\subsubsection{Training protocol}

\paragraph{Common initialization.}
The four solver--sampler combinations and projected CPVI start from
the same value checkpoint and initial mixture policy. Their critics have
the form
\[
 \hat v_\psi(x)=\tfrac12\langle q^2\rangle+
 \gamma_\psi^\top\phi(x)+0.1N_\psi(\widetilde x,\phi(x)),
\]
Here $\widetilde x=(q,\alpha/0.5,y/0.3,z/0.35)$, and $N_\psi$ has
two $192$-wide SiLU layers. The network has no physical-time input.
The public statistics supplied to the critic are
\begin{equation}
\begin{aligned}
 \phi(x)=\bigl(&1,\langle q^2\rangle,\langle\alpha^2\rangle,
 \langle y^2\rangle,\langle q\alpha\rangle,\langle qy\rangle,
 \langle\alpha y\rangle,\langle(\Pi_{\rm sec}q)^2\rangle,\\
 &\langle(\Pi_{\rm sec}\alpha)^2\rangle,
 \langle(\Pi_{\rm sec}y)^2\rangle,
 \langle(\Pi_{\rm sec}q)(\Pi_{\rm sec}\alpha)\rangle,
 \langle z^2\rangle,\langle qz\rangle\bigr).
\end{aligned}
\label{eq:app-public-statistics}
\end{equation}
Their common pools contain $768$ training and $96$ validation states
from reference-policy trajectories. To construct the initial policy
and subsequent Gibbs targets, use the dimensionless action energy
\begin{equation}
 \mathcal E_w(x,u)=\frac{c_\infty(x,u)+\LL^u w(x)}{\tau}
           -\log\varrho_\mu(u\mid x).
\label{eq:routing-action-energy}
\end{equation}
Thus the Gibbs law $q_w$ has density proportional to
$\exp[-\mathcal E_w(x,u)]$. Additive terms depending only on $x$
can be omitted when sampling or differentiating in $u$.

Initialization proceeds in three steps. First, at the mean training
state, run $384$ gradient-based restarts for $\mathcal E_g$, with
$g(x)=\tfrac12\langle q^2\rangle$, to identify twelve centers
$\bar u_j$; local curvature supplies a common isotropic variance.
Second, form $\pi_0$ from these centers and weights proportional to
$\exp[-\mathcal E_g(x,\bar u_j)]$. Third, fit the initial value using
eight relaxed trajectories per context under $\pi_0$, physical step
$0.25$, horizon $16$, and discounted quadratic tail $g$.
The reference in the KL penalty remains
$\mathcal N(0,0.20^2 I_{512})$ throughout.

\paragraph{Policy evaluation and flow improvement.}
Each round first freezes the deployed policy and fits its stationary
value, then freezes that value and updates the policy toward
$\hat q_k=q_{\hat v_k}$. Policy evaluation uses the stationary
residual $\LL^{\pi_k}\hat v_\psi+\ell^{\pi_k}-\beta\hat v_\psi$
and the deployed law's log density.

For DRDM, draw cached actions and antithetic neighboring states,
form two residual estimates, and minimize their product. Each of the
four evaluations uses $750$ critic updates, batches of $32$ states,
and $128$ cached actions per state. Each residual group draws eight
cache indices independently and uses physical difference step $0.005$.
The groups are independent conditional on the fixed cache.
For PINN, average drift and regularized running cost over the cached
actions, obtain the value gradient and Hessian--vector products by
automatic differentiation, and minimize the squared sampled residual.
The PINN runs share the initial checkpoint, state pools, critic
architecture, action cache size, and evaluation budget. They use four
Rademacher probes for the diffusion trace during training and sixteen
on validation states. Held-out public residuals select the critic
checkpoints.

The 12-mode flow actor uses an invertible affine-coupling map with a
change-of-variables density. For each improvement, first reweight the
base mixture by $\exp[-\mathcal E_{\hat v_k}(x,\bar u_j)]$, retaining
the discovered centers and variance. Then hold this base and
$\hat v_k$ fixed while training four affine coupling layers
conditioned on $x$ and $D_x\hat v_k(x)$. At each sampled training
state, draw $U$ from the current actor $\pi$ and minimize
\[
\E\!\left[\log p(U\mid x)+\mathcal E_{\hat v_k}(x,U)\right],
\qquad U\sim\pi(\cdot\mid x),
\]
where $p$ is that actor's change-of-variables density. This objective
differs from $\KL(\pi\Vert\hat q_k)$ only by a term fixed during
the actor fit. It trains the action map directly by reverse KL.
Each actor fit uses
$1500$ updates, and held-out free energies select its checkpoint.
The selected map and its density define the policy supplied to the
next evaluation and to deployment.

\paragraph{Projected CPVI.}
Projected CPVI updates the value and mixture moments simultaneously
from the old pair $(\hat v_k,\pi_k)$, following the coupled dynamics
of \citet{FW2025}. We use the cost-minimization convention: negate
their reward value and include $-\tau\log\varrho_\mu$ in the action
cost. This gives the same reference-relative entropy objective as
DRDM--flow.

First compute the action drift
$F_k(x,u)=-\tau\nabla_u\mathcal E_{\hat v_k}(x,u)$ from the frozen
old critic. Together with action diffusion $\sqrt{2\tau}\,I_{512}$,
it determines the moment update for every mixture component. At a
fixed context $x$, let $m_{k,j}$ and $v_{k,j}$ be component $j$'s
mean and diagonal-variance vector. With iteration step
$h_{\rm it}=0.25$, update
\begin{align*}
 m_{k+1,j}&=m_{k,j}+h_{\rm it}\E[F_k(x,U)],\\
 v_{k+1,j}&=v_{k,j}+2h_{\rm it}
 \bigl(\tau\mathbf 1+\E[(U-m_{k,j})\odot F_k(x,U)]\bigr).
\end{align*}
In these two expectations,
$U\sim\mathcal N(m_{k,j},\operatorname{diag}v_{k,j})$ with $x$ fixed,
and $\odot$ denotes coordinatewise multiplication. The labelled
components retain their initial weights throughout these updates.

Using the same old value and policy, form the value target
\[
\hat v_k+h_{\rm it}
\bigl(\LL^{\pi_k}\hat v_k+\ell^{\pi_k}-\beta\hat v_k\bigr).
\]
The term in parentheses is the stationary policy-evaluation residual;
its running cost uses the old mixture's normalized density $p_k$.
Fit the next critic to this target before forming the next pair of
updates. Each step uses six antithetic action draws per component,
four Rademacher probes for the state-diffusion trace, and twelve
minibatch updates for the value regression. The $320$ simultaneous
steps advance iteration time to $80$.

To deploy the resulting policy at new states, train a conditional
readout for $800$ updates to return the final means and diagonal
variances. The full mixture formed from these moments and the initial
weights supplies both action samples and the deployed density.
Across seeds the validation mean and
log-variance RMSEs range from $0.00206$ to $0.00231$ and $0.00309$ to
$0.00314$, respectively. An auxiliary particle calculation advances
$512$ particles at each of $12$
validation contexts under the same sequence of critic gradients
and compares average posterior modal masses at its endpoints.
The mean modal-mass change is at most $1.65\times10^{-4}$ in
total variation over the four seeds.

\paragraph{Stationary baseline implementations.}
The cluster baselines share the training and validation states for
each seed. Their value and policy updates differ as follows.

\emph{Projected Langevin.} Start both solver variants at the common
$(\hat v_0,\pi_0)$. At each improvement, freeze $\hat v_k$, start chains
from the Gaussian reference, and move actions along the score
$-\nabla_u\mathcal E_{\hat v_k}$ with Gaussian noise. Then fit a
conditional diagonal Gaussian mixture to the chain endpoints. This
fitted mixture is the policy whose density is used in the next
evaluation and in projected deployment. Both DRDM--projected
Langevin and PINN--projected Langevin use three improvements,
$128$ unadjusted score steps of size $0.004$, and $256$ endpoint
actions per context. Each DRDM policy evaluation uses $750$ updates.

\emph{Projected mirror descent.} Start at the same
$(\hat v_0,\pi_0)$, form an action density proportional to
$p_k^{0.75}\exp[-0.25\mathcal E_{\hat v_k}]$, and project its moments
to the same conditional mixture family. Here $p_k$ is the density
of $\pi_k$. Evaluate the projected policy before repeating the
improvement. The run uses ten improvements and $750$ DRDM updates
per evaluation.

\emph{Gaussian q-learning.} Initialize the value at the shared
checkpoint and the normalized diagonal-Gaussian policy at the
reference. Use reward $-c_\infty+\tau\log\varrho_\mu$ and parameterize
the little-q output as temperature times the policy log density.
Training applies the stationary discounted martingale derivative
test to on-policy transitions. It uses $20{,}000$ SGD updates and
four transitions of physical step $0.125$ per training trajectory.
The deployed policy remains in the same Gaussian family as in the
finite-horizon comparison.

\emph{Deep Splitting.} First regress the network
to the shared initial value on the training pool. Then perform discounted
fixed-point sweeps: freeze the next-state value and Hamiltonian driver,
form regression targets, and fit the new value before the next sweep.
It uses $80$ sweeps at physical step $0.125$, with $150$ regression
updates per sweep. The method computes the soft Hamiltonian by
importance sampling with the public mode-centered mixture.

For deployment, the final Deep Splitting value defines a
Gibbs target, which is projected to the conditional diagonal
Gaussian-mixture family used by the other projected baselines.
The Gibbs and mirror projections draw $64$ candidates per discovered
mode and fit the projected moments for $1000$ updates. All policy
metrics are computed from the resulting deployed normalized density.

\paragraph{Score--Hamiltonian PI.}
We use the \href{https://arxiv.org/html/2610.04285v1\#S3.SS2}{actor-first Gram updates}
and the \href{https://arxiv.org/html/2610.04285v1\#S5.SS3}{coordinatewise monotone actor}
of \citet{FW2026Amortized}. The policy is a monotone transformation
of independent Gaussian noise, with three tanh terms per action
coordinate; the critic is quadratic plus $96$ trainable tanh units.

Initialize the transport coefficients as constants and fit the
initial actor's value with $200$ numerical policy-evaluation steps.
The cited example instead initializes with the exact initial-policy
value. During each subsequent actor-first update, the current critic
and reference $\mu_x$ define the Gibbs law. In the cost convention,
the actor projects temperature times the Gibbs-minus-actor
log-density gradient onto its transport tangent space; the coupled
Gram update also updates the critic. Expectations use fresh joint
action draws so that the coupled action cost is retained.

The run takes $2000$ updates of size $0.1$, with actor and critic
damping $0.1$ and $0.01$, respectively. Each update uses $256$ states
from a $1024$-state discounted occupation pool, refreshed every
$50$ updates, and $16$ joint action draws per state. For each of four
seeds, Score--Hamiltonian PI runs on an RTX~3090 and is compared with
DRDM--flow on the same test states, rollout initial states, and
Brownian increments.

\subsubsection{Evaluation and reporting}
\label{app:routing-evaluation}

\paragraph{Policy and rollout errors.}
The comparison uses training seeds $26111150$, $26111160$,
$26111170$, and $26111180$, with a fresh pool of $128$ states per
seed and $512$ actions per test state for policy metrics.
To estimate modal masses, draw actions from the deployed policy and
evaluate the optimal mixture's posterior component probabilities at
these actions. Their averages estimate how much mass the deployed
policy assigns to each mode. Compare these estimates with the optimal
mixture weights at the same state to compute modal-mass total variation.
For a deployed stationary policy $\pi$, the reported discounted
closed-loop excess is
\begin{equation}
 \Delta_T(\pi)=\int_{\R^d}\E^\pi_{0,x}\!\left[\int_0^T e^{-\beta t}
 \tau\KL\bigl(\pi(\cdot\mid X_t)\Vert\pi^*(\cdot\mid X_t)\bigr)
 \dd t\right]\rho_{\rm init}(\dd x).
\label{eq:routing-closed-loop-excess}
\end{equation}
Here $\rho_{\rm init}$ is the benchmark's initial-state law and
$\E^\pi_{0,x}$ averages the relaxed state trajectory started at $x$.
To estimate this quantity, simulate $96$ paired relaxed state paths
through $T=16$ with shared Brownian increments. At each physical
step, use $32$ policy action draws to estimate the conditional KL
and the mean action in the relaxed drift, then accumulate the
discounted KL contribution along the path. Pair the fine Brownian
increments at step $0.125$ to obtain those at step $0.25$.
The interpretation follows from the manufactured cost and It\^o's
formula: \eqref{eq:routing-closed-loop-excess} is the running soft
cost through $T$ plus the discounted optimal tail, relative to the
mean initial optimal value.

\paragraph{Plotting protocol.}
Figure~\ref{fig:stationary-action-scatter} uses the first held-out
state for training seed $26111150$ and $1024$ action draws per panel.
All panels use the centered projection
\[
P^\top A_\xi^\top\bigl(\xi-0.08A_{\xi,8}\chi(x)\bigr),
\qquad
P_{j,:}=\sqrt{2/12}(\cos\theta_j,\sin\theta_j),
\quad \theta_j=2\pi(j-1)/12.
\]
Here $P$ has twelve rows, and the basis and state-dependent shift are
those of the manufactured law above. Contours and distances use the
projected optimal law $\pi^*$, with the contour levels and sliced
Wasserstein calculation in Appendix~\ref{app:experiment-evaluation}.

\paragraph{Timing protocol.}
Each cluster run uses one Tesla V100-SXM2-16GB GPU. Generation times
are medians of five synchronized batches, each drawing $128$
actions and evaluating their log densities at $16$ states after
one warm-up. CPVI training time includes the terminal readout fit
and the auxiliary particle calculation. Table~\ref{tab:stationary-multiseed-timing}
reports wall times under the stated method-specific budgets.
It excludes Score--Hamiltonian PI, which was run on an RTX~3090;
the main accuracy and cost comparison includes that method.

\begin{table}[!htbp]
\centering
\caption{Signed deviations in modal mass from the optimal policy,
in percentage points. At each test state, the estimated modal weights
average $\pi^*$'s posterior component probabilities over $512$ actions
drawn from the replayed policy. Each entry averages the estimated minus
exact weight over $128$ test states per seed and four training seeds.}
\label{tab:stationary-mode-bias}
\begingroup\small
\setlength{\tabcolsep}{2.5pt}
\begin{tabular}{@{}l*{12}{r}@{}}
\hline
Method & 1 & 2 & 3 & 4 & 5 & 6 & 7 & 8 & 9 & 10 & 11 & 12\\
\hline
DRDM--flow & $-0.5$ & $-5.7$ & $-4.7$ & $+1.3$ & $+5.5$ & $+0.7$ & $+9.3$ & $-6.0$ & $+0.5$ & $-0.9$ & $-8.6$ & $+9.1$\\
DRDM--projected Langevin & $+3.1$ & $+1.0$ & $+0.3$ & $+5.0$ & $-8.7$ & $+1.3$ & $+2.3$ & $+2.1$ & $+3.1$ & $-0.4$ & $-4.6$ & $-4.6$\\
PINN--flow & $-0.6$ & $-5.8$ & $-5.0$ & $+1.0$ & $+5.6$ & $+0.5$ & $+10.1$ & $-6.3$ & $+0.4$ & $-1.0$ & $-9.1$ & $+10.1$\\
PINN--projected Langevin & $+3.2$ & $+0.9$ & $+0.4$ & $+5.0$ & $-8.6$ & $+1.2$ & $+2.3$ & $+2.1$ & $+3.1$ & $-0.4$ & $-4.6$ & $-4.5$\\
Proj.\ CPVI & $-1.1$ & $-5.8$ & $-5.5$ & $+1.4$ & $+6.5$ & $+0.4$ & $+11.6$ & $-7.2$ & $+0.2$ & $-2.2$ & $-10.9$ & $+12.7$\\
Score--Hamiltonian PI & $+3.1$ & $+2.5$ & $+0.6$ & $+4.6$ & $-9.8$ & $+0.9$ & $+1.9$ & $+3.6$ & $+3.3$ & $-0.5$ & $-5.2$ & $-4.9$\\
Gaussian q-learning & $+3.4$ & $+2.3$ & $+1.0$ & $+4.7$ & $-9.0$ & $+1.0$ & $+1.5$ & $+2.8$ & $+3.4$ & $-0.7$ & $-5.2$ & $-5.1$\\
Deep Splitting & $-1.0$ & $-5.8$ & $-5.5$ & $+1.1$ & $+5.2$ & $+0.7$ & $+12.0$ & $-7.1$ & $+0.3$ & $-2.5$ & $-10.7$ & $+13.3$\\
Proj.\ mirror descent & $-0.9$ & $-5.8$ & $-5.4$ & $+1.3$ & $+6.1$ & $+0.6$ & $+11.0$ & $-6.9$ & $+0.3$ & $-1.8$ & $-10.4$ & $+11.8$\\
\hline
\end{tabular}
\endgroup
\end{table}

\begin{table}[htbp]
\centering
\caption{Cluster timings on the 12-mode stationary benchmark. Entries are means
and sample standard deviations over four seeds. Training excludes the
shared initialization ($30.04\pm0.24$\,s). Batch generation includes
the log-density evaluation for the deployed policy.}
\label{tab:stationary-multiseed-timing}
\vspace{.4em}
\begin{tabular}{lrr}
\hline
Method & Training (s) & $2048$ actions (ms)\\
\hline
DRDM--flow & $\mathbf{111.3}\pm10.0$ & $11.88\pm1.22$\\
DRDM--projected Langevin & $198.4\pm11.8$ & $5.46\pm0.19$\\
PINN--flow & $333.6\pm65.7$ & $11.02\pm1.36$\\
PINN--projected Langevin & $377.3\pm48.7$ & $5.15\pm0.62$\\
Projected CPVI & $223.1\pm12.8$ & $5.84\pm0.91$\\
Gaussian q-learning & $588.2\pm6.2$ & $\mathbf{0.71}\pm0.03$\\
Deep Splitting & $137.0\pm5.4$ & $5.29\pm0.14$\\
Projected mirror descent & $130.6\pm6.2$ & $5.48\pm0.40$\\
\hline
\end{tabular}
\end{table}

\paragraph{Finite-chain Langevin diagnostic.}
This diagnostic deploys the actual finite chain from the last
Langevin improvement, whose law differs from the fitted mixture used
in the projected comparison. To evaluate it without an exact chain
density, bound its entropy below by the entropy of the last independent
Gaussian increment and above by cross-entropy against the fitted
mixture. Insert these entropy bounds into the conditional KL in
\eqref{eq:routing-closed-loop-excess} and estimate the resulting
closed-loop bounds. Their mean estimates across seeds are $0.00178$
and $0.2220$; this range includes the Flow excess estimate. Timing
the actual finite chain gives $902\pm51$\,ms for $2048$ actions over
the four seeds.
\FloatBarrier

\subsection{Experiment 4: Unregularized HJB equation}
\label{app:hjb3a-solver-pairs}

\subsubsection{Problem setting}
This experiment uses the unregularized HJB-3a problem specified in
\eqref{eq:experiment-hjb3a-problem}, with $T=1$,
$d=m\in\{1000,10{,}000\}$, and initial state uniform on $[-1,1]^d$
and independent of the Brownian motion. The separable structure permits
one scalar value network to be shared across coordinates. In the shifted
coordinate $y=x_i-t$, write this network as $\hat v_{\rm sc}(t,y)$ and
let $a(t,y)$ be the conditional mean action. The full numerical value is
\[
\hat v(t,x)=\frac1d\sum_{i=1}^d\hat v_{\rm sc}(t,x_i-t).
\]
Both PDE solvers impose $\hat v_{\rm sc}(1,y)=\phi(y)$ by construction.

\subsubsection{Training protocol}
\paragraph{Evaluate the current feedback.}
Start with mean action zero. During each value fit, hold $a$ fixed and
solve the scalar policy-evaluation equation
\[
\partial_t\hat v_{\rm sc}
+2a\,\partial_y\hat v_{\rm sc}
+0.04\,\partial_{yy}\hat v_{\rm sc}+25a^2=0.
\]
PINN evaluates these derivatives by automatic differentiation and
minimizes the squared residual at sampled time--state pairs. DRDM
approximates the differential terms with antithetic Gaussian
successors and minimizes the product of two independent residual
estimates, using the construction in Appendix~\ref{app:solver-implementations}.

\paragraph{Improve the feedback and repeat.}
Freeze the fitted critic. Minimizing its quadratic action objective
gives the center
\[
a_c(t,y)=-0.04\,\partial_y\hat v_{\rm sc}(t,y).
\]
Flow fits a conditional Gaussian translation at each improvement.
Langevin uses the score of the Gaussian centered at $a_c$, with action
standard deviation $0.05$. Each chain starts from
$\mathcal N(0,0.05^2)$ and takes $16$ steps of size
$h_{\rm L}=0.05^2/4$. Holding the condition $(t,y)$ fixed, its mean obeys
\[
m_{j+1}=\tfrac34m_j+\tfrac14a_c,\qquad m_0=0,
\qquad m_{16}=[1-(3/4)^{16}]a_c.
\]
DRDM--Langevin and PINN--Langevin use this finite-chain mean in the next
policy evaluation and in closed-loop deployment. Thus the next value
fit uses the feedback returned by the sampler. We verified the formula
for $m_{16}$ against antithetic chain endpoints after each improvement.

\paragraph{Training budget.}
Each DRDM or PINN run makes four policy evaluations with $1600$ critic
updates per evaluation, interleaved with four policy improvements.
All six methods use $6400$ value-network updates per seed at each dimension.
The reported value is the last fitted scalar network, and the reported
controller is obtained by improving that network. There is no additional
fifth policy evaluation, so this value is not presented as a separate
evaluation of the final mean controller. The Gaussian target is used
only for the policy update; the unregularized cost contains neither an
entropy penalty nor an action-variance contribution.

\subsubsection{Evaluation and reporting}
\label{app:hjb3a-evaluation}
After development at $d=1000$, fix the settings and train three fresh
seeds at each dimension. Reconstruct the full value from the shared
scalar network and compute $E_V$ on $48$ common time--state pairs.
The reference is \eqref{eq:experiment-hjb3a-reference}, whose scalar
expectation is evaluated by $128$-point Gauss--Hermite quadrature.
Compute $J_0$ on $128$ paired trajectories with $384$ Euler steps,
using the deployed mean action in both the state dynamics and running
cost, including the finite-chain correction above. For the two flow
controllers, form a $95\%$ paired $t$-interval from the three seed-wise
differences in mean cost. Section~\ref{subsec:experiment-hjb3a} reports
the value plots and the means and sample standard deviations of the
metrics over the three training seeds.

\paragraph{Plotting protocol.}
Figure~\ref{fig:hjb3a-values} follows the spatial-slice presentation of
\citet[Figure~5]{CFZZ2024DFMartNet}: fix $t=0$ and set
$x=s\mathbf1_d$ for $s\in[-1,1]$. Evaluate the frozen final critics
and the Cole--Hopf reference on $201$ equally spaced values of $s$.
The blue curve shows the reference; red markers show the mean over
three training seeds at every eighth grid point. The normalized
additive structure makes the reference profile independent of $d$.
Figure~\ref{fig:hjb3a-value-errors} uses the original $48$ held-out
time--state pairs, grouped into eight states at each of
$t=0,0.25,0.5,0.75,0.9,0.97$. For each seed and time group, compute
the relative $L^2$ value error, then average these errors over seeds.
The vertical axis is logarithmic, and errors are reported in percent.

\FloatBarrier

\section{Technical estimates for policy iteration}\label{app:technical-proofs}
The transition smoothing, Gibbs normalizer, and Abel-kernel estimates
support the error propagation proof in Section~\ref{subsec:main-proofs}.

\subsection{Uniform transition smoothing}
\label{app:uniform-smoothing}
\begin{proof}[Proof of Lemma~\ref{lem:uniform-smoothing}]
We use the transition-density gradient estimate of
\citet[Theorem~1.2(ii)]{MPZ2021}, which allows coefficients that are
measurable in time. Every averaged drift satisfies
\[
|\bar b^\pi(t,x)|\le\int_{\mathcal U}|b(t,x,u)|\pi(\dd u\mid t,x)\le B.
\]
In particular, $|\bar b^\pi(t,0)|\le B$ and
$|\bar b^\pi(t,x)-\bar b^\pi(t,y)|\le2B$, so their drift condition
$(H^b_0)$ holds with a common constant $\kappa_1=2B$.
Their diffusion condition holds for $A^{1/2}$, which gives the same
generator as $\sigma$. Boundedness and uniform spatial Lipschitz
continuity of $\sigma$ give the corresponding properties for
$A=\sigma\sigma^\top$. Uniform ellipticity makes the matrix
square-root map Lipschitz on the range of $A$, so $A^{1/2}$ is
bounded and uniformly Lipschitz in space. Its bounds depend only on
$a_0$, the uniform bound on $\sigma$, and the spatial Lipschitz constant
of $\sigma$. These bounds also give a uniform spatial H\"older bound
of exponent $\alpha=1/2$.
Consequently all parameters in their Theorem~1.2(ii) are uniform in
the policy; no derivative bound on $\bar b^\pi$ enters that estimate.
Writing $p^\pi(t,x;s,y)$ for the Lebesgue density of
$\mathbb P^\pi_{t,x}(X_s\in\cdot)$ on $\R^d$ and
$\vartheta^\pi_{s,t}(x)$ for the regularized drift flow in that theorem,
its estimate reads
\[
|D_xp^\pi(t,x;s,y)|\le
\frac{C_1}{\sqrt{s-t}}(s-t)^{-d/2}
\exp\!\left(-\lambda_1
\frac{|y-\vartheta^\pi_{s,t}(x)|^2}{s-t}\right),
\]
with common $C_1,\lambda_1>0$. The Gaussian integral in $y$ is
$(\pi/\lambda_1)^{d/2}$ independently of its policy-dependent center.
Hence, taking $K_T=C_1(\pi/\lambda_1)^{d/2}$ gives
\[
\sup_x\int_{\R^d}|D_xp^\pi(t,x;s,y)|\dd y
\le K_T(s-t)^{-1/2}.
\]
Differentiating
$P^\pi_{t,s}f(x)=\int f(y)p^\pi(t,x;s,y)\dd y$ proves
\eqref{eq:uniform-transition-gradient}. Periodic coefficients and
functions can be lifted to $\R^d$, giving the same conclusion on a torus.
\end{proof}

\subsection{Gibbs normalizing constants}
\label{app:gibbs-normalization}
At a fixed $(t,x)$, use the reduced energies in \eqref{eq:action-energy},
with the common diffusion factor already canceled, and write
\[
\widehat Z_k=\int_{\mathcal U}e^{-E_{\hat{v}_k}(t,x,u)/\tau}\,\mu_{t,x}(\dd u),
\qquad Z_*=\int_{\mathcal U}e^{-E_{v_\tau}(t,x,u)/\tau}\,\mu_{t,x}(\dd u).
\]
Exponentiating \eqref{eq:proof-energy-perturbation} and integrating
against the same reference measure gives
\begin{equation}\label{eq:proof-partition-ratio}
\begin{gathered}
\exp[-B(d_k(t)+\gamma_k(t))/\tau]Z_*\le\widehat Z_k
\le \exp[B(d_k(t)+\gamma_k(t))/\tau]Z_*
\\
\Longrightarrow\quad
\left|\log\frac{\widehat Z_k}{Z_*}\right|
\le\frac{B(d_k(t)+\gamma_k(t))}\tau.
\end{gathered}
\end{equation}
The Gibbs densities with respect to $\mu_{t,x}$ are strictly positive.
Their quotient gives the following version of the log Radon--Nikodym derivative,
\[
\log\frac{\dd\hat{q}_k}{\dd\pi^*}(u\mid t,x)
=-\frac{E_{\hat{v}_k}(t,x,u)-E_{v_\tau}(t,x,u)}\tau
-\log\frac{\widehat Z_k(t,x)}{Z_*(t,x)}.
\]
Each term on the right is bounded in absolute value by
$B(d_k(t)+\gamma_k(t))/\tau$, proving
\eqref{eq:gibbs-log-ratio-gradient}.

\subsection{Abel-kernel convolution}
\label{app:abel-convolution}
For $K_{\mathrm{Abel}}(u)=u^{-1/2}$, Euler's beta integral
\citep[Eq.~(5.12.1)]{DLMFBeta2026} gives
\[
\int_0^h(h-r)^{a-1}r^{b-1}\dd r
=h^{a+b-1}\frac{\Gamma(a)\Gamma(b)}{\Gamma(a+b)},
\qquad a,b>0.
\]
Induction on the number of convolutions therefore yields
\begin{equation}\label{eq:abel-convolution}
K_{\mathrm{Abel}}^{*m}(h)=\frac{\pi^{m/2}}{\Gamma(m/2)}h^{m/2-1},
\qquad m\ge1.
\end{equation}
This is the half-order case of the Riemann--Liouville composition law
$I^\alpha I^\beta=I^{\alpha+\beta}$
\citep[Eqs.~(1.15.47)--(1.15.48)]{DLMFFractional2026}.
The same iterated-kernel argument underlies the weakly singular
Gronwall estimate of \citet[Lemma~7.1.1]{Henry1981}.
Here we use the explicit kernel powers in the iteration index;
the historical-error limit follows from their summability and the
finite-head/tail argument in the main proof.

\begingroup
\linespread{0.95}
\small
\urlstyle{same}
\setlength{\bibsep}{0pt}
\bibliographystyle{plainnat}
\bibliography{references}
\endgroup

\end{document}